\documentclass[11pt,a4paper,leqno]{article}
\usepackage[T1]{fontenc}
\usepackage[utf8]{inputenc}
\usepackage{lmodern}
\usepackage{microtype}
\usepackage{geometry}
\usepackage{amsmath,amssymb,amsthm,mathtools,mathrsfs,bm}
\usepackage{graphicx,float}
\usepackage{caption}
\usepackage{enumitem}
\usepackage{needspace}
\usepackage{hyperref}
\hypersetup{colorlinks=true,citecolor=black,linkcolor=black,urlcolor=black}
\numberwithin{equation}{section}
\newtheorem{theorem}{Theorem}[section]
\newtheorem{proposition}[theorem]{Proposition}
\newtheorem{corollary}[theorem]{Corollary}
\newtheorem{lemma}[theorem]{Lemma}
\newtheorem{property}{Property}[section]
\theoremstyle{definition}
\newtheorem{definition}{Definition}[section]
\newtheorem{numbereddefinition}[definition]{Definition}
\theoremstyle{remark}
\newtheorem{remark}{Remark}[section]

\allowdisplaybreaks[1]
\AtBeginDocument{%
  \setlength{\abovedisplayskip}{9pt plus 2pt minus 4pt}%
  \setlength{\belowdisplayskip}{9pt plus 2pt minus 4pt}%
}
\newcommand{\D}{\mathcal D}
\newcommand{\Ein}{\mathcal E^{1,2}}
\newcommand{\AWitt}{\mathbb A}

\newcommand{\A}{\mathcal A}
\newcommand{\Q}{\mathbb Q}

\newcommand{\Eell}{\mathrm E}
\newcommand{\R}{\mathbb R}

\begin{document}

{\centering\fontsize{13}{15.5}\selectfont\bfseries Integrability and closed spacelike extremals\\[3pt]of conformal arclength in the (2+1)-dimensional Einstein universe\par}
\vspace{1.15em}
{\centering\normalsize\scshape Emilio Musso\par}
\vspace{1.1em}

\begin{quote}
\small\noindent\textbf{Abstract.}
We study the conformal arclength functional for nondegenerate spacelike curves in the three-dimensional Einstein universe and classify its closed extremals up to conformal equivalence. We show that the extremals are projections of trajectories of a generically noncommutatively integrable Reeb flow. For extremals with nonconstant conformal invariants, closedness forces negative helicity and regular elliptic momentum. We prove that the restricted moment map over the corresponding invariant region defines a nontrivial principal $\mathbb T^3$-bundle, and construct explicit local trivializations that linearize the Reeb flow. This yields explicit parametrizations in terms of elliptic functions and integrals. We prove that the associated period map is a real-analytic diffeomorphism onto a bounded, convex planar domain whose rational points parametrize the conformal equivalence classes of closed extremals with nonconstant conformal invariants. We conclude with a study of the global geometry of closed extremals.
\end{quote}
\vspace{0.35em}

\begingroup
\renewcommand{\thefootnote}{}
\footnotetext{\par
Author partially supported by PRIN project ``Real and Complex Manifolds:
Topology, Geometry and Holomorphic Dynamics'' (protocols
2017JZ2SW5-004 and 2022AP8HZ9-003); by the GNSAGA of INdAM;
and by Politecnico di Torino through the ``Finanziamento Individuale per le Missioni
Istituzionali di Ateneo''.}
\endgroup

\section*{Introduction}

In homogeneous geometries, the moving frame method (see~\cite{Cartan1935,Chern1985,FelsOlverI,FelsOlverII}) provides a systematic framework for constructing canonical frames and differential invariants. For curves of constant type (see~\cite{Green1978,Jensen1977}), this construction can single out an invariant arclength element and a set of fundamental differential invariants. The arclength element defines a natural length functional and, together with the differential invariants, provides the ingredients for more general invariant variational functionals.

In conformal geometry, this construction yields the conformal arclength element and a set of fundamental differential invariants. The conformal differential geometry of curves has been developed in classical works by Takasu, Vessiot, and Fialkow~\cite{Takasu,Vessiot1925,Fialkow}, as well as in more recent studies employing invariant moving frames~\cite{MariBeffa} and tractor calculus~\cite{SilhanZadnik}.

In the present paper, we study the conformal arclength functional for nondegenerate spacelike curves in the three-dimensional Einstein universe, the homogeneous model of Lorentzian conformal geometry in dimension three. Our main goal is to classify its closed extremals up to conformal equivalence and investigate their geometry. The conformal arclength problem has been studied in M\"obius geometry~\cite{Musso1994,MagliaroMariRigoli2011,MussoNicolodi2017} and, for timelike curves, in Lorentzian conformal geometry~\cite{DzhalilovMussoNicolodi2016,EshkobilovMussoNicolodi2021}. Although the timelike and spacelike problems share several formal features, their differential geometries differ in an essential respect: generic timelike curves belong to a single geometric type, whereas generic spacelike curves fall into two types, distinguished by their osculating flags and encoded by the sign of their helicity. These two types give rise to distinct canonical Frenet systems and must both be considered in the variational analysis.

Closed spacelike extremals with constant conformal invariants are homogeneous: their images are simple closed curves on anti-de Sitter walls, and their conformal equivalence classes are parametrized by the rational numbers in $(0,1)$. The classification of extremals with nonconstant conformal invariants is more involved. Even when their conformal invariants are periodic, the corresponding curves need not be closed: closedness is equivalent to the finite-order condition on the monodromy of their canonical frames. In~\cite{DzhalilovMussoNicolodi2016}, a period map characterizes closedness for timelike extremals in a prescribed spectral region. For the spacelike problem, we first determine the spectral restrictions imposed by closedness itself and then study the global dependence of the monodromy on the parameters.

We formulate the variational problem in Griffiths' formalism~\cite{Griffiths}. Using Hsu's characterization of singular integral curves in terms of nonzero characteristic lifts~\cite{Hsu}, we prove the regularity of the admissible integral curves of the canonical Pfaffian system. We then apply the variational criteria of Bryant and Hsu~\cite{Bryant,Hsu} to identify the extremals with the projections of trajectories of a Reeb flow on the momentum space. We compute the associated contact moment map and prove weak noncommutative integrability for both helicities on the open dense regular locus of the moment map, and noncommutative integrability on the restricted momentum spaces of negative helicity, where the restricted moment maps are proper submersions with compact $\mathbb T^3$-fibres (cf.~\cite{Jovanovic2012,Jov2025}). We then use momentum reduction (cf.~\cite{GutierrezSagredo2025,MussoGrant2004}) to study the evolution of the conformal curvature and torsion, identifying the reduced spaces with unions of connected components of real algebraic phase curves and determining the induced Reeb dynamics. Consequently, when the conformal invariants are nonconstant and periodic, they can be expressed in terms of Jacobi elliptic functions.

For spacelike extremals with nonconstant conformal invariants, we prove that closedness forces negative helicity and regular elliptic momentum. The proof uses the frame monodromy to rule out positive helicity and, in the negative-helicity case, to exclude all remaining spectral possibilities, including the degenerate cases. Combining these necessary conditions with the classification of periodic reduced trajectories, we reduce the classification problem to the study of periodic trajectories of the negative-helicity Reeb flow in a distinguished invariant region of the momentum space.

We prove that, on this region, the restricted moment map defines a nontrivial principal $\mathbb T^3$-bundle (Theorem~\ref{thm:invariant-torus-fibration}). Its nontriviality obstructs a global choice of torus coordinates over the image of the restricted moment map. We construct explicit local trivializations that linearize the Reeb vector field (Theorem~\ref{thm:linearization-reeb}). In these coordinates, its restriction to each invariant torus generates a linear flow with angular frequencies depending only on the Casimir invariants. The construction also yields explicit parametrizations of the corresponding extremals, whether closed or not, in terms of Jacobi elliptic functions and elliptic integrals of the third kind.

The normalized angular frequencies of the linearized Reeb flow determine a period map. For extremals arising from the distinguished invariant region, closedness is equivalent to the rationality of their period coordinates. As in the M\"obius classification of~\cite{MussoNicolodi2017}, we prove that the period map is a real-analytic diffeomorphism onto a bounded, convex planar domain $\widetilde{\mathcal D}$ (Theorem~\ref{thm:period-map}). This yields a bijection between the rational points of $\widetilde{\mathcal D}$ and the conformal equivalence classes of all closed spacelike extremals with nonconstant conformal invariants.

We also investigate the global geometry of closed extremals and prove that none is contained in any chamber of the Einstein universe. For the nonhomogeneous closed extremals, we determine the full conformal symmetry groups and study their spherical projections, including rotational symmetries and geodesic inflections. We interpret the characteristic numbers associated with the rational period coordinates in terms of winding numbers and symmetry orders.

The paper is organized as follows. Section~1 develops the conformal geometry of nondegenerate spacelike curves in the Einstein universe and establishes the regularity of the admissible integral curves of the canonical Pfaffian system. Section~2 characterizes extremals through the Reeb flow, studies the moment map and integrability, analyzes the reduced dynamics, and classifies the homogeneous closed extremals. Section~3 establishes the spectral restrictions, constructs the invariant torus fibration, linearizes the Reeb flow on its fibres, derives explicit parametrizations of the corresponding extremals, and proves the global classification through the period map. Section~4 develops the geometric results and presents numerical examples.

For Jacobi elliptic functions and elliptic integrals, we adopt the notation of Abramowitz and Stegun~\cite{AS} and use the NIST Digital Library of Mathematical Functions (DLMF)~\cite{dlmf} as our main reference for identities and analytic properties. The Einstein summation convention is understood throughout.

\section{Geometric preliminaries}

\subsection{The Einstein universe and its conformal geometry}\label{subsec:einstein-universe}

We briefly recall the main features of the conformal Lorentzian geometry of the $(2+1)$-dimensional Einstein universe. For a comprehensive treatment, see~\cite{Barbot2008}.

Let $\Ein$ denote the three-dimensional submanifold of $\mathbb R^5$ defined by $x_1^2+x_2^2=1$ and $x_3^2+x_4^2+x_5^2=1$. The restriction to $\Ein$ of the quadratic form $-dx_1^2-dx_2^2+dx_3^2+dx_4^2+dx_5^2$ induces a Lorentzian metric. We endow $\Ein$ with the product orientation and the time-orientation determined by the unit timelike vector field $T=-x_2\partial_{x_1}+x_1\partial_{x_2}$, which we take to be future-directed. The resulting oriented and time-oriented Lorentzian manifold is called the \emph{$(2+1)$-dimensional Einstein universe}.

To describe its conformal geometry, we introduce the coordinates $Y=(y_1,\ldots,y_5)^t$ defined by
\[
y_1=(x_1+x_5)/\sqrt2,\quad y_2=x_3,\quad y_3=x_4,\quad y_4=x_2,\quad y_5=(x_5-x_1)/\sqrt2.
\]
The corresponding basis $(E_1,\ldots,E_5)$ of $\mathbb R^5$ is called the \emph{standard Witt basis}. We equip $\mathbb R^5$ with the scalar product of signature $(2,3)$ associated with the quadratic form introduced above, which, in the standard Witt basis, is given by
\[
\langle Y,Z\rangle=y_1z_5+y_5z_1+y_2z_2+y_3z_3-y_4z_4=g_{ij}y_i z_j,
\]
where $\mathbf G=(g_{ij})$ is the corresponding Gram matrix.

We equip $\mathbb R^5$ with the orientation defined by $dV=dy_1\wedge\cdots\wedge dy_5$ and denote by $\mathbb R^{2,3}$ the resulting oriented pseudo-Euclidean space. We use capital symbols for vectors in this space and lowercase Latin letters for their coordinates. Unless otherwise specified, vectors will be regarded as column vectors with respect to the standard Witt basis.

For each nonzero vector $Y\in\mathbb R^{2,3}$, let $[Y]=\{\lambda Y:\lambda>0\}$ denote the ray generated by $Y$. The map $\Ein\ni Y\mapsto[Y]$ identifies the Einstein universe with the manifold of null rays. Henceforth, we shall make no distinction between these two realizations. The null-ray model is a twofold covering of the projective model considered in~\cite{Barbot2008}.

\Needspace{12\baselineskip}
For any nonzero vector $V\in\mathbb R^{2,3}$, set
\[
\mathcal T_\pm(V)=\{[Y]\in\Ein:\pm\langle Y,V\rangle>0\},\qquad
\mathcal T_\infty(V)=\{[Y]\in\Ein:\langle Y,V\rangle=0\}.
\]
Then $\Ein=\mathcal T_+(V)\sqcup\mathcal T_\infty(V)\sqcup\mathcal T_-(V)$. Depending on whether $V$ is spacelike, null, or timelike, the open sets $\mathcal T_\pm(V)$ are conformally equivalent to anti-de Sitter, Minkowski, or de Sitter three-space, respectively. In the anti-de Sitter case, $\mathcal T_\infty(V)$ is a totally umbilical timelike torus; in the Minkowski case, a null cone with two singular vertices; and in the de Sitter case, the disjoint union of two totally umbilical spacelike spheres. For further details, we refer the reader to~\cite{DzhalilovMussoNicolodi2016}.

\begin{definition}\label{def:chambers-walls}
The open sets $\mathcal T_\pm(V)$ are called \emph{anti-de Sitter chambers}, \emph{Minkowski chambers}, or \emph{de Sitter chambers}, according as $V$ is spacelike, null, or timelike. If $V$ is spacelike or null, the set $\mathcal T_\infty(V)$ is called the corresponding \emph{wall}. If $V$ is timelike, the two connected components of $\mathcal T_\infty(V)$ are called the corresponding \emph{de Sitter walls}.
\end{definition}

The pseudo-orthogonal group associated with the standard Witt basis is
\[
\AWitt(2,3)=\left\{\mathbf A\in GL(5,\mathbb R):\mathbf A^t\mathbf G\mathbf A=\mathbf G\right\}.
\]
We denote by $\AWitt_+(2,3)$ and $\AWitt_-(2,3)$ the subsets of $\AWitt(2,3)$ consisting of matrices with determinant $+1$ and $-1$, respectively. Both $\AWitt(2,3)$ and $\AWitt_+(2,3)$ act effectively and transitively on the Einstein universe from the left by $\mathbf A\cdot[Y]=[\mathbf A Y]$. These actions preserve its conformal Lorentzian structure, and the action of $\AWitt_+(2,3)$ also preserves its orientation.

Each of $\AWitt_+(2,3)$ and $\AWitt_-(2,3)$ has two connected components. To distinguish them, consider the cone $\mathcal C=\{U\wedge V\ne0:\langle U,U\rangle=\langle V,V\rangle=\langle U,V\rangle=0\}\subset\Lambda^2\mathbb R^{2,3}$ consisting of nonzero totally isotropic bivectors. Let $\Phi:\Lambda^2\mathbb R^{2,3}\to\mathbb R$ be defined by $\Phi(U\wedge V)=dV(U,V,E_2,E_3,E_1+E_5)$. The functional $\Phi$ does not vanish on $\mathcal C$, and the sets $\mathcal C^\pm=\{\beta\in\mathcal C:\pm\Phi(\beta)>0\}$ are its two connected components. This choice of $\mathcal C^+$ agrees with the time-orientation of $\Ein$ fixed above.

Every element of $\AWitt(2,3)$ either preserves these two components or interchanges them. Accordingly, we denote by $\AWitt^\uparrow(2,3)$ the subgroup preserving $\mathcal C^+$, and by $\AWitt^\downarrow(2,3)$ the subset interchanging $\mathcal C^+$ and $\mathcal C^-$. The four connected components of $\AWitt(2,3)$ are therefore
\[
\AWitt_\pm^\uparrow(2,3)=\AWitt_\pm(2,3)\cap\AWitt^\uparrow(2,3),\qquad \AWitt_\pm^\downarrow(2,3)=\AWitt_\pm(2,3)\cap\AWitt^\downarrow(2,3).
\]
In particular, $\AWitt_+^\uparrow(2,3)$ is the identity component of $\AWitt(2,3)$, consisting of the transformations preserving both the orientation and the time-orientation of the Einstein universe.

The map $\mathbf A\mapsto(A_1,\ldots,A_5)$, where $A_1,\ldots,A_5$ are the columns of $\mathbf A$, identifies $\AWitt_+(2,3)$ with the manifold of positively oriented Witt bases of $\mathbb R^{2,3}$, that is, positively oriented bases satisfying $\langle A_i,A_j\rangle=g_{ij}$ ($1\le i,j\le5$). Such a basis is said to be \emph{time-oriented} if $A_1\wedge(A_2+A_4)\in\mathcal C^+$. Equivalently, $\mathbf A\in\AWitt_+^\uparrow(2,3)$.

The Lie algebra of $\AWitt(2,3)$ is
\[
\mathfrak a(2,3)=\{\mathbf X\in\mathfrak{gl}(5,\mathbb R):\mathbf X^t\mathbf G+\mathbf G\mathbf X=0\}.
\]
Let $\boldsymbol{\Omega}=(\omega^i{}_j)$ denote the left-invariant, $\mathfrak a(2,3)$-valued Maurer--Cartan form on $\AWitt_+(2,3)$. It is characterized by
\begin{equation}\label{eq:witt-frame-equations}
dA_j=\omega^i{}_j A_i,\qquad 1\le j\le5.
\end{equation}
The Maurer--Cartan form has the matrix representation
\[
\boldsymbol{\Omega}=
\begin{pmatrix}
\omega^1{}_1&-\omega^2{}_5&-\omega^3{}_5&\omega^4{}_5&0\\
\omega^2{}_1&0&-\omega^3{}_2&\omega^4{}_2&\omega^2{}_5\\
\omega^3{}_1&\omega^3{}_2&0&\omega^4{}_3&\omega^3{}_5\\
\omega^4{}_1&\omega^4{}_2&\omega^4{}_3&0&\omega^4{}_5\\
0&-\omega^2{}_1&-\omega^3{}_1&\omega^4{}_1&-\omega^1{}_1
\end{pmatrix}.
\]
Differentiating~\eqref{eq:witt-frame-equations} yields the Maurer--Cartan equations
\[
d\omega^i{}_j=-\omega^i{}_k\wedge\omega^k{}_j,\qquad 1\le i,j\le5.
\]

In the following sections, we shall use both Witt and pseudo-orthonormal bases of $\mathbb R^{2,3}$. We conclude by introducing the latter and comparing the corresponding matrix representations of the pseudo-orthogonal group.

An ordered basis $(F_1,\ldots,F_5)$ of $\mathbb R^{2,3}$ is called \emph{pseudo-orthonormal} if $\langle F_i,F_j\rangle=h_{ij}$, where $\mathbf G_0=(h_{ij})=\operatorname{diag}(-1,-1,1,1,1)$. Such a basis is said to be \emph{positively oriented} if $dV(F_1,\ldots,F_5)=1$, and \emph{time-oriented} if $\Phi(F_1\wedge F_2)>0$. In particular, the basis
\[
B_1=(E_1-E_5)/\sqrt2,\quad
B_2=E_4,\quad B_3=E_2,\quad B_4=E_3,\quad
B_5=(E_1+E_5)/\sqrt2
\]
is positively oriented and time-oriented. We call it the \emph{standard pseudo-orthonormal basis} of $\mathbb R^{2,3}$.

Let $\mathsf E^i_j$ denote the $5\times5$ elementary matrix with entry $1$ in position $(i,j)$ and all other entries zero. The change-of-basis matrix from the standard pseudo-orthonormal basis to the standard Witt basis is
\begin{equation}\label{eq:change-of-basis-matrix}
\mathbf P=(\mathsf E^1_1-\mathsf E^5_1+\mathsf E^1_5+\mathsf E^5_5)/\sqrt2+\mathsf E^4_2+\mathsf E^2_3+\mathsf E^3_4.
\end{equation}
Then $\mathbf P\in SL(5,\mathbb R)$ and $\mathbf P^t\mathbf G\mathbf P=\mathbf G_0$. The pseudo-orthogonal group in the standard pseudo-orthonormal basis is
\[
O(2,3)=\{\mathbf M\in GL(5,\mathbb R):\mathbf M^t\mathbf G_0\mathbf M=\mathbf G_0\}.
\]
We also set $SO(2,3)=O(2,3)\cap SL(5,\mathbb R)$. The two matrix representations are related by conjugation: $\AWitt(2,3)=\mathbf P\,O(2,3)\,\mathbf P^{-1}$ and $\AWitt_+(2,3)=\mathbf P\,SO(2,3)\,\mathbf P^{-1}$. Accordingly, the four connected components of $O(2,3)$ are $O_\pm^\uparrow(2,3)=\mathbf P^{-1}\AWitt_\pm^\uparrow(2,3)\mathbf P$ and $O_\pm^\downarrow(2,3)=\mathbf P^{-1}\AWitt_\pm^\downarrow(2,3)\mathbf P$.

\subsection{Conformal geometry of spacelike curves}

Let $\gamma:I\to\Ein$ be an oriented spacelike curve, and let $\Gamma:I\to\mathbb R^{2,3}$ be an arbitrary null lift of $\gamma$. Since $\gamma$ is spacelike, its second osculating space $\mathcal O_2(\gamma)_t=\operatorname{span}\{\Gamma(t),\Gamma'(t),\Gamma''(t)\}$ is a nondegenerate three-dimensional subspace of signature $(1,2)$, independent of the choice of $\Gamma$. Consequently, its orthogonal complement $\mathcal N_\gamma(t)=\mathcal O_2(\gamma)_t^\perp$ is a Lorentzian two-plane. Let $\pi_{\mathcal N}(t):\mathbb R^{2,3}\to\mathcal N_\gamma(t)$ denote the orthogonal projection. The quartic differential
\[
Q_\gamma
=\frac{\langle\pi_{\mathcal N}(\Gamma'''),\pi_{\mathcal N}(\Gamma''')\rangle}
{\langle\Gamma',\Gamma'\rangle}\,dt^4
\]
is independent of the choice of the null lift and invariant under conformal transformations. The sign of $Q_\gamma$ determines a conformally invariant classification of the points of $\gamma$. A point is called \emph{elliptic}, \emph{hyperbolic}, or \emph{parabolic} according as $Q_\gamma$ is positive, negative, or zero.

\begin{definition}\label{def:nondegenerate-helicity}
A spacelike curve is called \emph{generic}, or \emph{nondegenerate}, if it has no parabolic points. The \emph{helicity} of a generic spacelike curve $\gamma$ is defined by $\varepsilon=\operatorname{sgn}(Q_\gamma)\in\{-1,1\}$. For an oriented generic spacelike curve, a parameter $s$, compatible with its orientation, is called a \emph{natural parameter} if $Q_\gamma=\varepsilon\,ds^4$. Since $Q_\gamma$ never vanishes, such a parameter exists and is unique up to translations. We also set $a_\varepsilon=(1+\varepsilon)/2$ and $b_\varepsilon=(1-\varepsilon)/2$, so that $(a_1,b_1)=(1,0)$ and $(a_{-1},b_{-1})=(0,1)$. Unless otherwise stated, all spacelike curves will be assumed to be generic and parametrized by their natural parameter.
\end{definition}

\begin{definition}\label{def:curve-congruence-equivalence}
Two oriented curves $\gamma:I\to\Ein$ and $\beta:J\to\Ein$ are \emph{congruent} if $I=J$ and there exists $\mathbf A\in\AWitt_+(2,3)$ such that $\gamma(s)=\mathbf A\cdot\beta(s)$ for every $s\in I$. They are \emph{equivalent} if there exist $\mathbf A\in\AWitt_+(2,3)$ and an orientation-preserving diffeomorphism $f:I\to J$ such that $\gamma(s)=\mathbf A\cdot\beta(f(s))$ for every $s\in I$. If both curves are parametrized by their natural parameters, then $f'(s)=1$, so that $f(s)=s+s_0$. Thus, for naturally parametrized curves, equivalence amounts to congruence up to a translation of the natural parameter.
\end{definition}

\medskip
\begin{proposition}\label{prop:canonical-moving-frame}
Let $\gamma:I\to\Ein$ be an oriented spacelike curve of helicity $\varepsilon$, parametrized by its natural parameter $s$. There exist unique smooth functions $\kappa_\gamma,\tau_\gamma:I\to\mathbb R$ and a unique smooth map $\mathbf A_\gamma=(A_1,\ldots,A_5):I\to\AWitt_+(2,3)$ such that $[A_1]=\gamma$ and
\begin{equation}\label{eq:frenet-system}
d\mathbf A_\gamma=\mathbf A_\gamma\mathbf K_\varepsilon(\kappa_\gamma,\tau_\gamma)\,ds,
\end{equation}
where $\mathbf K_\varepsilon:\mathbb R^2\to\mathfrak a(2,3)$ is defined by
\begin{equation}\label{eq:frenet-matrix}
\mathbf K_\varepsilon(\kappa,\tau)=
\begin{pmatrix}
0&-\kappa&-a_\varepsilon&b_\varepsilon&0\\
1&0&0&0&\kappa\\
0&0&0&\tau&a_\varepsilon\\
0&0&\tau&0&b_\varepsilon\\
0&-1&0&0&0
\end{pmatrix}.
\end{equation}
Conversely, given an interval $I\subset\mathbb R$, a sign $\varepsilon\in\{-1,1\}$, and two smooth functions $f,g\in C^\infty(I)$, there exists an oriented spacelike curve $\gamma:I\to\Ein$, parametrized by its natural parameter, of helicity $\varepsilon$, such that $\kappa_\gamma=f$ and $\tau_\gamma=g$. Such a curve is unique up to $\AWitt_+(2,3)$-congruence.
\end{proposition}

\begin{proof}
Throughout the proof, a prime denotes differentiation with respect to $s$.

Let $A_1$ be the unique positive null lift of $\gamma$ satisfying $\langle A_1',A_1'\rangle=1$, and set $A_2=A_1'$. Differentiating the identities $\langle A_1,A_1\rangle=0$ and $\langle A_2,A_2\rangle=1$, we obtain $\langle A_1,A_1''\rangle=-1$ and $\langle A_2,A_1''\rangle=0$.

Set $\kappa_\gamma=\langle A_1'',A_1''\rangle/2$ and define $A_5=-A_1''-\kappa_\gamma A_1$. Then $A_5$ is null, $\langle A_1,A_5\rangle=1$, and $\langle A_2,A_5\rangle=0$. Thus $(A_1,A_2,A_5)$ is a Witt basis of $\mathcal O_2(\gamma)$, and $A_1'=A_2$, $A_2'=-\kappa_\gamma A_1-A_5$.

Observe that $\pi_{\mathcal N}(A_5')=-\pi_{\mathcal N}(A_1''')$. The normalization of the natural parameter then gives $\langle\pi_{\mathcal N}(A_5'),\pi_{\mathcal N}(A_5')\rangle=\varepsilon$.

If $\varepsilon=1$, set $A_3=\pi_{\mathcal N}(A_5')$; if $\varepsilon=-1$, set $A_4=\pi_{\mathcal N}(A_5')$. In either case, the remaining vector is uniquely determined by the inner-product and orientation conditions requiring that $\mathbf A_\gamma=(A_1,A_2,A_3,A_4,A_5)$ belong to $\AWitt_+(2,3)$. This defines a smooth moving frame satisfying $\pi_{\mathcal N}(A_5')=a_\varepsilon A_3+b_\varepsilon A_4$.

Differentiating the scalar products between the frame vectors, we obtain $A_3'=-a_\varepsilon A_1+\tau_\gamma A_4$, $A_4'=b_\varepsilon A_1+\tau_\gamma A_3$, and $A_5'=\kappa_\gamma A_2+a_\varepsilon A_3+b_\varepsilon A_4$, where $\tau_\gamma$ is a smooth function. Together with the equations for $A_1'$ and $A_2'$, these identities give the Frenet system~\eqref{eq:frenet-system}.

For uniqueness, let $\widetilde{\mathbf A}$ be another $\AWitt_+(2,3)$-valued lift of $\gamma$ satisfying a Frenet system of the prescribed form. The first Frenet equation and the normalization of the positive null lift imply $\widetilde A_1=A_1$ and $\widetilde A_2=A_2$. The second equation, together with the nullity of $\widetilde A_5$, then gives $\widetilde A_5=A_5$. The fifth equation determines $\widetilde A_3=A_3$ when $\varepsilon=1$, and $\widetilde A_4=A_4$ when $\varepsilon=-1$. The remaining vector is fixed by the inner-product and orientation conditions. Thus $\widetilde{\mathbf A}=\mathbf A_\gamma$, and the uniqueness of $\kappa_\gamma$ and $\tau_\gamma$ follows from the Frenet equations.

Conversely, given $\varepsilon\in\{-1,1\}$ and $f,g\in C^\infty(I)$, consider the linear system $\mathbf A'=\mathbf A\mathbf K_\varepsilon(f,g)$. For any initial condition in $\AWitt_+(2,3)$, this system admits a unique solution $\mathbf A=(A_1,\ldots,A_5)$ with values in $\AWitt_+(2,3)$. The Frenet equations show that $\gamma=[A_1]$ is a spacelike curve parametrized by its natural parameter, with helicity $\varepsilon$, conformal curvature $f$, and conformal torsion $g$. Uniqueness up to the action of $\AWitt_+(2,3)$ follows from the uniqueness of solutions of the linear system.
\end{proof}

\begin{definition}\label{def:canonical-moving-frame}
The map $\mathbf A_\gamma$ is called the \emph{canonical moving frame} of $\gamma$, and $\mathbf K_\varepsilon(\kappa_\gamma,\tau_\gamma)$ its \emph{Frenet matrix}. The functions $\kappa_\gamma$ and $\tau_\gamma$ are called the \emph{conformal curvature} and \emph{conformal torsion} of $\gamma$, respectively. The map $\mathcal A_\gamma:I\to\AWitt_+(2,3)\times\mathbb R^3$, defined by $\mathcal A_\gamma(s)=(\mathbf A_\gamma(s),\kappa_\gamma(s),\tau_\gamma(s),\tau_\gamma'(s))$, is called the \emph{extended canonical frame} of $\gamma$.
\end{definition}

Let $\gamma$ be a generic spacelike curve of helicity $\varepsilon$, and set $\widetilde\gamma(s)=\gamma(-s)$. Its canonical frame is $\mathbf A_{\widetilde\gamma}(s)=\mathbf A_\gamma(-s)D_\varepsilon$, where $D_{+1}=\operatorname{diag}(1,-1,-1,1,1)\in\AWitt_+^\uparrow(2,3)$ and $D_{-1}=\operatorname{diag}(1,-1,1,-1,1)\in\AWitt_+^\downarrow(2,3)$. Thus, reversal of orientation preserves the component of the canonical frame when $\varepsilon=1$, whereas it interchanges the two time-orientation components when $\varepsilon=-1$. In the negative-helicity case, precisely one of the two orientations of the underlying curve is therefore positive.

\begin{definition}\label{def:positive-orientation}
An oriented generic spacelike curve of negative helicity is said to be \emph{positively oriented} if its canonical frame is time-oriented, that is, if $\mathbf A_\gamma(s)\in\AWitt_+^\uparrow(2,3)$.
\end{definition}

\begin{definition}\label{def:monodromic}
An oriented spacelike curve $\gamma:\mathbb R\to\Ein$, parametrized by its natural parameter, is called \emph{monodromic} if its conformal curvature $\kappa_\gamma$ and torsion $\tau_\gamma$ admit a common positive period. Suppose that these invariants are not both constant, and let $T>0$ be their minimal common period. The periodicity of the Frenet system implies the existence of a unique element $\mathbf M_\gamma\in\AWitt_+(2,3)$ such that $\mathbf A_\gamma(s+T)=\mathbf M_\gamma\mathbf A_\gamma(s)$. The element $\mathbf M_\gamma$ is called the \emph{monodromy} of $\gamma$.
\end{definition}

\begin{remark}\label{rem:monodromic-properties}
An oriented spacelike curve $\gamma:\mathbb R\to\Ein$, parametrized by its natural parameter, is monodromic if and only if there exist $P>0$ and $\mathbf M\in\AWitt_+(2,3)$ such that $\gamma(s+P)=\mathbf M\cdot\gamma(s)$ for all $s$. Indeed, the periodicity of the conformal invariants implies the existence of such an element by the Frenet system. Conversely, whenever this relation holds, the equivariance and uniqueness of the canonical moving frame imply that $\kappa_\gamma$ and $\tau_\gamma$ are $P$-periodic. For a monodromic curve with conformal invariants not both constant, the monodromy satisfies $\mathbf M_{\mathbf B\cdot\gamma}=\mathbf B\mathbf M_\gamma\mathbf B^{-1}$. Thus, its conjugacy class is invariant under $\AWitt_+(2,3)$-congruence. Finally, such a curve is closed if and only if its monodromy has finite order. Indeed, every period of $\gamma$ is an integer multiple of $T$, and the uniqueness of the canonical moving frame implies that $\gamma(s+kT)=\gamma(s)$ if and only if $\mathbf M_\gamma^k=\mathbf I$. More precisely, if $\mathbf M_\gamma$ has order $N$, then $NT$ is the least positive period of $\gamma$.
\end{remark}

\subsection{The momentum space and the canonical Pfaffian system}\label{subsec:canonical-pfaffian}

The $13$-dimensional manifold
$\mathcal Y=\AWitt_+(2,3)\times\mathbb R^3$, with coordinates
$(\kappa,\tau,\tau_1)$ on the second factor, provides the natural ambient
space for the extended canonical frames introduced in
Definition~\ref{def:canonical-moving-frame}. We shall refer to
$\mathcal Y$ as the \emph{momentum space}.

For $\varepsilon\in\{-1,1\}$, we define a rank-ten Pfaffian system
$\mathcal D_\varepsilon\subset T^*\mathcal Y$ by
$\mathcal D_\varepsilon=\operatorname{span}\{\eta_\varepsilon^1,\ldots,\eta_\varepsilon^{10}\}$, where
\begin{equation}\label{eq:canonical-pfaffian-forms}
\left\{
\begin{array}{@{}l@{}}
(\eta_\varepsilon^1,\eta_\varepsilon^5,\eta_\varepsilon^6,\eta_\varepsilon^7,\eta_\varepsilon^9)
=(\omega^1{}_1,\omega^2{}_3,\omega^2{}_4,\omega^3{}_1,\omega^4{}_1),\\[5pt]
(\eta_\varepsilon^2,\eta_\varepsilon^3,\eta_\varepsilon^4)
=(\omega^1{}_2+\kappa\omega,\ \omega^1{}_3+a_\varepsilon\omega,\ \omega^1{}_4-b_\varepsilon\omega),\\[5pt]
(\eta_\varepsilon^8,\eta_\varepsilon^{10})
=(\omega^3{}_4-\tau\omega,\ d\tau-\tau_1\omega).
\end{array}
\right.
\end{equation}
Here $\omega^i{}_j$ denote the components of the Maurer--Cartan form of $\AWitt_+(2,3)$,
and $\omega=\omega^2{}_1$. The independence condition is $\omega\neq0$.
The Pfaffian system $\mathcal D_\varepsilon$, together with this independence
condition, is called the \emph{canonical Pfaffian system}.

\begin{property}
For each $\varepsilon\in\{-1,1\}$, the integral curves of the canonical
Pfaffian system $\mathcal D_\varepsilon$ satisfying the independence
condition $\omega\neq0$ are, after parametrization by $\omega=ds$,
precisely the extended canonical frames of nondegenerate spacelike
curves of helicity $\varepsilon$.
\end{property}

\begin{proof}
Let $\mathcal A=(\mathbf A,\kappa,\tau,\tau_1):I\to\mathcal Y$ be a
smooth curve. The Pfaffian equations
$\mathcal A^*\eta_\varepsilon^j=0$, $1\leq j\leq10$, are equivalent to $\mathbf A^{-1}d\mathbf A=\mathbf K_\varepsilon(\kappa,\tau)\,\omega$
and $d\tau=\tau_1\omega$, where $\omega$ denotes its pullback to $I$. If the independence
condition holds, there exists a parameter $s$, unique up to translations,
such that $\omega=ds$. By Proposition~\ref{prop:canonical-moving-frame},
with respect to this parameter, $\mathcal A$ is the extended canonical
frame of the nondegenerate spacelike curve $\gamma=[A_1]$ of helicity
$\varepsilon$. Conversely, every such extended canonical frame satisfies
the Pfaffian equations and the independence condition.
\end{proof}

The regularity of integral curves of the canonical Pfaffian system will play
an important role in the variational analysis developed in the next section.
We adopt Hsu's characterization of regularity in terms of characteristic lifts
\cite[Theorem~6]{Hsu}, which he proved equivalent to his original definition
based on the holonomy map \cite[equation~(2.6), the subsequent definition and
Appendix~B]{Hsu}. See also
\cite[Proposition~5.6 and Theorem~5.7]{CGR} and
\cite[Proposition~2.2]{dPS}. For completeness, we begin by recalling the
relevant definitions.

\bigskip

Let $N$ be a smooth manifold and let $\Sigma\in\Omega^2(N)$. The
\emph{Cartan system} associated with $\Sigma$ is defined by
$\mathcal C_\Sigma=\operatorname{span}\{\iota_X\Sigma:X\in\mathfrak X(N)\}$.
A smooth curve $\sigma:I\to N$ is called a \emph{characteristic curve} of
$\Sigma$ if it is an integral curve of $\mathcal C_\Sigma$, or equivalently,
if $\iota_{\sigma'}\Sigma=0$.

Now let $D\subset T^*M$ be a Pfaffian system of constant rank $r$, and let
$\zeta_0$ denote the restriction to $D$ of the tautological $1$-form of $T^*M$.
In terms of a local frame $\theta^1,\ldots,\theta^r$ of $D$ and fiber
coordinates $\lambda_1,\ldots,\lambda_r$, we have
$\zeta_0=\lambda_j\theta^j$, with summation over $j=1,\ldots,r$.
Contracting $d\zeta_0$ with the vertical vector fields $\partial_{\lambda_j}$
shows that every characteristic curve of $d\zeta_0$ projects onto an integral
curve of $D$. Such a characteristic curve is called a \emph{characteristic lift}
of its projection.

\begin{definition}
An integral curve $\sigma:I\to M$ of $D$ is called \emph{regular} if, on
every compact subinterval $J\subset I$, the zero lift is the only characteristic
lift of $\sigma|_J$. Otherwise, $\sigma$ is called \emph{singular}.
\end{definition}

We now apply this definition to the canonical Pfaffian system.

\begin{proposition}\label{prop:regularity}
Every integral curve $\mathcal A$ of the canonical Pfaffian system
$\mathcal D_\varepsilon$ satisfying the independence condition $\omega\neq0$
is regular.
\end{proposition}

\begin{proof}
Let $\zeta_0=\lambda_j\eta_\varepsilon^j$ be the restriction of the tautological
$1$-form to the total space of $\mathcal D_\varepsilon$.
The $1$-forms $\bigl(\omega,\eta_\varepsilon^1,\ldots,\eta_\varepsilon^9,
d\kappa,d\tau,d\tau_1\bigr)$ constitute a global coframe on $\mathcal Y$.
We denote its dual frame by $\bigl(\partial_\omega,
\partial_{\eta_\varepsilon^1},\ldots,\partial_{\eta_\varepsilon^9},
\partial_\kappa,\partial_\tau,\partial_{\tau_1}\bigr)$.
On the total space of $\mathcal D_\varepsilon$, we complete this coframe
and its dual frame by adjoining $d\lambda_1,\ldots,d\lambda_{10}$ and
$\partial_{\lambda_1},\ldots,\partial_{\lambda_{10}}$, respectively.

In the following computations, the symbol $\equiv$ denotes equality modulo the $2$-forms $\eta_\varepsilon^i\wedge\eta_\varepsilon^j$, $1\leq i,j\leq10$. The Maurer--Cartan equations give the following structure equations:
{\small
\begin{equation}\label{eq:pfaffian-structure-equations}
\left\{
\begin{alignedat}{2}
d\eta_\varepsilon^1&\equiv\omega\wedge(\eta_\varepsilon^2+a_\varepsilon\eta_\varepsilon^7-b_\varepsilon\eta_\varepsilon^9),
&\;d\eta_\varepsilon^{8\phantom{0}}&\equiv\omega\wedge(\tau\eta_\varepsilon^1+b_\varepsilon\eta_\varepsilon^7-a_\varepsilon\eta_\varepsilon^9+\eta_\varepsilon^{10}),\\[3pt]
d\eta_\varepsilon^2&\equiv-\omega\wedge(d\kappa+2\kappa\eta_\varepsilon^1+a_\varepsilon\eta_\varepsilon^5+b_\varepsilon\eta_\varepsilon^6),
&\;d\eta_\varepsilon^{7\phantom{0}}&\equiv-\omega\wedge(\eta_\varepsilon^5+\tau\eta_\varepsilon^9),\\[3pt]
d\eta_\varepsilon^3&\equiv\omega\wedge(-2a_\varepsilon\eta_\varepsilon^1+\tau\eta_\varepsilon^4+\kappa\eta_\varepsilon^5-b_\varepsilon\eta_\varepsilon^8),
&\;d\eta_\varepsilon^{9\phantom{0}}&\equiv\omega\wedge(\eta_\varepsilon^6-\tau\eta_\varepsilon^7),\\[3pt]
d\eta_\varepsilon^4&\equiv\omega\wedge(2b_\varepsilon\eta_\varepsilon^1+\tau\eta_\varepsilon^3+\kappa\eta_\varepsilon^6+a_\varepsilon\eta_\varepsilon^8),
&\;d\eta_\varepsilon^{10}&\equiv\omega\wedge(d\tau_1+\tau_1\eta_\varepsilon^1),\\[3pt]
d\eta_\varepsilon^5&\equiv\omega\wedge(-\eta_\varepsilon^3+\tau\eta_\varepsilon^6+\kappa\eta_\varepsilon^7),
&\;d\eta_\varepsilon^{6\phantom{0}}&\equiv\omega\wedge(-\eta_\varepsilon^4+\tau\eta_\varepsilon^5-\kappa\eta_\varepsilon^9),
\end{alignedat}
\right.
\end{equation}
}
Then
\begin{equation}\label{eq:dzeta0}
d\zeta_0\equiv d\lambda_j\wedge\eta_\varepsilon^j+\omega\wedge\upsilon,
\end{equation}
where
\begin{equation}\label{eq:upsilon}
\begin{aligned}
\upsilon={}&(-2\kappa\lambda_2-2a_\varepsilon\lambda_3+2b_\varepsilon\lambda_4
+\tau\lambda_8+\tau_1\lambda_{10})\eta_\varepsilon^1\\
&+\lambda_1\eta_\varepsilon^2+(\tau\lambda_4-\lambda_5)\eta_\varepsilon^3
+(\tau\lambda_3-\lambda_6)\eta_\varepsilon^4\\
&+(-a_\varepsilon\lambda_2+\kappa\lambda_3+\tau\lambda_6-\lambda_7)\eta_\varepsilon^5\\
&+(-b_\varepsilon\lambda_2+\kappa\lambda_4+\tau\lambda_5+\lambda_9)\eta_\varepsilon^6\\
&+(a_\varepsilon\lambda_1+\kappa\lambda_5+b_\varepsilon\lambda_8-\tau\lambda_9)\eta_\varepsilon^7
+(-b_\varepsilon\lambda_3+a_\varepsilon\lambda_4)\eta_\varepsilon^8\\
&+(-b_\varepsilon\lambda_1-\kappa\lambda_6-\tau\lambda_7-a_\varepsilon\lambda_8)\eta_\varepsilon^9\\
&+\lambda_8\eta_\varepsilon^{10}-\lambda_2\,d\kappa+\lambda_{10}\,d\tau_1.
\end{aligned}
\end{equation}

Let $\mathfrak A$ be a characteristic lift of an integral curve $\mathcal A$
satisfying $\omega\neq0$. We use the natural parameter $s$, so that
$\mathcal A^*\omega=ds$, and denote differentiation with respect to $s$ by
a prime. Since the tangent vector to $\mathfrak A$ annihilates the forms
$\eta_\varepsilon^j$, the terms involving
$\eta_\varepsilon^i\wedge\eta_\varepsilon^j$ do not contribute
to the characteristic equations.

By contracting with $\partial_\kappa$ and $\partial_{\tau_1}$, we obtain
$\lambda_2=\lambda_{10}=0$. Using
$\eta_\varepsilon^{10}=d\tau-\tau_1\omega$, contraction with $\partial_\tau$
gives $\lambda_{10}'=-\lambda_8$. Since $\lambda_{10}$ vanishes identically,
we obtain $\lambda_8=0$. The contractions with
$\partial_{\eta_\varepsilon^2}$, $\partial_{\eta_\varepsilon^8}$, and
$\partial_{\eta_\varepsilon^1}$ give
$\lambda_2'=-\lambda_1$,
$\lambda_8'=b_\varepsilon\lambda_3-a_\varepsilon\lambda_4$, and
$\lambda_1'=2a_\varepsilon\lambda_3-2b_\varepsilon\lambda_4$.
Since $\lambda_2=\lambda_8=0$, these relations imply $\lambda_1=0$,
$b_\varepsilon\lambda_3-a_\varepsilon\lambda_4=0$, and
$a_\varepsilon\lambda_3-b_\varepsilon\lambda_4=0$.
As $a_\varepsilon^2-b_\varepsilon^2=\varepsilon\neq0$, we conclude that
$\lambda_1=\lambda_3=\lambda_4=0$.

Finally, the contractions with $\partial_{\eta_\varepsilon^3}$ and
$\partial_{\eta_\varepsilon^4}$ give
$\lambda_3'=\lambda_5-\tau\lambda_4$ and
$\lambda_4'=\lambda_6-\tau\lambda_3$, respectively. Those with
$\partial_{\eta_\varepsilon^5}$ and $\partial_{\eta_\varepsilon^6}$ give
$\lambda_5'=\lambda_7-\kappa\lambda_3-\tau\lambda_6$ and
$\lambda_6'=-\lambda_9-\kappa\lambda_4-\tau\lambda_5$, respectively.
Since $\lambda_3=\lambda_4=0$, it follows successively that
$\lambda_5=\lambda_6=\lambda_7=\lambda_9=0$.

Thus $\mathfrak A$ is the zero lift, and $\mathcal A$ is regular.
\end{proof}

\medskip

\section{The conformal arclength functional}\label{sec:conformal-variational-system}

Using the formalism of~\cite{Griffiths}, together with the canonical Pfaffian system and its regularity established in Section~1, we identify the extremals of the conformal arclength functional with the trajectories of the Reeb flow on the momentum space. We then construct the moment map, analyze the reduced phase curves, classify closed extremals with constant conformal invariants, and establish noncommutative integrability.

\subsection{The conformal arclength functional}\label{subsec:arclength-functional}

\begin{definition}
The \emph{conformal arclength functional} of a nondegenerate space-like curve $\gamma:I\to\Ein$ is given by
\[
\mathcal L(\gamma)=\int_I ds_\gamma,
\]
whenever the integral is finite.
\end{definition}

\begin{definition}
A nondegenerate space-like curve is called a \emph{conformal extremal} if its first variation vanishes under every compactly supported variation through nondegenerate space-like curves.
\end{definition}

In taking variations, the parameter is kept fixed but is not required to be a natural parameter for the varied curves. When the total conformal arclength is infinite, first variations are understood locally, by integrating over a compact subinterval containing the support of the variation.

We briefly recall the formalism for constrained variational problems. Let $M$ be a smooth manifold, $D\subset T^*M$ a Pfaffian system of constant rank, and $\varphi\in\Omega^1(M)$. The \emph{constrained variational problem} associated with $(D,\varphi)$ consists in studying the critical points of
\[
\sigma\longmapsto\int_\sigma\varphi
\]
on the integral curves of $D$, with respect to compactly supported variations through such curves.

Following~\cite{Griffiths}, we associate with this problem the affine bundle $Z=\varphi+D\subset T^*M$. Let $\zeta$ denote the restriction to $Z$ of the tautological $1$-form of $T^*M$. A characteristic lift of an integral curve of $D$ is a characteristic curve of $d\zeta$ projecting onto it. Combining~\cite[Proposition~3]{Bryant} with~\cite[Theorem~10]{Hsu}, we obtain that a regular integral curve of $D$ is critical for the constrained variational problem if and only if it admits a characteristic lift to $Z$. When it exists, the characteristic lift is unique by regularity.

In the present setting, the constrained variational problem associated with the conformal arclength functional is defined by
\[
\A\longmapsto\int_{\A}\omega
\]
on the integral curves of the canonical Pfaffian system $\mathcal D_\varepsilon$ satisfying the independence condition $\omega\neq0$.

\begin{proposition}\label{prop:variational-correspondence}
Let $\gamma:I\to\Ein$ be a nondegenerate space-like curve of helicity $\varepsilon$, and let $\A_\gamma$ be its extended canonical frame. Then $\gamma$ is a conformal extremal if and only if $\A_\gamma$ is critical for the constrained variational problem $(\mathcal D_\varepsilon,\omega)$.
\end{proposition}

\begin{proof}
Let $\gamma_t$ be a compactly supported variation of $\gamma$ through nondegenerate space-like curves, and write $\A_{\gamma_t}=(\mathbf A_{\gamma_t},\kappa_{\gamma_t},\tau_{\gamma_t},\tau_{1,\gamma_t})$, where $\tau_{1,\gamma_t}=d\tau_{\gamma_t}/ds_{\gamma_t}$.
The smooth dependence of the canonical frames and invariants on the jets of the curves shows that $\A_{\gamma_t}$ defines a compactly supported variation of $\A_\gamma$ through integral curves of $\mathcal D_\varepsilon$ satisfying $\omega\neq0$.

Conversely, every such variation of $\A_\gamma$ projects onto a compactly supported variation of $\gamma$ through nondegenerate space-like curves, whose extended canonical frames are precisely the curves $\A_{\gamma_t}$. Moreover, $\A_{\gamma_t}^*\omega=ds_{\gamma_t}$, so the two functionals agree on every compact subinterval containing the support of the variation. Consequently, their critical points correspond.
\end{proof}

We now introduce a contact structure on the momentum space.

\begin{lemma}\label{lem:contact-reeb}
For each $\varepsilon\in\{-1,1\}$, let
$\alpha_\varepsilon\in\Omega^1(\mathcal Y)$ be the $1$-form
\begin{equation}\label{eq:contact-form}
\begin{aligned}
\alpha_\varepsilon={}&\omega+\frac12\bigl[
b_\varepsilon(\eta_\varepsilon^4+\tau\eta_\varepsilon^5)
-a_\varepsilon(\eta_\varepsilon^3+\tau\eta_\varepsilon^6)\\
&+\tau_1(b_\varepsilon\eta_\varepsilon^7+a_\varepsilon\eta_\varepsilon^9)
-(\kappa+\tau^2)(a_\varepsilon\eta_\varepsilon^7+b_\varepsilon\eta_\varepsilon^9)
\bigr].
\end{aligned}
\end{equation}
Then $\alpha_\varepsilon$ is a contact form on $\mathcal Y$.
Its Reeb vector field is
\begin{equation}\label{eq:reeb-vector-field}
R_\varepsilon=
\partial_\omega
-3\tau\tau_1\partial_\kappa
+\tau_1\partial_\tau
-\tau(2\kappa+\tau^2)\partial_{\tau_1}.
\end{equation}
\end{lemma}

\begin{proof}
A direct computation using the Maurer--Cartan equations shows that the
$C^\infty(\mathcal Y)$-module
$\{\iota_Xd\alpha_\varepsilon:X\in\mathfrak X(\mathcal Y)\}$ is generated by
the twelve linearly independent $1$-forms
$\eta_\varepsilon^1,\ldots,\eta_\varepsilon^9$,
$d\kappa+3\tau\tau_1\omega$,
$d\tau-\tau_1\omega$, and
$d\tau_1+\tau(2\kappa+\tau^2)\omega$.
Consequently, $d\alpha_\varepsilon$ has constant rank twelve, and its
kernel is spanned by the vector field $R_\varepsilon$ defined in
\eqref{eq:reeb-vector-field}. Since $\alpha_\varepsilon(R_\varepsilon)=1$,
the $1$-form $\alpha_\varepsilon$ is a contact form and
$R_\varepsilon$ is its Reeb vector field.
\end{proof}

\begin{theorem}\label{thm:conformal-extremals-reeb}
A nondegenerate spacelike curve $\gamma:I\to\Ein$ of helicity
$\varepsilon$, parametrized by conformal arclength, is a conformal
extremal if and only if its extended canonical frame
$\A_\gamma:I\to\mathcal Y$ is an integral curve of the Reeb vector
field $R_\varepsilon$. Equivalently, the conformal extremals of
helicity $\varepsilon$, parametrized by conformal arclength, are
precisely the projections onto $\Ein$ of the integral curves of
$R_\varepsilon$.
\end{theorem}

\begin{proof}
Let $\mathcal Z_\varepsilon=\omega+\mathcal D_\varepsilon\subset T^*\mathcal Y$
be Griffiths' affine bundle, and let $\zeta_\varepsilon$ denote the
restriction of the tautological $1$-form to $\mathcal Z_\varepsilon$.
On this space, write $\zeta_\varepsilon=\omega+\lambda_j\eta_\varepsilon^j$.
The Maurer--Cartan equations give $d\omega\equiv-\omega\wedge\eta_\varepsilon^1$. Together with \eqref{eq:dzeta0}, this yields $d\zeta_\varepsilon\equiv d\lambda_j\wedge\eta_\varepsilon^j+\omega\wedge(\upsilon-\eta_\varepsilon^1)$, where $\upsilon$ is the $1$-form introduced in the proof of Proposition~\ref{prop:regularity}.

By Propositions~\ref{prop:variational-correspondence} and
\ref{prop:regularity} and the variational criterion in~\cite{Bryant,Hsu},
a nondegenerate spacelike curve $\gamma$ is a conformal extremal if and
only if its extended canonical frame $\A_\gamma$ admits a characteristic
lift to $\mathcal Z_\varepsilon$.

Let $\mathfrak A$ be a characteristic curve of $d\zeta_\varepsilon$
satisfying $\omega\neq0$. Contraction with the vertical vector fields
shows that its projection is an integral curve of $\mathcal D_\varepsilon$.
Thus, along $\mathfrak A$, the Pfaffian equations
$\eta_\varepsilon^j=0$, $1\leq j\leq10$, are satisfied, including
$d\tau=\tau_1\omega$.

Proceeding as in the proof of Proposition~\ref{prop:regularity}, the
remaining characteristic equations determine the multipliers
$\lambda_j=\Lambda_j^\varepsilon(\kappa,\tau,\tau_1)$,
$1\leq j\leq10$, where
\[
\begin{aligned}
\Lambda_1^\varepsilon
&=\Lambda_2^\varepsilon
=\Lambda_8^\varepsilon
=\Lambda_{10}^\varepsilon=0,\\
\Lambda_3^\varepsilon
&=-a_\varepsilon/2,
\quad
\Lambda_4^\varepsilon=b_\varepsilon/2,
\quad
\Lambda_5^\varepsilon=b_\varepsilon\tau/2,
\quad
\Lambda_6^\varepsilon=-a_\varepsilon\tau/2,\\
\Lambda_7^\varepsilon
&=\bigl[b_\varepsilon\tau_1-a_\varepsilon(\kappa+\tau^2)\bigr]/2,
\quad
\Lambda_9^\varepsilon=\bigl[a_\varepsilon\tau_1-b_\varepsilon(\kappa+\tau^2)\bigr]/2.
\end{aligned}
\]
The remaining compatibility conditions are equivalent, since
$a_\varepsilon^2-b_\varepsilon^2=\varepsilon\neq0$, to
\[
d\kappa=-3\tau\tau_1\omega,
\qquad
d\tau_1=-\tau(2\kappa+\tau^2)\omega.
\]

The multipliers define a smooth global section
$j_\varepsilon:\mathcal Y\to\mathcal Z_\varepsilon$ by
\[
j_\varepsilon(\mathbf A,\kappa,\tau,\tau_1)
=\bigl(\mathbf A,\kappa,\tau,\tau_1,
\Lambda^\varepsilon(\kappa,\tau,\tau_1)\bigr).
\]
By the expressions of the multipliers, the pullback of the tautological form satisfies $j_\varepsilon^*\zeta_\varepsilon=\alpha_\varepsilon$, where $\alpha_\varepsilon$ is the contact form defined in \eqref{eq:contact-form}.

The characteristic equations, together with Lemma~\ref{lem:contact-reeb},
show that the characteristic curves of $d\zeta_\varepsilon$ satisfying
$\omega\neq0$ are precisely, after natural parametrization, the images
under $j_\varepsilon$ of the integral curves of $R_\varepsilon$.
The variational characterization follows.
\end{proof}

\begin{corollary}\label{cor:EL-equations}
Let $\gamma:I\to\Ein$ be a nondegenerate spacelike curve of helicity
$\varepsilon$, parametrized by conformal arclength. Then $\gamma$ is
a conformal extremal if and only if its conformal curvature
$\kappa_\gamma$ and torsion $\tau_\gamma$ satisfy the reduced equations
\begin{equation}\label{eq:reeb-dynamics}
\kappa_\gamma'=-3\tau_\gamma\tau_\gamma',
\qquad
\tau_\gamma''=-\tau_\gamma(2\kappa_\gamma+\tau_\gamma^2).
\end{equation}
Conversely, for each $\varepsilon\in\{-1,1\}$, every smooth solution
of \eqref{eq:reeb-dynamics} determines a conformal extremal of helicity
$\varepsilon$, uniquely up to $\AWitt_+(2,3)$-congruence.
\end{corollary}

\begin{proof}
By Theorem~\ref{thm:conformal-extremals-reeb}, $\gamma$ is a
conformal extremal if and only if its extended canonical frame is an
integral curve of $R_\varepsilon$. In view of
\eqref{eq:reeb-vector-field}, this is equivalent to
\eqref{eq:reeb-dynamics}. The converse follows from
Proposition~\ref{prop:canonical-moving-frame}.
\end{proof}

\subsection{The moment map, reduced spaces, and phase curves}\label{subsec:contact-moment-map}

The natural left action of $\AWitt_+(2,3)$ on the momentum space
$\mathcal Y$ preserves the contact form $\alpha_\varepsilon$ and hence admits
a contact moment map
$\widehat J_\varepsilon:\mathcal Y\to\mathfrak a(2,3)^*$
(see \cite{Loose,Willett}).

Since $\mathfrak a(2,3)$ is semisimple, its Killing form is nondegenerate.
We identify $\mathfrak a(2,3)^*$ with $\mathfrak a(2,3)$ by means of the
Killing form normalized by
$\langle\mathbf X,\mathbf Y\rangle=\operatorname{tr}(\mathbf X\mathbf Y)/4$.
Under this identification, $\widehat J_\varepsilon$ is represented by a map
$\mathbf J_\varepsilon:\mathcal Y\to\mathfrak a(2,3)$.
Our aim is to determine $\mathbf J_\varepsilon$ explicitly.

\begin{definition}
For each $\varepsilon\in\{-1,1\}$, the \emph{Hamiltonian matrix} of the
conformal arclength problem is the matrix-valued map
$\mathbf H_\varepsilon:\mathbb R^3\to\mathfrak a(2,3)$ defined by
\begin{equation}\label{eq:hamiltonian-matrix}
\mathbf H_\varepsilon(\kappa,\tau,\tau_1)=
\begin{pmatrix}
0&1&u_\varepsilon&v_\varepsilon&0\\
0&0&-b_\varepsilon\tau&-a_\varepsilon\tau&-1\\
-a_\varepsilon&b_\varepsilon\tau&0&0&-u_\varepsilon\\
b_\varepsilon&-a_\varepsilon\tau&0&0&v_\varepsilon\\
0&0&a_\varepsilon&b_\varepsilon&0
\end{pmatrix},
\end{equation}
where $u_\varepsilon=b_\varepsilon\tau_1-a_\varepsilon(\kappa+\tau^2)$ and
$v_\varepsilon=a_\varepsilon\tau_1-b_\varepsilon(\kappa+\tau^2)$.
\end{definition}

\begin{proposition}\label{prop:moment-map}
Under the identification of $\mathfrak a(2,3)^*$ with
$\mathfrak a(2,3)$ defined above, the contact moment map is represented by
\begin{equation}\label{eq:moment-map}
\mathbf J_\varepsilon(\mathcal A)
=\mathbf A\mathbf H_\varepsilon(\kappa,\tau,\tau_1)\mathbf A^{-1},
\qquad \mathcal A=(\mathbf A,\kappa,\tau,\tau_1).
\end{equation}
\end{proposition}

\begin{proof}
Let $\mathbf X\in\mathfrak a(2,3)$, and denote by $X_{\mathcal Y}$ the
corresponding fundamental vector field. By definition of the contact
moment map, $\widehat J_\varepsilon(\mathcal A)(\mathbf X)
=\alpha_\varepsilon(X_{\mathcal Y}(\mathcal A))$.
For $\mathcal A=(\mathbf A,\kappa,\tau,\tau_1)$, we have
$X_{\mathcal Y}(\mathcal A)=(\mathbf X\mathbf A,0,0,0)$, and hence
$\boldsymbol\Omega(X_{\mathcal Y})=\mathbf A^{-1}\mathbf X\mathbf A$.
Substituting into \eqref{eq:contact-form} and using the cyclic property
of the trace, we obtain
\[
\alpha_\varepsilon(X_{\mathcal Y})
=\frac14\operatorname{tr}(\mathbf H_\varepsilon\mathbf A^{-1}\mathbf X\mathbf A)
=\frac14\operatorname{tr}(\mathbf A\mathbf H_\varepsilon\mathbf A^{-1}\mathbf X).
\]
The assertion follows from the identification of $\mathfrak a(2,3)^*$
with $\mathfrak a(2,3)$.
\end{proof}

\begin{remark}
Since the action of $\AWitt_+(2,3)$ preserves the contact form
$\alpha_\varepsilon$, the moment map $\mathbf J_\varepsilon$ is constant
along the integral curves of the Reeb vector field $R_\varepsilon$.
Hence, by the preceding proposition, if $\gamma:I\to\Ein$ is a conformal
extremal of helicity $\varepsilon$, parametrized by conformal arclength,
with conformal curvature $\kappa_\gamma$ and torsion $\tau_\gamma$, then
we have
\begin{equation}\label{eq:conserved-momentum}
\mathbf J_\varepsilon(\mathcal A_\gamma(s))
=\mathbf A_\gamma(s)\mathbf H_\gamma(s)\mathbf A_\gamma(s)^{-1}
=\mathbf J_\gamma\in\mathfrak a(2,3),
\end{equation}
where $\mathbf H_\gamma(s)=\mathbf H_\varepsilon(\kappa_\gamma(s),\tau_\gamma(s),\tau_\gamma'(s))$,
$\mathcal A_\gamma$ is the extended canonical frame of $\gamma$, and
$\mathbf J_\gamma$ is constant.
\end{remark}

\begin{definition}
The constant element $\mathbf J_\gamma\in\mathfrak a(2,3)$ appearing
in \eqref{eq:conserved-momentum} is called the \emph{momentum} of the
conformal extremal $\gamma$.
\end{definition}

Differentiating \eqref{eq:conserved-momentum} and using
$\mathbf A_\gamma'=\mathbf A_\gamma\mathbf K_\varepsilon(\kappa_\gamma,\tau_\gamma)$,
we obtain the Lax equation
\begin{equation}\label{eq:lax}
\mathbf H_\gamma'+[\mathbf K_\varepsilon(\kappa_\gamma,\tau_\gamma),\mathbf H_\gamma]=0.
\end{equation}

The Lie algebra $\mathfrak a(2,3)$ has rank two, and its algebra of
$\operatorname{Ad}$-invariant polynomials is generated by two algebraically
independent Casimir invariants of degrees two and four. For each
$\varepsilon\in\{-1,1\}$, we choose
\begin{equation}\label{eq:casimirs}
h_\varepsilon(\mathbf X)=\frac{\varepsilon}{4}\operatorname{tr}(\mathbf X^2),\quad
E_\varepsilon(\mathbf X)=\frac18\bigl(\operatorname{tr}(\mathbf X^2)\bigr)^2-\frac14\operatorname{tr}(\mathbf X^4)+\varepsilon.
\end{equation}
The associated \emph{Casimir map} is defined by
\begin{equation}\label{eq:casimir-map}
\mathfrak c_\varepsilon:\mathfrak a(2,3)\longrightarrow\mathbb R^2,
\qquad\mathbf X\longmapsto
\bigl(h_\varepsilon(\mathbf X),E_\varepsilon(\mathbf X)\bigr).
\end{equation}
Evaluating these invariants on the Hamiltonian matrix
\eqref{eq:hamiltonian-matrix}, we obtain
\begin{equation}\label{eq:casimir-values}
h_\varepsilon(\mathbf H_\varepsilon)=\kappa+3\tau^2/2,\qquad
E_\varepsilon(\mathbf H_\varepsilon)=\tau_1^2+2\kappa\tau^2+2\tau^4.
\end{equation}

Since these invariants are preserved under conjugation and the momentum
$\mathbf J_\gamma$ is constant along every conformal extremal, we obtain
the following two conservation laws for the Euler--Lagrange system
\eqref{eq:reeb-dynamics}:
\begin{equation}\label{eq:first-integrals}
h_\gamma=\kappa_\gamma+3\tau_\gamma^2/2,\qquad
E_\gamma=(\tau_\gamma')^2+2h_\gamma\tau_\gamma^2-\tau_\gamma^4,
\end{equation}
where $h_\gamma=h_\varepsilon(\mathbf J_\gamma)$ and
$E_\gamma=E_\varepsilon(\mathbf J_\gamma)$.

\begin{definition}
The constants $h_\gamma$ and $E_\gamma$ appearing in
\eqref{eq:first-integrals} are called the \emph{parameters} of the
conformal extremal $\gamma$.
\end{definition}

\begin{proposition}\label{prop:moment-map-regularity}
The moment map $\mathbf J_\varepsilon:\mathcal Y\to\mathfrak a(2,3)$
has the following properties:
\begin{enumerate}[label=(\roman*),itemsep=0pt,parsep=0pt]
\item For every $\mathcal A\in\mathcal Y$, the matrix
$\mathbf J_\varepsilon(\mathcal A)$ is cyclic and hence a regular element
of $\mathfrak a(2,3)$. Its centralizer is the two-dimensional abelian
Lie subalgebra spanned by $\mathbf J_\varepsilon(\mathcal A)$ and
$\mathbf J_\varepsilon(\mathcal A)^3$.
\item The critical locus of $\mathbf J_\varepsilon$ is
$\mathcal Y_*=\{(\mathbf A,\kappa,\tau,\tau_1)\in\mathcal Y:
\tau_1=0,\ \tau(2\kappa+\tau^2)=0\}$.
More precisely, $\operatorname{rank}d\mathbf J_\varepsilon(\mathcal A)=10$
if $\mathcal A\in\mathcal Y\setminus\mathcal Y_*$, and
$\operatorname{rank}d\mathbf J_\varepsilon(\mathcal A)=9$
if $\mathcal A\in\mathcal Y_*$.
\end{enumerate}
\end{proposition}

\begin{proof}
Since conjugation preserves cyclicity and conjugates centralizers, and
the equivariance of $\mathbf J_\varepsilon$ implies that the rank of
its differential is independent of $\mathbf A$, it suffices to work
at $\mathbf A=\mathbf I$, where
$\mathbf J_\varepsilon(\mathbf I,\kappa,\tau,\tau_1)
=\mathbf H_\varepsilon(\kappa,\tau,\tau_1)$.

Let $\mathfrak a_{\mathbb C}(2,3)=\mathfrak a(2,3)\otimes_{\mathbb R}\mathbb C$.
By the regularity criterion in~\cite[Theorem~0.1]{Kostant1963},
an element $\mathbf X\in\mathfrak a_{\mathbb C}(2,3)$ is regular if and
only if the differentials of the fundamental invariant polynomials are
linearly independent. Since their gradients span the same subspace as
$\mathbf X$ and $\mathbf X^3$, this is equivalent to the linear independence
of these two matrices. Moreover, by \cite[Theorem~1]{Konyaev2010}, regularity
in $\mathfrak a_{\mathbb C}(2,3)$ is equivalent to cyclicity in
$\mathfrak{gl}(5,\mathbb C)$. Recall that a matrix
$\mathbf X\in\mathfrak{gl}(5,\mathbb C)$ is \emph{cyclic} if there exists a
vector $v\in\mathbb C^5$ such that
$v,\mathbf Xv,\mathbf X^2v,\mathbf X^3v,\mathbf X^4v$ form a basis of
$\mathbb C^5$. Equivalently, its minimal polynomial has degree five.
These results apply to the real form $\mathfrak a(2,3)$. Indeed, for every
$\mathbf X\in\mathfrak a(2,3)$, we have
$\mathfrak z_{\mathfrak a_{\mathbb C}(2,3)}(\mathbf X)
=\mathfrak z_{\mathfrak a(2,3)}(\mathbf X)\otimes_{\mathbb R}\mathbb C$,
and cyclicity is preserved under complexification, since the minimal
polynomial is unchanged by extension of scalars.

Direct calculation gives $(\mathbf H_\varepsilon)_{21}=0$ and
$(\mathbf H_\varepsilon^3)_{21}=\varepsilon$. Since
$(\mathbf H_\varepsilon)_{12}=1$, the matrices $\mathbf H_\varepsilon$
and $\mathbf H_\varepsilon^3$ are linearly independent. Hence
$\mathbf H_\varepsilon$ is both regular and cyclic. Since the
centralizer of a cyclic endomorphism consists of polynomials in that
endomorphism, skew-adjointness implies that
$\mathfrak z_{\mathfrak a(2,3)}(\mathbf H_\varepsilon)
=\operatorname{span}_{\mathbb R}\{\mathbf H_\varepsilon,\mathbf H_\varepsilon^3\}$.
This proves (i).

For (ii), consider the Casimir map $\mathfrak c_\varepsilon$ defined in \eqref{eq:casimir-map}.
Since their gradients, with respect to the trace pairing, are
$\frac{\varepsilon}{2}\mathbf X$ and
$2\varepsilon h_\varepsilon(\mathbf X)\mathbf X-\mathbf X^3$,
respectively, they are linearly independent at $\mathbf H_\varepsilon$
by (i). Consequently,
$\ker(d\mathfrak c_\varepsilon|_{\mathbf H_\varepsilon})=[\mathfrak a(2,3),\mathbf H_\varepsilon]$,
which has dimension eight.

At $(\mathbf I,\kappa,\tau,\tau_1)$, the differential of
$\mathbf J_\varepsilon$ is given by
$d\mathbf J_\varepsilon(\mathbf X,\mathbf v)
=[\mathbf X,\mathbf H_\varepsilon]+d\mathbf H_\varepsilon(\mathbf v)$,
where $\mathbf X\in\mathfrak a(2,3)$ and $\mathbf v\in\mathbb R^3$.
The first term spans the eight-dimensional subspace
$[\mathfrak a(2,3),\mathbf H_\varepsilon]$. Since
$d\mathfrak c_\varepsilon|_{\mathbf H_\varepsilon}$ is surjective, the contribution of
the second term modulo this subspace has dimension
$\operatorname{rank}d(\mathfrak c_\varepsilon\circ\mathbf H_\varepsilon)$.
Consequently,
$\operatorname{rank}d\mathbf J_\varepsilon
=8+\operatorname{rank}d(\mathfrak c_\varepsilon\circ\mathbf H_\varepsilon)$.
Using the expressions of the Casimir invariants in
\eqref{eq:casimir-values}, we find
\[
d(\mathfrak c_\varepsilon\circ\mathbf H_\varepsilon)=
\begin{pmatrix}
1&3\tau&0\\
2\tau^2&4\kappa\tau+8\tau^3&2\tau_1
\end{pmatrix}.
\]
The first row never vanishes, and the matrix has rank one precisely
when $\tau_1=0$ and $\tau(2\kappa+\tau^2)=0$. This establishes the
description of $\mathcal Y_*$ and the asserted ranks of
$d\mathbf J_\varepsilon$.
\end{proof}

\begin{definition}\label{def:phase-curve}
For each $\boldsymbol{\mu}\in\mathbf J_\varepsilon(\mathcal Y)$, let $Z_{\boldsymbol{\mu}}$ denote
its centralizer in $\AWitt_+(2,3)$. The \emph{reduced space} and the
\emph{phase curve} at momentum $\boldsymbol{\mu}$ are, respectively,
\[
\mathcal O_{\boldsymbol{\mu}}^\varepsilon=\mathbf J_\varepsilon^{-1}(\boldsymbol{\mu})/Z_{\boldsymbol{\mu}},
\qquad
\mathcal C_{\boldsymbol{\mu}}^\varepsilon=\mathbf H_\varepsilon(\mathbb R^3)\cap
\operatorname{Ad}\bigl(\AWitt_+(2,3)\bigr)\boldsymbol{\mu}.
\]
Here $Z_{\boldsymbol{\mu}}$ acts freely and properly on the momentum fiber by left
multiplication on the frame component.
\end{definition}

Since the orbit map $\mathbf A\mapsto\mathbf A^{-1}\boldsymbol{\mu}\mathbf A$ admits local smooth sections, it follows from (2.5) that $[(\mathbf A,\kappa,\tau,\tau_1)]\in\mathcal O_{\boldsymbol{\mu}}^\varepsilon\mapsto\mathbf H_\varepsilon(\kappa,\tau,\tau_1)\in\mathcal C_{\boldsymbol{\mu}}^\varepsilon$ is a homeomorphism. Moreover, if $\boldsymbol{\mu}$ is a regular value of $\mathbf J_\varepsilon$, Proposition~\ref{prop:moment-map-regularity} implies that this homeomorphism is a diffeomorphism between smooth one-dimensional manifolds, not necessarily connected.

Let $(h,E)=\mathfrak c_\varepsilon(\boldsymbol{\mu})$. By \eqref{eq:casimir-values}, every element of $\mathcal C_{\boldsymbol{\mu}}^\varepsilon$ can be written uniquely as $\mathbf H_\varepsilon(h-3x^2/2,x,y)$, where $(x,y)$ belongs to the plane quartic
\begin{equation}\label{eq:reduced-plane-quartic}
\mathcal S_{(h,E)}=\left\{(x,y)\in\mathbb R^2:\ y^2=x^4-2hx^2+E\right\}.
\end{equation}
The planar image of $\mathcal C_{\boldsymbol{\mu}}^\varepsilon$ is precisely
\[
\widetilde{\mathcal C}_{\boldsymbol{\mu}}^\varepsilon
=\left\{(x,y)\in\mathcal S_{(h,E)}:\mathbf H_\varepsilon(h-3x^2/2,x,y)\in\operatorname{Ad}(\AWitt_+(2,3))\boldsymbol{\mu}\right\}
\subseteq\mathcal S_{\mathfrak c_\varepsilon(\boldsymbol{\mu})}.
\]
Since all Hamiltonian matrices are regular, $\widetilde{\mathcal C}_{\boldsymbol{\mu}}^\varepsilon$ is a union of connected components of $\mathcal S_{\mathfrak c_\varepsilon(\boldsymbol{\mu})}$. The reduced Reeb field on $\widetilde{\mathcal C}_{\boldsymbol{\mu}}^\varepsilon$ is \mbox{$\overline{\mathcal R}_{\varepsilon,\boldsymbol{\mu}}=y\,\partial_x+2x(x^2-h)\,\partial_y$}. For related reduction procedures, see~\cite{GutierrezSagredo2025} and~\cite{MussoGrant2004}.

\vspace{0.7\baselineskip}

\begin{figure}[H]
\centering
\includegraphics[width=0.315\linewidth]{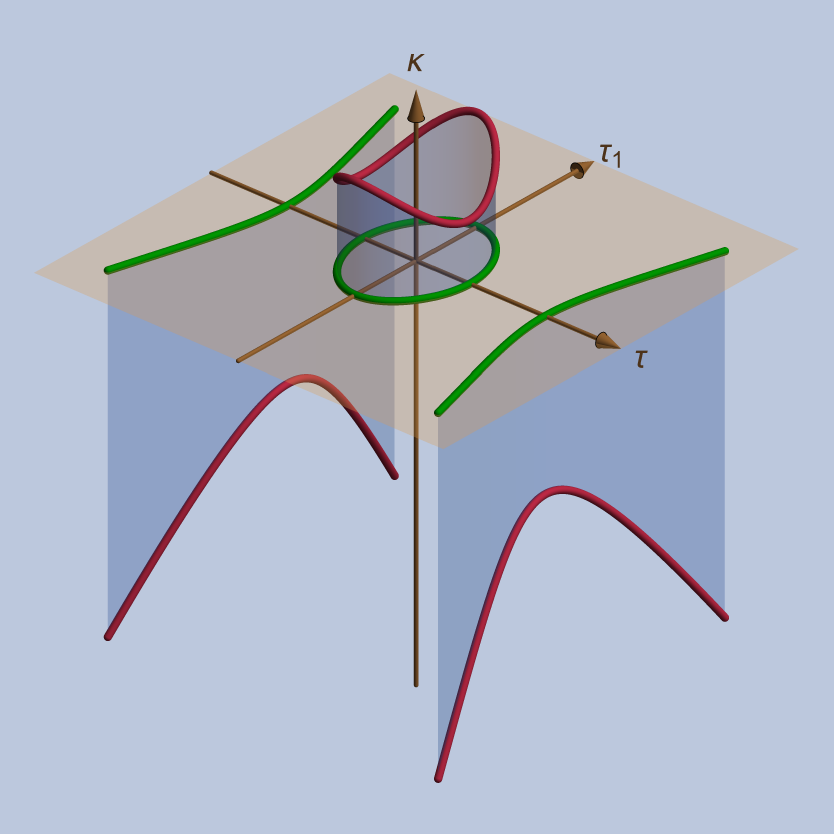}\hfill
\includegraphics[width=0.315\linewidth]{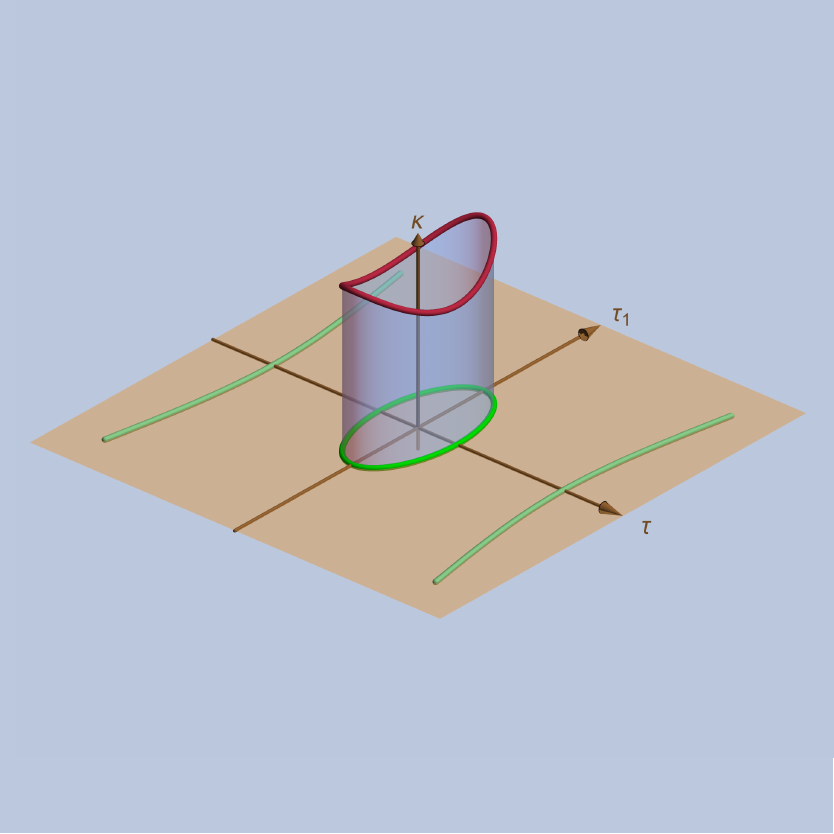}\hfill
\includegraphics[width=0.315\linewidth]{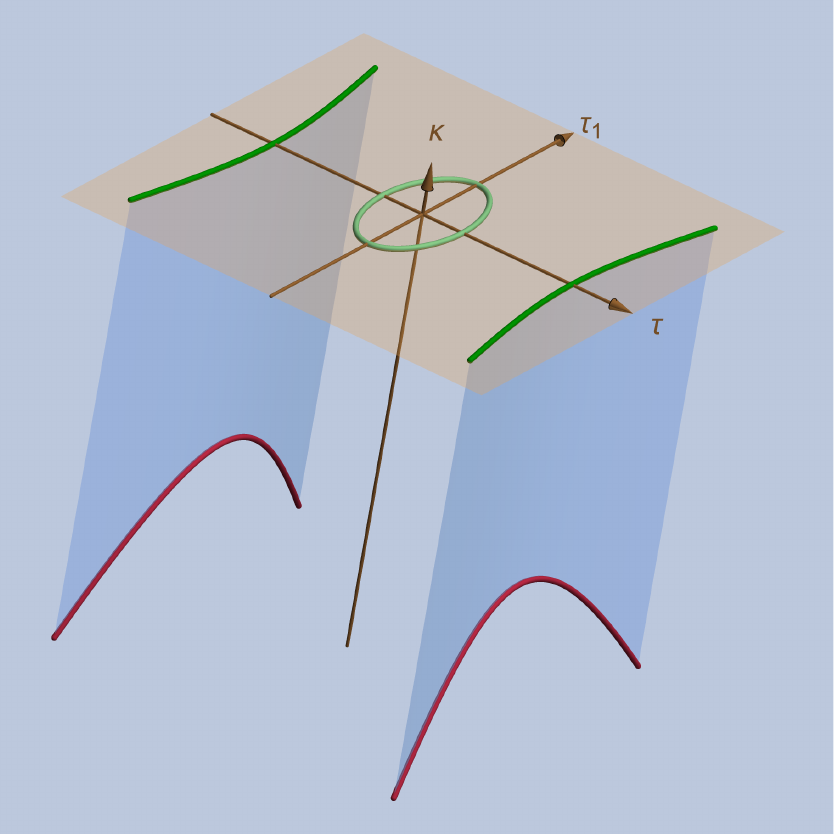}
\caption{Phase curves for negative helicity in coordinates
$(\tau,\tau_1,\kappa)$. Left: $(h,E)=(1,1/2)$ and
$\boldsymbol\mu=\mathbf H_-(1,0,1/\sqrt2)$; the phase curve has
three connected components, one compact and two unbounded.
Middle: $(h,E)=(2,1)$ and $\boldsymbol\mu=\mathbf H_-(2,0,1)$;
only the compact component is selected.
Right: $(h,E)=(2,1)$ and $\boldsymbol\mu=\mathbf H_-(-4,2,1)$;
only the two unbounded components are selected.
The spatial phase curves are shown in burgundy, their projections
onto $\kappa=0$ in dark green, and the unselected components of the
quartic in sage green. The translucent blue surfaces indicate the projections.}
\label{fig:reduced-phase-curves}
\end{figure}

Figure~\ref{fig:reduced-phase-curves} illustrates three possible configurations
of the phase curve for negative helicity when the quartic has three
connected components: the entire quartic, the compact oval alone,
or the two unbounded branches alone.

For regular values of $\mathbf J_\varepsilon$, the reduced vector field
is nowhere vanishing, and each connected component of $\widetilde{\mathcal C}_{\boldsymbol{\mu}}^\varepsilon$
is the image of a maximal reduced trajectory. The corresponding conformal
invariants are recovered from $\tau(s)=x(s)$ and
$\kappa(s)=h-3x(s)^2/2$. In particular, compact components yield periodic
torsion and curvature.

We now consider compact components of the phase curves. A nonconstant
periodic reduced trajectory must lie on a compact nonsingular oval of
the quartic. Such an oval exists precisely when $h>0$ and $0<E<h^2$.
For such $(h,E)$, we define
\begin{equation}\label{eq:elliptic-parameters}
r_\pm(h,E)=\bigl(h\pm\sqrt{h^2-E}\bigr)^{1/2},\qquad m(h,E)=r_-(h,E)^2/r_+(h,E)^2.
\end{equation}
The unique compact oval is given by
\begin{equation}\label{eq:reduced-phase-curve}
\widehat{\mathcal S}_{(h,E)}
=\left\{(x,y)\in\mathbb R^2:
|x|\le r_-,\quad
 y^2=(r_-^2-x^2)(r_+^2-x^2)\right\}.
\end{equation}
For either helicity, $\widehat{\mathcal S}_{(h,E)}$ is a periodic orbit
of the planar reduced dynamics through $(0,\sqrt E)$. The Lax equation then shows that it is a connected component of $\widetilde{\mathcal C}_{\boldsymbol{\mu}}^\varepsilon$, where $\boldsymbol{\mu}=\mathbf H_\varepsilon(h,0,\sqrt E)$.

The equilibria and compact trajectories of the reduced dynamics yield
the following classification of the conformal curvature and torsion of
monodromic conformal extremals.

\begin{proposition}\label{prop:monodromic-classification}
Let $\gamma$ be a monodromic conformal extremal with conformal
curvature $\kappa_\gamma$, torsion $\tau_\gamma$, and parameters
$(h_\gamma,E_\gamma)$. Then exactly one of the following cases occurs.
\begin{enumerate}[label=\textup{(\roman*)},leftmargin=*,itemsep=0pt,parsep=0pt]
\item The conformal torsion vanishes identically. In this case,
$\tau_\gamma=0$, $\kappa_\gamma=h_\gamma$, and $E_\gamma=0$,
where $h_\gamma\in\mathbb R$ is arbitrary.

\item The conformal torsion is a nonzero constant. In this case,
$\tau_\gamma=c$, $\kappa_\gamma=-c^2/2$, $h_\gamma=c^2$, and
$E_\gamma=c^4$, where $c\in\mathbb R\setminus\{0\}$.

\item The conformal torsion is nonconstant. Then $h_\gamma>0$,
$0<E_\gamma<h_\gamma^2$, and, up to a translation of the natural
parameter,
\begin{equation}\label{eq:periodic-torsion-curvature}
\tau_\gamma(s)=r_-\operatorname{sn}(r_+s\mid m),\quad \kappa_\gamma(s)=h_\gamma-\frac32r_-^2\operatorname{sn}^2(r_+s\mid m).
\end{equation}
Here $r_\pm$ and $m$ are defined by~\eqref{eq:elliptic-parameters}
for $(h,E)=(h_\gamma,E_\gamma)$. The least positive period of the torsion is
\begin{equation}\label{eq:fundamental-period}
T_\gamma=4K(m)/r_+,
\end{equation}
where $K(m)$ is the complete elliptic integral of the first kind (cf.~\cite[formula~17.3.1]{AS}).
The curvature has least positive period $T_\gamma/2$.
\end{enumerate}
\end{proposition}

\begin{proof}
For constant torsion, the Euler--Lagrange equations give
$\tau_\gamma(2\kappa_\gamma+\tau_\gamma^2)=0$, yielding (i) and (ii). For nonconstant periodic torsion, the reduced trajectory is a compact
orbit that cannot meet an equilibrium. It therefore traverses the
compact oval $\widehat{\mathcal S}_{(h_\gamma,E_\gamma)}$. The preceding
discussion gives the restrictions on $(h_\gamma,E_\gamma)$. The explicit
formulas and periods follow from the differential equation and periodicity
of the Jacobi elliptic sine; see~\cite[\S\S 22.13 and 22.4]{dlmf}.
Uniqueness of solutions shows that all periodic trajectories on the oval
differ by a translation of the natural parameter.
\end{proof}

\subsection{Closed conformal extremals with constant conformal torsion}\label{subsec:constant-invariants}

We conclude by classifying the closed solutions in the constant-invariant cases (i) and (ii) of Proposition~\ref{prop:monodromic-classification}.

\begin{theorem}\label{thm:closed-constant-invariants}
Let $\gamma:\mathbb R\to\Ein$ be a conformal extremal with constant conformal curvature $\kappa$ and conformal torsion $\tau$. Then $\gamma$ is closed if and only if $\varepsilon=-1$, $\tau=0$, and $\kappa=(m^2+n^2)/(2mn)$, for some relatively prime positive integers $m<n$. In this case, $\gamma$ lies on an anti-de Sitter wall of $\Ein$ (see Definition~\ref{def:chambers-walls}) and is a simple closed curve of homology class $(m,n)$.
\end{theorem}

\begin{proof}
The Euler--Lagrange equations reduce to $\tau(2\kappa+\tau^2)=0$, so that either $\tau=0$, or $\tau\ne0$ and $\kappa=-\tau^2/2$. The Frenet matrix $\mathbf K_\varepsilon:=\mathbf K_\varepsilon(\kappa,\tau)$, defined in~\eqref{eq:frenet-matrix}, is constant. Up to congruence, the canonical frame is $\mathbf A(s)=\exp(s\mathbf K_\varepsilon)$. Since $\gamma$ is closed, its canonical frame is $T$-periodic for some $T>0$. Hence $\exp(T\mathbf K_\varepsilon)=\mathbf I$, which implies that $\mathbf K_\varepsilon$ is semisimple and has purely imaginary eigenvalues, possibly zero. The characteristic polynomial of $\mathbf K_\varepsilon$ is $z\bigl(z^4+(2\kappa-\tau^2)z^2-2\kappa\tau^2-\varepsilon\bigr)$.

When $\tau\ne0$, the quartic factor reduces to $(z^2-\tau^2)^2-\varepsilon$, which is incompatible with a purely imaginary spectrum. Thus $\tau=0$, and the remaining polynomial forces $\varepsilon=-1$ and $\kappa\ge1$. Since $\mathbf K_-$ is not semisimple when $\kappa=1$, we conclude that $\kappa>1$. Set $d=\sqrt{\kappa^2-1}$, $\omega_1=\sqrt{\kappa-d}$, and $\omega_2=\sqrt{\kappa+d}$, with $0<\omega_1<\omega_2$. \Needspace{6\baselineskip}Then $\omega_1\omega_2=1$, and the eigenvalues of $\mathbf K_-$ are $0,\pm i\omega_1,\pm i\omega_2$. Define
\[
B_{2j-1}=\sqrt{\frac d2}\,E_1+\frac{(-1)^j}{\sqrt{2d}}\,E_5,\quad B_{2j}=\frac{(-1)^j(\omega_j^2E_2+E_4)}{\omega_j\sqrt{2d}},\quad B_5=-E_3,\quad j=1,2.
\]
A direct computation shows that $\mathcal B=(B_1,\ldots,B_5)$ is a positively oriented pseudo-orthonormal basis satisfying $A_1(0)=(B_1+B_3)/\sqrt{2d}$ and
\[
\mathbf K_{\mathcal B}
:=\mathbf B^{-1}\mathbf K_-\mathbf B
=\omega_1(\mathsf E^2_1-\mathsf E^1_2)
+\omega_2(\mathsf E^4_3-\mathsf E^3_4),
\]
where $\mathbf B$ denotes the matrix whose columns are the coordinates of $B_1,\ldots,B_5$ with respect to the standard Witt basis of $\mathbb R^{2,3}$.

Thus $\Gamma(s):=\sqrt{2d}\,A_1(s)$ is a null lift of $\gamma$, and
\[
\Gamma(s)=\cos(\omega_1s)B_1+\sin(\omega_1s)B_2+\cos(\omega_2s)B_3+\sin(\omega_2s)B_4.
\]

Hence $\gamma$ is closed if and only if the two frequencies are commensurable. Writing $\omega_1/\omega_2=m/n$, with $0<m<n$ and $\gcd(m,n)=1$, and using $\omega_1\omega_2=1$, we obtain $\omega_1=\sqrt{m/n}$, $\omega_2=\sqrt{n/m}$, and $\kappa=(m^2+n^2)/(2mn)$. Set $\theta=s/\sqrt{mn}$. The corresponding anti-de Sitter wall $\mathcal T_\infty(B_5)$ is parametrized by the null rays $[\cos u\,B_1+\sin u\,B_2+\cos v\,B_3+\sin v\,B_4]$, with $u,v\in\mathbb R/2\pi\mathbb Z$. In the homogeneous coordinates determined by $\mathcal B$, the curve is represented by
\[
\gamma_{m,n}(\theta)
=\bigl[\cos(m\theta),\sin(m\theta),\cos(n\theta),\sin(n\theta),0\bigr].
\]

Since $m$ and $n$ are relatively prime, this parametrization is an embedding of $S^1$ into the anti-de Sitter wall, with homology class $(m,n)$ relative to the homology generators determined by the two oriented coordinate circles.
\end{proof}

\begin{remark}
More generally, every conformal extremal with constant conformal invariants and helicity $\varepsilon$ is, up to conformal congruence, the orbit of the null ray $[E_1]$ under the one-parameter subgroup generated by $\mathbf K_\varepsilon$. Thus explicit parametrizations can be obtained by computing the matrix exponential $\exp(s\mathbf K_\varepsilon)$. They involve trigonometric and hyperbolic functions, possibly multiplied by polynomial factors in the non-semisimple cases. A complete treatment would require a case-by-case analysis of the Jordan normal forms of $\mathbf K_\varepsilon$, an elementary but lengthy procedure that we do not pursue here.
\end{remark}

\subsection{Noncommutative integrability}\label{subsec:noncommutative-integrability}

\begin{definition}\label{def:weak-noncommutative-integrability}
Let $\alpha$ be a contact form on a manifold $M$, and let $R_\alpha$ be its Reeb vector field. We say that $\alpha$ is \emph{weakly noncommutatively integrable} if there exist an open dense subset $M_{\mathrm{reg}}\subset M$ and an $r$-dimensional abelian Lie algebra $\mathfrak a$ of strict infinitesimal symmetries of $\alpha$ such that:
\begin{enumerate}[label=\textup{(\roman*)},leftmargin=*,itemsep=0pt,parsep=0pt]
\item $R_\alpha\in\mathfrak a$ and $\dim\mathfrak a_x=r$ for every $x\in M_{\mathrm{reg}}$, where $\mathfrak a_x=\{X(x):X\in\mathfrak a\}$;
\item the distribution $x\mapsto\mathfrak a_x$ is fibrating on $M_{\mathrm{reg}}$, meaning that it is the vertical distribution of a smooth surjective submersion $\pi_B:M_{\mathrm{reg}}\to B$.
\end{enumerate}
If $\pi_B$ can be chosen proper, we say that $\alpha$ is \emph{noncommutatively integrable}.
\end{definition}

\begin{remark}
The above definition is tailored to the Reeb flow of a fixed contact form: the abelian Lie algebra consists of strict infinitesimal symmetries and contains the Reeb vector field. In the proper case, it is a Reeb-type specialization of the noncommutative contact-integrability framework of~\cite{Jovanovic2012}. The qualifier \emph{weakly} indicates that properness is not required.
\end{remark}

\begin{remark}
Suppose that $\pi_B$ is proper, and let $\mathcal F_{\mathfrak a}$ be the foliation tangent to the distribution $x\mapsto\mathfrak a_x$. The Ehresmann fibration theorem~\cite{Ehresmann1951} implies that the leaf space $M_{\mathrm{reg}}/\mathcal F_{\mathfrak a}$, with its quotient topology, is a smooth manifold and that the natural projection $\pi_{\mathcal F_{\mathfrak a}}:M_{\mathrm{reg}}\to M_{\mathrm{reg}}/\mathcal F_{\mathfrak a}$ is a smooth surjective submersion. Moreover, the induced map $\overline{\pi}_B:M_{\mathrm{reg}}/\mathcal F_{\mathfrak a}\to B$ is a finite covering. Thus, in the proper case, one can replace $B$ by the leaf space $M_{\mathrm{reg}}/\mathcal F_{\mathfrak a}$; the fibration is then intrinsically determined by $\mathfrak a$.
\end{remark}

\Needspace{6\baselineskip}
\begin{theorem}\label{thm:weak-noncommutative-integrability}
For each $\varepsilon\in\{-1,1\}$, the contact form $\alpha_\varepsilon$ on $\mathcal Y$ is weakly noncommutatively integrable.
\end{theorem}

\begin{proof}
Proposition~\ref{prop:moment-map-regularity} identifies $\mathcal Y_*$ as the critical locus of $\mathbf J_\varepsilon$ and shows that $\mathbf J_\varepsilon$ has rank $10$ on $\mathcal Y\setminus\mathcal Y_*$. Set $\mathcal B_\varepsilon:=\mathbf J_\varepsilon(\mathcal Y\setminus\mathcal Y_*)$. Since $\mathbf J_\varepsilon$ has maximal rank on this open dense set, $\mathcal B_\varepsilon$ is open, and the restricted moment map $\mathbf J_\varepsilon|_{\mathcal Y\setminus\mathcal Y_*}:\mathcal Y\setminus\mathcal Y_*\to\mathcal B_\varepsilon$ is a smooth surjective submersion.

Let $\overline h_\varepsilon,\overline E_\varepsilon:\mathcal Y\to\mathbb R$ be the functions defined by
\begin{equation}\label{eq:casimir-pullbacks}
\overline h_\varepsilon=h_\varepsilon\circ\mathbf J_\varepsilon,\qquad
\overline E_\varepsilon=E_\varepsilon\circ\mathbf J_\varepsilon.
\end{equation}
For a first integral $f$ of $R_\varepsilon$, let $X_f$ be its contact Hamiltonian vector field, defined by $\alpha_\varepsilon(X_f)=f$ and $\iota_{X_f}d\alpha_\varepsilon=-df$. For first integrals $f$ and $g$ of $R_\varepsilon$, we use the Jacobi bracket convention $[f,g]=d\alpha_\varepsilon(X_f,X_g)$. The Casimir property gives $R_\varepsilon(\overline h_\varepsilon)=R_\varepsilon(\overline E_\varepsilon)=0$ and $[\overline h_\varepsilon,\overline E_\varepsilon]=0$. Hence $X_{\overline h_\varepsilon}$ and $X_{\overline E_\varepsilon}$ are strict infinitesimal symmetries of $\alpha_\varepsilon$ and commute. Every strict infinitesimal symmetry commutes with $R_\varepsilon$. Therefore
\[
\mathfrak a_\varepsilon
=\operatorname{span}_{\mathbb R}
\left\{R_\varepsilon,X_{\overline h_\varepsilon},X_{\overline E_\varepsilon}\right\}
\]
is an abelian Lie algebra of strict infinitesimal symmetries.

From the definition of the moment map and the Casimir property, we have $d\mathbf J_\varepsilon(R_\varepsilon)=d\mathbf J_\varepsilon(X_{\overline h_\varepsilon})=d\mathbf J_\varepsilon(X_{\overline E_\varepsilon})=0$. On $\mathcal Y\setminus\mathcal Y_*$, the differentials $d\overline h_\varepsilon$ and $d\overline E_\varepsilon$ are independent; consequently, their contact Hamiltonian vector fields are independent modulo $R_\varepsilon$. Since $\dim\mathcal Y=13$ and $\operatorname{rank}d\mathbf J_\varepsilon=10$, it follows that
\[
\mathfrak a_{\varepsilon,x}
=\operatorname{span}
\left\{R_\varepsilon(x),X_{\overline h_\varepsilon}(x),X_{\overline E_\varepsilon}(x)\right\}
=\ker d\!\left(\mathbf J_\varepsilon|_{\mathcal Y\setminus\mathcal Y_*}\right)_x.
\]
Thus the distribution $x\mapsto\mathfrak a_{\varepsilon,x}$ has constant rank $3$ and is the vertical distribution of the restricted moment map. This proves the theorem.
\end{proof}

\begin{definition}\label{def:regular-elliptic}
An element $\boldsymbol{\mu}\in\mathfrak a(2,3)$ is called \emph{regular elliptic} if its spectrum is $\{0,\pm i\omega_1,\pm i\omega_2\}$, where $0<\omega_1<\omega_2$.
\end{definition}

\begin{definition}\label{def:restricted-momentum-spaces}
For each $\sigma\in\{\uparrow,\downarrow\}$, let $\widehat{\mathcal Y}_{-1}^\sigma\subset\mathcal Y\setminus\mathcal Y_*$ consist of the points $y=(\mathbf A,\kappa,\tau,\tau_1)$ such that $\mathbf A\in\mathbb A_+^\sigma(2,3)$, the moment $\boldsymbol{\mu}=\mathbf J_-(y)$ is regular elliptic, the plane quartic $\mathcal S_{\mathfrak c_-(\boldsymbol{\mu})}$ is smooth with a compact oval $\widehat{\mathcal S}_{\mathfrak c_-(\boldsymbol{\mu})}\subset\mathcal S_{\mathfrak c_-(\boldsymbol{\mu})}$, and the phase point $(\tau,\tau_1)$ belongs to $\widehat{\mathcal S}_{\mathfrak c_-(\boldsymbol{\mu})}$. We call $\widehat{\mathcal Y}_{-1}^\uparrow$ and $\widehat{\mathcal Y}_{-1}^\downarrow$ the \emph{restricted momentum spaces}, and denote the restrictions of $\mathbf J_-$ to them by $\mathbf J_-^\uparrow$ and $\mathbf J_-^\downarrow$, respectively; these are the \emph{restricted moment maps}.
\end{definition}

\begin{remark}\label{rem:restricted-momentum-spaces}
The geometric conditions in Definition~\ref{def:restricted-momentum-spaces} are equivalent to
\begin{equation}\label{eq:restricted-momentum-spaces}
\widehat{\mathcal Y}_{-1}^\sigma=\mathbb A_+^\sigma(2,3)\times\mathcal P,\qquad \sigma\in\{\uparrow,\downarrow\},
\end{equation}
where $\mathcal P\subset\mathbb R^3$ is the open, nonempty set of triples satisfying
\[
\kappa+3\tau^2/2>1,\qquad 0<\tau_1^2+2\kappa\tau^2+2\tau^4<(\kappa+3\tau^2/2)^2-1,\qquad 2\kappa+\tau^2>0.
\]
\end{remark}

\Needspace{6\baselineskip}
\begin{theorem}\label{thm:restricted-noncommutative-integrability}
For each $\sigma\in\{\uparrow,\downarrow\}$, the restricted contact form $\left.\alpha_{-1}\right|_{\widehat{\mathcal Y}_{-1}^\sigma}$ is noncommutatively integrable.
\end{theorem}

\begin{proof}
We prove the assertion on $\widehat{\mathcal Y}_{-1}^\uparrow$. In what follows, we use $\mathbb A_+^\uparrow(2,3)$ throughout. Set
\begin{equation}\label{eq:regular-elliptic-domain}
\mathcal D=\left\{(h,E)\in\mathbb R^2:h>1,\;0<E<h^2-1\right\}.
\end{equation}
For $(h,E)\in\mathcal D$, the phase point $(0,\sqrt E)$ belongs to the compact oval $\widehat{\mathcal S}_{(h,E)}$ of the corresponding reduced phase curve. Consequently, the image of the restricted moment map is
\begin{equation}\label{eq:restricted-moment-map-image}
\mathcal U_+:=\mathbf J_-^\uparrow\bigl(\widehat{\mathcal Y}_{-1}^\uparrow\bigr)=\left\{\mathbf A\,\mathbf H_-(h,0,\sqrt E)\,\mathbf A^{-1}:\mathbf A\in\mathbb A_+^\uparrow(2,3),\;(h,E)\in\mathcal D\right\}.
\end{equation}
The Casimir map $\mathfrak c_-:\mathcal U_+\to\mathcal D$ is the map defined in \eqref{eq:casimir-map} on the present domain.
By the preceding theorem, $\mathbf J_-^\uparrow:\widehat{\mathcal Y}_{-1}^\uparrow\to\mathcal U_+$ is a smooth surjective submersion with vertical distribution $x\mapsto\mathfrak a_{-1,x}$. It remains to prove that this map is proper.

Fix $\boldsymbol{\mu}_0\in\mathcal U_+$. Define the smooth functions $\omega,\nu:\mathcal D\to\mathbb R$ by
\begin{equation}\label{eq:normal-form-frequencies}
\omega(h,E)=(h+(h^2-E-1)^{1/2})^{1/2},\qquad \nu(h,E)=(h-(h^2-E-1)^{1/2})^{1/2}.
\end{equation}
On a sufficiently small contractible neighborhood $V\subset\mathcal U_+$ of $\boldsymbol{\mu}_0$, the real spectral subspaces admit smooth pseudo-orthonormal frames. Hence one can choose smooth conjugators $\mathbf C(\boldsymbol{\mu})\in\mathbb A_+^\uparrow(2,3)$ such that
\begin{equation}\label{eq:normal-form-conjugacy}
\mathbf C(\boldsymbol{\mu})^{-1}\boldsymbol{\mu}\mathbf C(\boldsymbol{\mu})=\widehat{\mathbf D}_{\mathfrak c_-(\boldsymbol{\mu})}:=\mathbf P\,\operatorname{diag}\bigl(\omega(\mathfrak c_-(\boldsymbol{\mu}))J_0,\!\nu(\mathfrak c_-(\boldsymbol{\mu}))J_0,\!0\bigr)\,\mathbf P^{-1},\quad J_0=\left(\begin{smallmatrix}0&-1\\1&0\end{smallmatrix}\right).
\end{equation}
Here $\mathbf P$ is the change-of-basis matrix defined in~\eqref{eq:change-of-basis-matrix}. The eigenvalues of this normal form, and hence of $\boldsymbol{\mu}$, are $0,\pm i\omega(\mathfrak c_-(\boldsymbol{\mu})),\pm i\nu(\mathfrak c_-(\boldsymbol{\mu}))$.

The centralizer of $\widehat{\mathbf D}_{\mathfrak c_-(\boldsymbol{\mu})}$ in $\mathbb A_+^\uparrow(2,3)$ is the two-dimensional torus $\widehat K=\mathbf P K\mathbf P^{-1}$, where
\begin{equation}\label{eq:compact-centralizer-K}
K=\left\{\operatorname{diag}(Q_1,Q_2,1):Q_1,Q_2\in\mathrm{SO}(2)\right\}.
\end{equation}
Accordingly, $Z_{\mathbb A_+^\uparrow(2,3)}(\boldsymbol{\mu})=\mathbf C(\boldsymbol{\mu})\widehat K\mathbf C(\boldsymbol{\mu})^{-1}$. For $\boldsymbol{\mu}\in V$, put $r_\pm^\mu=r_\pm(\mathfrak c_-(\boldsymbol{\mu}))$. The compact oval $\widehat{\mathcal S}_{\mathfrak c_-(\boldsymbol{\mu})}$ of the reduced phase curve admits the smooth parametrization
\[
\theta\longmapsto\bigl(r_-^\mu\cos\theta,\ r_-^\mu\sin\theta\bigl((r_+^\mu)^2-(r_-^\mu)^2\cos^2\theta\bigr)^{1/2}\bigr),\qquad\theta\in\mathbb R/2\pi\mathbb Z.
\]
The oval therefore varies smoothly with $\boldsymbol{\mu}$. For each point $(\tau,\tau_1)$ of this oval, with $\kappa=h_-(\boldsymbol{\mu})-\frac32\tau^2$, the centralizer $Z_{\mathbb A_+^\uparrow(2,3)}(\boldsymbol{\mu})$ acts freely and transitively, by left multiplication, on the set $\{\mathbf A\in\mathbb A_+^\uparrow(2,3):\mathbf A\mathbf H_-(\kappa,\tau,\tau_1)\mathbf A^{-1}=\boldsymbol{\mu}\}$. Using the smooth local conjugators $\mathbf C(\boldsymbol{\mu})$, these sets form a principal $\widehat K$-bundle over $V\times S^1$. Since $V$ is contractible and $\widehat K\cong\mathbb T^2$ is connected, this bundle is trivial. It follows that $(\mathbf J_-^\uparrow)^{-1}(V)\cong V\times S^1\times\mathbb T^2$. In particular, the restricted moment map is locally trivial with compact fibre $\mathbb T^3$, and hence is proper.

The same conclusion holds on $\widehat{\mathcal Y}_{-1}^\downarrow$. Indeed, the symmetry $\Phi_{\mathbf S}(\mathbf A,\kappa,\tau,\tau_1)=(\mathbf S\mathbf A,\kappa,\tau,\tau_1)$, for a fixed $\mathbf S\in\mathbb A_+^\downarrow(2,3)$, carries the upward restricted momentum space to the downward one, preserves $\alpha_{-1}$, and satisfies $\mathbf J_-\circ\Phi_{\mathbf S}=\operatorname{Ad}_{\mathbf S}\circ\mathbf J_-$. Thus properness transfers to the downward branch.
\end{proof}

\begin{remark}
The fibres of each restricted moment map are connected, since they are torus bundles over compact ovals with connected structure group $\mathbb T^2$. They are therefore the leaves of the distribution $x\mapsto\mathfrak a_{-1,x}$, and $\mathcal U_+$ identifies naturally with the leaf space on the upward restricted momentum space. Theorems~\ref{thm:invariant-torus-fibration} and~\ref{thm:linearization-reeb} will give an explicit trivializing cover of this restricted moment map by three open subsets of $\mathcal U_+$, together with local coordinates that linearize the Reeb flow. The corresponding construction on the downward branch is obtained through $\Phi_{\mathbf S}$.
\end{remark}

\section{Closed conformal extremals with nonconstant conformal torsion}\label{sec:integration}

Throughout this section, we consider only conformal extremals with
nonconstant periodic conformal torsion. This assumption will be understood without further mention.

\subsection{Spectral restrictions on the momenta of closed conformal extremals}\label{subsec:spectral-restrictions}

\begin{proposition}\label{prop:spectral-restrictions}
Let $\gamma$ be a conformal extremal with nonconstant periodic conformal torsion.
If $\gamma$ is closed, then its helicity is negative, its parameters satisfy
$h_\gamma>1$ and $0<E_\gamma<h_\gamma^2-1$, and its momentum
$\mathbf J_\gamma$ is regular elliptic.
\end{proposition}

\begin{proof}
Throughout the proof, we use the abbreviated notation
$\kappa=\kappa_\gamma$, $\tau=\tau_\gamma$,
$\mathbf A=\mathbf A_\gamma$, $h=h_\gamma$, and $E=E_\gamma$, and let $\varepsilon\in\{-1,1\}$ denote the helicity of $\gamma$.
By Proposition~\ref{prop:monodromic-classification}, the reduced phase trajectory
$s\mapsto(\tau(s),\tau'(s))$ parametrizes the oval $\widehat{\mathcal S}_{(h,E)}$.
In particular, $h>0$, $0<E<h^2$, and
$\tau(s)^2\le r_-(h,E)^2=h-\sqrt{h^2-E}<h$.

Let $T$ be the fundamental period of $\tau$.
Up to a conformal congruence, we may assume $\mathbf A(0)=\mathbf I$,
so that the monodromy of $\gamma$ is $\mathbf M_\gamma=\mathbf A(T)$.
Since $\gamma$ is closed, the monodromy has finite order.
In particular, $\mathbf M_\gamma$ is semisimple and all its eigenvalues have modulus one.

The characteristic polynomial of the momentum $\mathbf J_\gamma$ is $zp_\varepsilon(z)$,
where $p_\varepsilon(z)=z^4-2\varepsilon hz^2+E-\varepsilon$.

To relate the spectrum of the momentum
$\mathbf J_\gamma=\mathbf A\mathbf H_\gamma\mathbf A^{-1}$
to the monodromy, we construct an eigenvector of $\mathbf H_\gamma$ in canonical coordinates.
Let $z\ne0$ be a root of $p_\varepsilon$ and suppose that
$D_z(s):=z^2-\varepsilon\tau(s)^2\ne0$ for every $s$.
Let $V_z$ be the $T$-periodic nowhere-zero complex vector-valued
function given by
\[
\resizebox{\linewidth}{!}{$\displaystyle
V_z=\bigl((\kappa+\tau^2-\varepsilon z^2)D_z,\,-z-\varepsilon\tau\tau',\,a_\varepsilon zD_z-b_\varepsilon(\tau+z\tau'),\,a_\varepsilon(\tau+z\tau')+b_\varepsilon zD_z,\,D_z\bigr)^t.
$}
\]
A direct computation gives
\begin{equation}\label{eq:Vz}
\mathbf H_\gamma V_z=zV_z,\qquad
V_z'+\mathbf K_\gamma V_z=\rho_zV_z,\qquad
\rho_z=\frac{z}{D_z}+\frac12\frac{D_z'}{D_z}.
\end{equation}
where $\mathbf K_\gamma=\mathbf A^{-1}\mathbf A'$. Let $Q_z$ be the solution of
\begin{equation}\label{eq:Qz}
Q_z'+\rho_zQ_z=0,\qquad Q_z(0)=1.
\end{equation}
Then $(Q_zV_z)'+\mathbf K_\gamma(Q_zV_z)=0$, and hence
$\mathbf A(s)Q_z(s)V_z(s)=V_z(0)$.
Since $V_z$ is $T$-periodic,
$\mathbf M_\gamma V_z(T)=Q_z(T)^{-1}V_z(T)$.
Thus $V_z(T)$ is an eigenvector of $\mathbf M_\gamma$, with eigenvalue
$\chi_z=Q_z(T)^{-1}$.

Since $D_z$ is $T$-periodic and nowhere zero,
equations \eqref{eq:Vz} and \eqref{eq:Qz} yield
\begin{equation}\label{eq:chi}
\log|\chi_z|=\operatorname{Re}\biggl(z\int_0^T\frac{ds}{D_z(s)}\biggr).
\end{equation}

Suppose first that $\varepsilon=1$.
The polynomial $p_+$ has the positive real root
$z_+=(h+\sqrt{h^2-E+1})^{1/2}$.
Since $z_+^2>h$, while $\tau(s)^2\le r_-(h,E)^2<h$,
we have $z_+^2-\tau(s)^2>0$ for every $s$.
Thus \eqref{eq:chi} yields $\log|\chi_{z_+}|>0$,
contradicting the fact that every eigenvalue of $\mathbf M_\gamma$ has modulus one.
Therefore a closed conformal extremal must have negative helicity.

Assume now that $\varepsilon=-1$.
The characteristic polynomial of the momentum is
$\det(z\mathbf I-\mathbf J_\gamma)=z(z^4+2hz^2+E+1)$.
Set $\Delta=h^2-E-1$.
We first exclude the cases $\Delta<0$ and $\Delta=0$.

Suppose that $\Delta<0$.
The polynomial $p_-$ has a root $z=a+ib$ with $ab\ne0$.
Since $z^2+\tau(s)^2\ne0$ for every $s$, formula \eqref{eq:chi} applies.
Moreover,
\[
\operatorname{Re}\frac{z}{z^2+\tau^2}
=\frac{a(|z|^2+\tau^2)}{|z^2+\tau^2|^2},
\]
which has constant nonzero sign. Consequently,
\[
\log|\chi_z|
=\int_0^T\frac{a(|z|^2+\tau^2)}{|z^2+\tau^2|^2}\,ds\ne0.
\]
Thus $|\chi_z|\ne1$, again a contradiction.

Suppose next that $\Delta=0$.
In this case $E=h^2-1$, with $h>1$, so that $p_-(z)=(z^2+h)^2$.
Set $\lambda=i\sqrt h$.
Since $\tau^2\le h-1$, we have $D_\lambda=\tau^2-h\le-1$.
Thus the expressions defining $V_z$ and $\rho_z$ extend holomorphically
to a complex neighborhood of $\lambda$.
For $\varepsilon=-1$, their expressions satisfy
\begin{equation}\label{eq:Vz-polynomial-identities}
\resizebox{.88\linewidth}{!}{$\displaystyle
(\mathbf H_\gamma-z\mathbf I)V_z=-zp_-(z)E_1,\qquad D_z\bigl(V_z'+\mathbf K_\gamma V_z-\rho_zV_z\bigr)=p_-(z)(z^2E_2+\tau zE_3),
$}
\end{equation}
for every complex $z$ such that $D_z\ne0$.
These identities follow by direct substitution, using the Euler--Lagrange equations
\eqref{eq:reeb-dynamics} and the conservation laws \eqref{eq:first-integrals}.

Since $\lambda$ is a double root of $p_-$, we have
$p_-(\lambda)=p_-'(\lambda)=0$.
Define $W_\lambda:=\left.\partial_zV_z\right|_{z=\lambda}$ and $g_\lambda:=\left.\partial_z\rho_z\right|_{z=\lambda}$. Explicitly,
\[
W_\lambda^t=\bigl(\lambda(\tau^2-2h),-1,-\tau',\tau^2-3h,2\lambda\bigr),\quad
g_\lambda=\frac{h+\tau^2}{(h-\tau^2)^2}-\lambda\biggl(\frac{1}{h-\tau^2}\biggr)'.
\]
Differentiating \eqref{eq:Vz-polynomial-identities} with respect to $z$
and evaluating at $z=\lambda$, we obtain
\begin{equation}\label{eq:Wlambda}
(\mathbf H_\gamma-\lambda\mathbf I)W_\lambda=V_\lambda,\qquad
W_\lambda'+\mathbf K_\gamma W_\lambda=\rho_\lambda W_\lambda+g_\lambda V_\lambda.
\end{equation}

Let $Q_\lambda$ be the solution of \eqref{eq:Qz} with $z=\lambda$, and set
\[
F_\lambda(s)=\int_0^s g_\lambda(u)\,du=\int_0^s\frac{h+\tau(u)^2}{(h-\tau(u)^2)^2}\,du-\lambda\biggl(\frac{1}{h-\tau(s)^2}-\frac{1}{h-\tau(0)^2}\biggr).
\]
Since $(\mathbf H_\gamma-\lambda\mathbf I)W_\lambda=V_\lambda\ne0$,
the vectors $V_\lambda$ and $W_\lambda$ are linearly independent.
It follows from \eqref{eq:Vz} and \eqref{eq:Wlambda}
that both $Q_\lambda V_\lambda$ and
$Q_\lambda(W_\lambda-F_\lambda V_\lambda)$ are solutions of the equation $Y'+\mathbf K_\gamma Y=0$.
Using the $T$-periodicity of $V_\lambda$ and $W_\lambda$,
we may therefore write
\[
\mathbf A(s)Q_\lambda(s)V_\lambda(s)=V_\lambda(T),\quad
\mathbf A(s)Q_\lambda(s)\bigl(W_\lambda(s)-F_\lambda(s)V_\lambda(s)\bigr)=W_\lambda(T).
\]
Evaluating at $s=T$ and putting $\chi_\lambda=Q_\lambda(T)^{-1}$,
we obtain
\begin{equation}\label{eq:monodromy-jordan}
\mathbf M_\gamma V_\lambda(T)=\chi_\lambda V_\lambda(T),\quad
\mathbf M_\gamma W_\lambda(T)=\chi_\lambda W_\lambda(T)+\chi_\lambda F_\lambda(T)V_\lambda(T).
\end{equation}
By the $T$-periodicity of $\tau$,
\[
F_\lambda(T)=\int_0^T\frac{h+\tau^2}{(h-\tau^2)^2}\,ds>0.
\]
Since $\chi_\lambda\ne0$, $F_\lambda(T)>0$, and $V_\lambda(T)$ and
$W_\lambda(T)$ are linearly independent, \eqref{eq:monodromy-jordan}
shows that the restriction of $\mathbf M_\gamma$ to their span has
a nontrivial Jordan block.
Hence $\mathbf M_\gamma$ is not semisimple,
contradicting the finite-order assumption.

It follows that $0<E<h^2-1$, and hence $h>1$.
Since $\Delta=h^2-E-1>0$, we have $0<\sqrt\Delta<h$.
Consequently, the four nonzero eigenvalues of $\mathbf J_\gamma$ are
$\pm i(h+\sqrt\Delta)^{1/2}$ and $\pm i(h-\sqrt\Delta)^{1/2}$.
These eigenvalues are distinct and purely imaginary.
Since zero is also simple, $\mathbf J_\gamma$ is regular elliptic.
\end{proof}

\subsection{The invariant torus fibration}\label{subsec:invariant-torus-fibration}

By Proposition 3.1, every closed conformal extremal $\gamma$ with nonconstant periodic torsion has negative helicity and momentum parameters in $\mathcal D$. We equip $\gamma$ with its canonical positive orientation (see Definition~\ref{def:positive-orientation}) and parametrize it by its natural parameter $s$. By Theorem 2.3, its extended canonical frame is an integral curve of the Reeb vector field on the restricted momentum space $\widehat{\mathcal Y}_{-1}^\uparrow$ (see Definition 2.9). Thus, the study of these extremals leads to the study of the restricted Reeb flow. Subsection 2.4 shows that its moment map is a proper fibration with fibres $\mathbb T^3$. We now refine this result by proving that the fibration is a nontrivial principal $\mathbb T^3$-bundle and constructing explicit trivializations over three open subsets of $\mathcal U_+$. These trivializations will be used in Subsection 3.3 to analyze the Reeb flow on the invariant tori.

\begin{theorem}\label{thm:invariant-torus-fibration}
The restricted moment map $\mathbf J_-^\uparrow:\widehat{\mathcal Y}_{-1}^\uparrow\to\mathcal U_+$ is a nontrivial smooth principal $\mathbb T^3$-bundle, whose fibres are invariant under the Reeb flow. The bundle admits explicit smooth trivializations over a cover of $\mathcal U_+$ by three open sets.
\end{theorem}

\begin{proof}
Let $\mathcal P\subset\mathbb R^3$ be the open set introduced in Remark 2.5, and define the polynomial map
\begin{equation}\label{eq:torus-polynomial-p}
 p:\mathcal P\to\mathcal D,\qquad
 (\kappa,\tau,\tau_1)\mapsto
 \bigl(\kappa+\frac32\tau^2,\,
 \tau_1^2+2\kappa\tau^2+2\tau^4\bigr).
\end{equation}
For each $(h,E)\in\mathcal D$, set $\mathbf X_{h,E}:=\mathbf H_-(h,0,\sqrt E)$. Then $\mathfrak c_-(\mathbf X_{h,E})=(h,E)$ by \eqref{eq:casimir-map}.

We first construct a common spectral normal form $\mathbf D_{h,E}$, whose centralizer is a fixed torus $K\simeq\mathbb T^2$. Using the distinguished momenta, we obtain local spectral normalizations $\mathbf B_\sigma$ on a three-set cover of $\mathcal U_+$ and a global normalization $\mathbf Q$ on $\mathcal P$. Their comparison, together with the periodic reduced ovals, yields the principal $\mathbb T^3$-bundle structure. Nontriviality will follow from the hairy ball theorem.

Let $\mathcal F_+^\uparrow(2,3)$ denote the manifold of positively oriented and time-oriented pseudo-orthonormal frames of $\mathbb R^{2,3}$, represented as matrices in the standard Witt basis. Define the smooth function $d:\mathcal D\to(0,\infty)$ by $d(h,E)=\omega(h,E)^2-\nu(h,E)^2$, where $\omega,\nu$ are the functions defined in (2.21). For each $(h,E)\in\mathcal D$, abbreviate $\omega(h,E)$, $\nu(h,E)$, and $d(h,E)$ to $\omega,\nu,d$, respectively.\par\Needspace{6\baselineskip}\noindent For $j=1,\ldots,5$, define the smooth maps $B_j^*:\mathcal D\to\mathbb R^{2,3}$ by the formulas below, writing $B_j^*=B_j^*(h,E)$:
\begin{equation}\label{eq:torus-star-vectors}
\left\{\begin{aligned}
B_1^*&=\frac{\sqrt d}{2}E_1-\frac{E_5}{\sqrt d},
&\qquad B_3^*&=\frac{\sqrt d}{2}E_1+\frac{E_5}{\sqrt d},\\
B_2^*&=\frac{E_2+\sqrt E\,E_3+\omega^2E_4}{\omega\sqrt d},
&\qquad B_4^*&=-\frac{E_2+\sqrt E\,E_3+\nu^2E_4}{\nu\sqrt d},\\
B_5^*&=\frac{\sqrt E\,E_2-E_3}{\nu\omega}.
\end{aligned}\right.
\end{equation}
For each $(h,E)\in\mathcal D$, these vectors form a positively oriented and time-oriented pseudo-orthonormal basis. Let $\mathbf B_*:\mathcal D\to\mathcal F_+^\uparrow(2,3)$ be the smooth map defined by $\mathbf B_*(h,E)=(B_1^*(h,E),\ldots,B_5^*(h,E))$. It satisfies
\begin{equation}\label{eq:torus-distinguished-normalization}
\mathbf B_*(h,E)^{-1}\mathbf X_{h,E}\mathbf B_*(h,E)
=\mathbf D_{h,E}:=\operatorname{diag}(\omega J_0,\nu J_0,0).
\end{equation}
Here $J_0$ is as in~\eqref{eq:normal-form-conjugacy}. Since $0<\nu<\omega$, the centralizer of $\mathbf D_{h,E}$ in $O_+^\uparrow(2,3)$ is precisely the two-dimensional torus $K$ defined in~\eqref{eq:compact-centralizer-K}. Its elements can be parametrized by
\[
R_1(\phi_1,\phi_2)=\operatorname{diag}\bigl(R(\phi_1),R(\phi_2),1\bigr),\qquad
R(\phi)=\begin{pmatrix}\cos\phi&-\sin\phi\\\sin\phi&\cos\phi\end{pmatrix}.
\]
We next construct smooth spectral normalizations on an open cover of $\mathcal U_+$. Define the smooth functions $\boldsymbol{\omega},\boldsymbol{\nu}:\mathcal U_+\to(0,\infty)$ by $\boldsymbol{\omega}(\mathbf X)=\omega(\mathfrak c_-(\mathbf X))$ and $\boldsymbol{\nu}(\mathbf X)=\nu(\mathfrak c_-(\mathbf X))$. For each $\mathbf X\in\mathcal U_+$, its eigenvalues are $0$, $\pm i\boldsymbol{\omega}(\mathbf X)$, and $\pm i\boldsymbol{\nu}(\mathbf X)$, and the corresponding real spectral projectors of $\mathbf X$ are
\[
\left\{\begin{aligned}
P_\omega(\mathbf X)
&=\frac{\mathbf X^2\bigl(\mathbf X^2+\boldsymbol{\nu}(\mathbf X)^2I\bigr)}
{\boldsymbol{\omega}(\mathbf X)^2\bigl(\boldsymbol{\omega}(\mathbf X)^2-\boldsymbol{\nu}(\mathbf X)^2\bigr)},\\[9pt]
P_\nu(\mathbf X)
&=-\frac{\mathbf X^2\bigl(\mathbf X^2+\boldsymbol{\omega}(\mathbf X)^2I\bigr)}
{\boldsymbol{\nu}(\mathbf X)^2\bigl(\boldsymbol{\omega}(\mathbf X)^2-\boldsymbol{\nu}(\mathbf X)^2\bigr)},\\[9pt]
P_0(\mathbf X)
&=\frac{\bigl(\mathbf X^2+\boldsymbol{\nu}(\mathbf X)^2I\bigr)
\bigl(\mathbf X^2+\boldsymbol{\omega}(\mathbf X)^2I\bigr)}
{\boldsymbol{\nu}(\mathbf X)^2\boldsymbol{\omega}(\mathbf X)^2}.
\end{aligned}\right.
\]
Their images give an orthogonal decomposition into a negative-definite plane, a positive-definite plane, and a positive line. For $j=1,3,4,5$, define
\[
\resizebox{\textwidth}{!}{$\displaystyle
f_j(\mathbf X)=\left\langle P_0(\mathbf X)B_j^*(\mathfrak c_-(\mathbf X)),P_0(\mathbf X)B_j^*(\mathfrak c_-(\mathbf X))\right\rangle,\quad
 g_j(\mathbf X)=\left\langle P_\nu(\mathbf X)B_j^*(\mathfrak c_-(\mathbf X)),P_\nu(\mathbf X)B_j^*(\mathfrak c_-(\mathbf X))\right\rangle.
$}
\]
Choose the index pairs $(i_1,j_1)=(4,3)$, $(i_2,j_2)=(5,4)$, and $(i_3,j_3)=(3,5)$, and set
\begin{equation}\label{eq:torus-open-cover}
\mathcal V_\sigma^+=\{\mathbf X\in\mathcal U_+:f_{j_\sigma}(\mathbf X)>0,\ g_{i_\sigma}(\mathbf X)>0\},\qquad \sigma=1,2,3.
\end{equation}

These three open sets cover $\mathcal U_+$. Indeed, fix $\mathbf X\in\mathcal U_+$. At least one of $f_3,f_4,f_5$ is positive, since the positive spectral line cannot be orthogonal to the positive three-plane spanned by $B_3^*(\mathfrak c_-(\mathbf X)),B_4^*(\mathfrak c_-(\mathbf X)),B_5^*(\mathfrak c_-(\mathbf X))$. On the other hand, $\ker P_\nu(\mathbf X)=\operatorname{im}P_\omega(\mathbf X)\oplus\operatorname{im}P_0(\mathbf X)$ has positive index one. Consequently, at most one of $g_3,g_4,g_5$ can vanish. Moreover, if $g_j=0$, then $B_j^*(\mathfrak c_-(\mathbf X))\in\ker P_\nu(\mathbf X)$, and the positivity of $B_j^*(\mathfrak c_-(\mathbf X))$ implies $P_0(\mathbf X)B_j^*(\mathfrak c_-(\mathbf X))\ne0$, hence $f_j>0$. Suppose first that $f_3>0$. If $g_4>0$, then $\mathbf X\in\mathcal V_1^+$. Otherwise, $g_4=0$, which implies $f_4>0$ and $g_5>0$, hence $\mathbf X\in\mathcal V_2^+$. The same argument applies cyclically when $f_4>0$ or $f_5>0$. Therefore, $\mathcal U_+=\mathcal V_1^+\cup\mathcal V_2^+\cup\mathcal V_3^+$.

For each $\sigma=1,2,3$, let $Z_\omega,Z_\nu^\sigma,Z_0^\sigma:\mathcal V_\sigma^+\to\mathbb R^{2,3}$ be the smooth maps defined by
\[
\left\{\begin{aligned}
Z_\omega(\mathbf X)&=\frac{P_\omega(\mathbf X)B_1^*(\mathfrak c_-(\mathbf X))}{\sqrt{1+f_1(\mathbf X)+g_1(\mathbf X)}},\\[8pt]
Z_\nu^\sigma(\mathbf X)&=\frac{P_\nu(\mathbf X)B_{i_\sigma}^*(\mathfrak c_-(\mathbf X))}{\sqrt{g_{i_\sigma}(\mathbf X)}},\\[8pt]
Z_0^\sigma(\mathbf X)&=\frac{P_0(\mathbf X)B_{j_\sigma}^*(\mathfrak c_-(\mathbf X))}{\sqrt{f_{j_\sigma}(\mathbf X)}}.
\end{aligned}\right.
\]
These maps are well defined. Indeed, $f_1,g_1\geq0$, so $1+f_1(\mathbf X)+g_1(\mathbf X)>0$, while the remaining denominators are strictly positive by the definition of $\mathcal V_\sigma^+$. Moreover, the orthogonality of the spectral decomposition and $\langle B_1^*(\mathfrak c_-(\mathbf X)),B_1^*(\mathfrak c_-(\mathbf X))\rangle=-1$ imply
\[
-\langle P_\omega(\mathbf X)B_1^*(\mathfrak c_-(\mathbf X)),P_\omega(\mathbf X)B_1^*(\mathfrak c_-(\mathbf X))\rangle=1+f_1(\mathbf X)+g_1(\mathbf X).
\]
Consequently, $\langle Z_\omega,Z_\omega\rangle=-1$ and $\langle Z_\nu^\sigma,Z_\nu^\sigma\rangle=\langle Z_0^\sigma,Z_0^\sigma\rangle=1$.

Define
\[
\mathbf B_\sigma(\mathbf X)=\bigl(
Z_\omega(\mathbf X),
\frac{\mathbf XZ_\omega(\mathbf X)}{\boldsymbol{\omega}(\mathbf X)},
Z_\nu^\sigma(\mathbf X),
\frac{\mathbf XZ_\nu^\sigma(\mathbf X)}{\boldsymbol{\nu}(\mathbf X)},
\epsilon_\sigma(\mathbf X)Z_0^\sigma(\mathbf X)\bigr),
\]
choosing $\epsilon_\sigma(\mathbf X)\in\{-1,1\}$ to ensure positive orientation. The spectral identities show that the columns form a pseudo-orthonormal basis. The negative spectral plane carries the orientation induced by $\mathbf X/\boldsymbol{\omega}(\mathbf X)$, which agrees with the prescribed time-orientation because $\mathbf X$ is conjugate to $\mathbf X_{h,E}$, where $(h,E)=\mathfrak c_-(\mathbf X)$, by an element of $\mathbb A_+^\uparrow(2,3)$. Since $\epsilon_\sigma$ is locally constant, each map $\mathbf B_\sigma:\mathcal V_\sigma^+\to\mathcal F_+^\uparrow(2,3)$ is smooth and satisfies
\begin{equation}\label{eq:torus-local-spectral-identity}
\mathbf B_\sigma(\mathbf X)^{-1}\mathbf X\mathbf B_\sigma(\mathbf X)=\mathbf D_{\mathfrak c_-(\mathbf X)}.
\end{equation}

We now construct a global smooth spectral normalization for the map $\mathbf H_-:\mathcal P\to\mathfrak a(2,3)$. For each $\mathtt u\in\mathcal P$, we write $\omega(\mathtt u)$, $\nu(\mathtt u)$, and $d(\mathtt u)$ for $\omega(p(\mathtt u))$, $\nu(p(\mathtt u))$, and $\omega(p(\mathtt u))^2-\nu(p(\mathtt u))^2$, respectively, where $p$ is defined by (3.7). We seek a smooth map $\mathbf Q:\mathcal P\to\mathcal F_+^\uparrow(2,3)$ satisfying
\begin{equation}\label{eq:torus-global-spectral-identity}
\mathbf Q(\mathtt u)^{-1}\mathbf H_-(\mathtt u)\mathbf Q(\mathtt u)=\mathbf D_{p(\mathtt u)}.
\end{equation}
For $\mathtt u=(\kappa,\tau,\tau_1)\in\mathcal P$, set $(h,E)=p(\mathtt u)$ and $\delta_a=\sqrt{a^2-\tau^2}$ for $a=\omega,\nu$. Let $W_a,V_a:\mathcal P\to\mathbb R^{2,3}$, $a=\omega,\nu$, and $W_0:\mathcal P\to\mathbb R^{2,3}$ be defined by
\begin{equation}\label{eq:torus-global-vectors}
\left\{\begin{aligned}
W_a&=\frac{1}{a\sqrt d}\bigl[\delta_a\bigl(a^2-h+\tau^2/2\bigr)E_1+\frac{\tau\tau_1}{\delta_a}E_2-\frac{\tau}{\delta_a}E_3-\delta_a E_5\bigr],\\
V_a&=\frac{1}{\sqrt d}\bigl[\frac{E_2+\tau_1E_3}{\delta_a}+\delta_a E_4\bigr],\\
W_0&=\frac{1}{\nu\omega}\bigl[\bigl(h-\tau^2/2\bigr)\tau E_1+\tau_1E_2-E_3+\tau E_5\bigr].
\end{aligned}\right.
\end{equation}
By \eqref{eq:reduced-plane-quartic} and \eqref{eq:elliptic-parameters}--\eqref{eq:reduced-phase-curve}, the fibres of $p$ are naturally identified with the compact ovals $\widehat{\mathcal S}_{(h,E)}$. Indeed, the condition $2\kappa+\tau^2>0$ is equivalent to $\tau^2<h$ and excludes the noncompact branches of the quartic \eqref{eq:reduced-plane-quartic}. Moreover, $\tau^2\le r_-(h,E)^2<\nu^2<\omega^2$, so the maps above are smooth on $\mathcal P$.

Consider the smooth map $\mathbf Q:\mathcal P\to\operatorname{Mat}_{5\times5}(\mathbb R)$ defined by
\begin{equation}\label{eq:torus-global-normalization}
\mathbf Q(\mathtt u)=\bigl(W_\omega(\mathtt u),V_\omega(\mathtt u),-W_\nu(\mathtt u),-V_\nu(\mathtt u),W_0(\mathtt u)\bigr).
\end{equation}
A direct computation shows that $\mathbf Q(\mathtt u)$ has pseudo-orthonormal columns and satisfies the normalization identity \eqref{eq:torus-global-spectral-identity}. Moreover, $\mathbf Q(h,0,\sqrt E)=\mathbf B_*(h,E)$ for every $(h,E)\in\mathcal D$. Since each fibre of $p$ is connected, positive orientation and time-orientation are preserved throughout $\mathcal P$. Consequently, $\mathbf Q$ takes values in $\mathcal F_+^\uparrow(2,3)$ and provides the required global smooth spectral normalization.

We now compare the two spectral normalizations. For $(\mathbf A,\mathtt u)\in(\mathbf J_-^\uparrow)^{-1}(\mathcal V_\sigma^+)$, the identity $p(\mathtt u)=\mathfrak c_-\bigl(\mathbf J_-^\uparrow(\mathbf A,\mathtt u)\bigr)$, together with the definition of the moment map (cf. (2.5)) and the normalization identities \eqref{eq:torus-local-spectral-identity} and \eqref{eq:torus-global-spectral-identity}, shows that $C_\sigma(\mathbf A,\mathtt u):=\mathbf B_\sigma\bigl(\mathbf J_-^\uparrow(\mathbf A,\mathtt u)\bigr)^{-1}\mathbf A\mathbf Q(\mathtt u)$ commutes with $\mathbf D_{p(\mathtt u)}$. Furthermore, both spectral frames have the same orientation and time-orientation, so $C_\sigma(\mathbf A,\mathtt u)\in O_+^\uparrow(2,3)$. Hence $C_\sigma(\mathbf A,\mathtt u)\in K$.

Consider the fibre product $\mathcal V_\sigma^+\times_{\mathcal D}\mathcal P=\{(\mathbf X,\mathtt u)\in\mathcal V_\sigma^+\times\mathcal P:\mathfrak c_-(\mathbf X)=p(\mathtt u)\}$ of $\mathfrak c_-$ and $p$. Then
\begingroup\fontsize{10.5}{12.5}\selectfont
\[
\Xi_\sigma\!:\,(\mathbf A,\mathtt u)\in(\mathbf J_-^\uparrow)^{-1}(\mathcal V_\sigma^+)
\longmapsto\bigl(\mathbf J_-^\uparrow(\mathbf A,\mathtt u),\mathtt u,C_\sigma(\mathbf A,\mathtt u)\bigr)
\in\bigl(\mathcal V_\sigma^+\times_{\mathcal D}\mathcal P\bigr)\times K
\]
\endgroup
is a smooth diffeomorphism with inverse $\Xi_\sigma^{-1}(\mathbf X,\mathtt u,C)=\bigl(\mathbf B_\sigma(\mathbf X)C\mathbf Q(\mathtt u)^{-1},\mathtt u\bigr)$.

We next identify the fibre product $\mathcal V_\sigma^+\times_{\mathcal D}\mathcal P$ with $\mathcal V_\sigma^+\times(\mathbb R/2\pi\mathbb Z)$. By \eqref{eq:elliptic-parameters}--\eqref{eq:fundamental-period}, for each $(h,E)\in\mathcal D$, the fibre $p^{-1}(h,E)$ is parametrized by the smooth periodic curve
\[
\mathtt u_{h,E}(s)=\bigl(
 h-\frac32r_-^2\operatorname{sn}^2(r_+s\mid m),\,
 r_-\operatorname{sn}(r_+s\mid m),\,
 r_-r_+\operatorname{cn}(r_+s\mid m)\operatorname{dn}(r_+s\mid m)
 \bigr),
\]
where $r_\pm$ and $m$ are defined in \eqref{eq:elliptic-parameters}. This parametrization depends smoothly on $(h,E)$, satisfies $\mathtt u_{h,E}(0)=(h,0,\sqrt E)$, and has fundamental period $T(h,E)=4K(m(h,E))/r_+(h,E)$. The map
\[
\mathtt q\!:\,(h,E,\vartheta)\in\mathcal D\times(\mathbb R/2\pi\mathbb Z)
\longmapsto\mathtt u_{h,E}\bigl(T(h,E)\vartheta/(2\pi)\bigr)\in\mathcal P
\]
is a smooth bijection. Since $p\circ\mathtt q=\operatorname{pr}_{\mathcal D}$ and $\partial_\vartheta\mathtt q\ne0$, its differential is everywhere invertible, so $\mathtt q$ is a diffeomorphism. It follows that
\[
\Lambda_\sigma\!:\,(\mathbf X,\vartheta)\in\mathcal V_\sigma^+\times(\mathbb R/2\pi\mathbb Z)
\longmapsto\bigl(\mathbf X,\mathtt q(\mathfrak c_-(\mathbf X),\vartheta)\bigr)
\in\mathcal V_\sigma^+\times_{\mathcal D}\mathcal P
\]
is a smooth diffeomorphism. Using this identification, we obtain the diffeomorphism
\begingroup\fontsize{10.5}{12.5}\selectfont
\[
\widetilde\Xi_\sigma\!:\,(\mathbf A,\mathtt u)\in(\mathbf J_-^\uparrow)^{-1}(\mathcal V_\sigma^+)
\longmapsto\bigl(\mathbf J_-^\uparrow(\mathbf A,\mathtt u),\operatorname{pr}_{\mathbb R/2\pi\mathbb Z}(\mathtt q^{-1}(\mathtt u)),C_\sigma(\mathbf A,\mathtt u)\bigr)
\in\mathcal V_\sigma^+\times(\mathbb R/2\pi\mathbb Z)\times K
\]
\endgroup
Identifying $K\simeq\mathbb T^2$, reordering the factors, and taking the inverse, we obtain the smooth trivializations
\begingroup\fontsize{9.5}{11.5}\selectfont
\begin{equation}\label{eq:torus-trivialization}
F_\sigma\!:\,(\mathbf X,\phi_1,\phi_2,\vartheta)\in\mathcal V_\sigma^+\times\mathbb T^3
\longmapsto\bigl(\mathbf B_\sigma(\mathbf X)R_1(\phi_1,\phi_2)\mathbf Q(\mathtt q(\mathfrak c_-(\mathbf X),\vartheta))^{-1},\mathtt q(\mathfrak c_-(\mathbf X),\vartheta)\bigr)
\in(\mathbf J_-^\uparrow)^{-1}(\mathcal V_\sigma^+).
\end{equation}
\endgroup

On each nonempty overlap $\mathcal V_\rho^+\cap\mathcal V_\sigma^+$, the spectral normalizations $\mathbf B_\rho$ and $\mathbf B_\sigma$ differ by an element of $K\simeq\mathbb T^2$. We may therefore write $\mathbf B_\rho(\mathbf X)^{-1}\mathbf B_\sigma(\mathbf X)=R_1\bigl(k_{\rho\sigma}(\mathbf X)\bigr)$ for a smooth map $k_{\rho\sigma}:\mathcal V_\rho^+\cap\mathcal V_\sigma^+\to\mathbb T^2$. Since the parametrization $\mathtt q$ is global, the transition functions of the trivializations \eqref{eq:torus-trivialization} take the form
\[
F_\rho^{-1}\circ F_\sigma(\mathbf X,\phi,\vartheta)
=\bigl(\mathbf X,\phi+k_{\rho\sigma}(\mathbf X),\vartheta\bigr),
\]
with $\phi=(\phi_1,\phi_2)$. Thus all transition functions act by fibrewise translations of $\mathbb T^3$, establishing the principal bundle structure. The invariance of the fibres under the Reeb flow follows from conservation of the moment map.

It remains to prove nontriviality. Fix $(h,E)\in\mathcal D$, and abbreviate $B_j^*(h,E)$ to $B_j^*$, $j=1,\ldots,5$ (see (3.8)). Consider the orthogonal decomposition
\[
\mathbb V_-:=\operatorname{span}(B_1^*,B_2^*),\qquad
\mathbb V_+:=\operatorname{span}(B_3^*,B_4^*,B_5^*),\qquad
\mathbb R^{2,3}=\mathbb V_-\oplus\mathbb V_+.
\]
By (3.8), $\mathbb V_-$ is negative-definite and $\mathbb V_+$ is positive-definite. We orient $\mathbb V_+$ by the ordered orthonormal basis $(B_3^*,B_4^*,B_5^*)$. For each vector $\mathrm N$ in the unit sphere $S^2\subset\mathbb V_+$, let $\mathbf X_{\mathrm N}\in\mathfrak a(2,3)$ be the matrix, in the standard Witt basis, of the skew-adjoint endomorphism defined by
\[
\mathbf X_{\mathrm N}(aB_1^*+bB_2^*+V)
=\omega(aB_2^*-bB_1^*)+\nu\,\mathrm N\times V,
\]
where $a,b\in\mathbb R$, $V\in\mathbb V_+$, and $\times$ denotes the cross product associated with the Euclidean metric and the chosen orientation of $\mathbb V_+$. By (3.9), we have $\mathbf X_{B_5^*}=\mathbf H_-(h,0,\sqrt E)$.

For any $\mathrm N\in S^2$, choose $R\in SO(\mathbb V_+)$ such that $R(B_5^*)=\mathrm N$, and extend $R$ to $\mathbb R^{2,3}$ by the identity on $\mathbb V_-$. Let $g_R$ denote the matrix of this isometry in the standard Witt basis. Since $SO(\mathbb V_+)$ is connected, $g_R\in\mathbb A_+^\uparrow(2,3)$. The equivariance of the cross product gives $\mathbf X_{\mathrm N}=g_R\mathbf H_-(h,0,\sqrt E)g_R^{-1}$. Thus $\mathbf X_{\mathrm N}\in\mathcal U_+$ for every $\mathrm N\in S^2$, and $\mathrm N\mapsto\mathbf X_{\mathrm N}$ defines a smooth map from $S^2$ into $\mathfrak c_-^{-1}(h,E)$.

We now identify the positive-definite spectral plane of $\mathbf X_{\mathrm N}$. Since $0<\nu<\omega$, the restriction of $\mathbf X_{\mathrm N}^2+\nu^2I$ to $\mathbb V_-$ is invertible. On the other hand, for $V\in\mathbb V_+$, the vector triple product identity gives $\mathrm N\times(\mathrm N\times V)=\langle\mathrm N,V\rangle\mathrm N-V$, and therefore $(\mathbf X_{\mathrm N}^2+\nu^2I)V=\nu^2\langle\mathrm N,V\rangle\mathrm N$. Consequently, $\ker(\mathbf X_{\mathrm N}^2+\nu^2I)=\mathrm N^\perp\cap\mathbb V_+=T_{\mathrm N}S^2$.

Suppose, for contradiction, that the principal $\mathbb T^3$-bundle is trivial. Then the restricted moment map admits a global smooth section $\mathfrak s:\mathcal U_+\to\widehat{\mathcal Y}_{-1}^\uparrow$. For each $\mathrm N\in S^2$, write $\mathfrak s(\mathbf X_{\mathrm N})=(\mathbf A(\mathrm N),\mathtt u(\mathrm N))$ and set $\mathbf S(\mathrm N)=\mathbf A(\mathrm N)\mathbf Q(\mathtt u(\mathrm N))$, with $\mathbf Q$ as in \eqref{eq:torus-global-normalization}. Since $\mathfrak s$ is a section, $\mathbf J_-^\uparrow(\mathbf A(\mathrm N),\mathtt u(\mathrm N))=\mathbf X_{\mathrm N}$. From (2.5), together with the normalization identity \eqref{eq:torus-global-spectral-identity} and $\mathfrak c_-(\mathbf X_{\mathrm N})=(h,E)$, we obtain $\mathbf X_{\mathrm N}\mathbf S(\mathrm N)=\mathbf S(\mathrm N)\mathbf D_{h,E}$.

Let $V_3(\mathrm N)$ and $V_4(\mathrm N)$ be the third and fourth columns of $\mathbf S(\mathrm N)$. By the preceding identity and the definition of $\mathbf D_{h,E}$, we have $\mathbf X_{\mathrm N}V_3(\mathrm N)=\nu V_4(\mathrm N)$ and $\mathbf X_{\mathrm N}V_4(\mathrm N)=-\nu V_3(\mathrm N)$. Hence both vectors belong to $\ker(\mathbf X_{\mathrm N}^2+\nu^2I)=T_{\mathrm N}S^2$. Moreover, $\mathbf S(\mathrm N)$ is a positively oriented and time-oriented pseudo-orthonormal frame. Thus $V_3(\mathrm N)$ and $V_4(\mathrm N)$ form an orthonormal basis of $T_{\mathrm N}S^2$, depending smoothly on $\mathrm N$. In particular, $\mathrm N\mapsto V_3(\mathrm N)$ is a smooth nowhere-vanishing tangent vector field on $S^2$, contradicting the hairy ball theorem.

Therefore the restricted moment map admits no global smooth section, and the principal $\mathbb T^3$-bundle is nontrivial.
\end{proof}

\subsection{Linearization of the Reeb flow}\label{subsec:linearization-reeb}

The trivializations constructed in Theorem~\ref{thm:invariant-torus-fibration}
provide angular coordinates on the invariant tori. We now introduce new
angular coordinates in which the Reeb flow becomes linear.

For $(h,E)\in\mathcal D$, recall that the normalized periodic solution of
the reduced Euler--Lagrange system is
\[
\mathtt u_{h,E}(s)=\bigl(\kappa_{h,E}(s),\tau_{h,E}(s),\tau_{h,E}'(s)\bigr),
\]
with $\kappa_{h,E}(s)=h-\frac32\tau_{h,E}(s)^2$ and $\tau_{h,E}$ given
by \eqref{eq:periodic-torsion-curvature}. For $a=\omega,\nu$, define
\begin{equation}\label{eq:angular-functions}
\resizebox{.84\linewidth}{!}{$\displaystyle
\theta_a(s;h,E)=-\int_0^s\frac{a(h,E)}{a(h,E)^2-\tau_{h,E}(\xi)^2}\,d\xi,\qquad
\Theta_a(h,E)=\frac{\theta_a(T(h,E);h,E)}{T(h,E)}
$}
\end{equation}
and set $\widetilde\theta_a(s;h,E)=\theta_a(s;h,E)-s\Theta_a(h,E)$.
Since $\tau_{h,E}$ has fundamental period $T(h,E)$, the functions
$\widetilde\theta_a(\,\cdot\,;h,E)$, representing the periodic parts of
the angular evolution, are smooth and $T(h,E)$-periodic, with
$\widetilde\theta_a(0;h,E)=0$. Here $s\in\mathbb R$ is the natural
parameter, while its residue class modulo $T(h,E)$ parametrizes the
reduced periodic orbit. In particular, the functions $\widetilde\theta_a$
descend to $\mathbb R/T(h,E)\mathbb Z$.

Define $\widehat T:\mathcal D\to(0,\infty)$ by
$\widehat T(h,E)=T(h,E)/(2\pi)$ (see \eqref{eq:fundamental-period}).
We now replace the trivializations $F_\sigma$ in
\eqref{eq:torus-trivialization} with new trivializations
$\widetilde F_\sigma$, defined by
\begin{equation}\label{eq:modified-torus-trivialization}
\resizebox{.84\linewidth}{!}{$\displaystyle
\widetilde F_\sigma(\mathbf X,\phi_1,\phi_2,\vartheta)=F_\sigma\bigl(\mathbf X,\phi_1+\widetilde\theta_\omega(\widehat T(\mathfrak c_-(\mathbf X))\vartheta;\mathfrak c_-(\mathbf X)),\phi_2+\widetilde\theta_\nu(\widehat T(\mathfrak c_-(\mathbf X))\vartheta;\mathfrak c_-(\mathbf X)),\vartheta\bigr).
$}
\end{equation}

\begin{theorem}\label{thm:linearization-reeb}
Let $\mathcal R$ denote the Reeb vector field on
$\widehat{\mathcal Y}_{-1}^\uparrow$. For $\sigma=1,2,3$, the
trivializations $\widetilde F_\sigma$ linearize the Reeb vector field.
More precisely,
\begin{equation}\label{eq:reeb-linearization}
\mathcal R=(\widetilde F_\sigma)_*\bigl((\Theta_\omega\circ\mathfrak c_-)\partial_{\phi_1}+(\Theta_\nu\circ\mathfrak c_-)\partial_{\phi_2}+(\widehat T\circ\mathfrak c_-)^{-1}\partial_{\vartheta}\bigr)
\end{equation}
on $(\mathbf J_-^\uparrow)^{-1}(\mathcal V_\sigma^+)$. Here the
coefficient functions are understood as pulled back to
$\mathcal V_\sigma^+\times\mathbb T^3$ by the first projection, and
the coordinate vector fields on $\mathbb T^3$ are regarded as vertical
vector fields on this product. In particular, the angular frequencies
depend on $\mathbf X$ only through its Casimir invariants.
\end{theorem}

\begin{proof}
Consider the smooth map defined by
\[
\bm{\mathcal Q}:\ (s,\mathbf X)\in\mathbb R\times\mathcal U_+\longmapsto\mathbf Q\bigl(\mathtt u_{\mathfrak c_-(\mathbf X)}(s)\bigr)\in\mathcal F_+^\uparrow(2,3).
\]
Fix $\sigma\in\{1,2,3\}$ and $\mathbf X\in\mathcal V_\sigma^+$. In
what follows, primes denote differentiation with respect to $s$, with
$\mathbf X$ fixed. Differentiating the spectral identity
$\mathbf H_-\bigl(\mathtt u_{\mathfrak c_-(\mathbf X)}(s)\bigr)\bm{\mathcal Q}(s;\mathbf X)
=\bm{\mathcal Q}(s;\mathbf X)\mathbf D_{\mathfrak c_-(\mathbf X)}$ and using the
Lax equation \eqref{eq:lax}, we find that
\[
\mathbf L(s;\mathbf X)=\bm{\mathcal Q}(s;\mathbf X)^{-1}\bigl(\bm{\mathcal Q}'(s;\mathbf X)+\mathbf K_-(\kappa_{\mathfrak c_-(\mathbf X)}(s),\tau_{\mathfrak c_-(\mathbf X)}(s))\bm{\mathcal Q}(s;\mathbf X)\bigr)
\]
belongs to the Lie algebra of $K$. Hence
\[
\mathbf L(s;\mathbf X)=\operatorname{diag}\bigl(c_\omega(s;\mathbf X)J_0,\,c_\nu(s;\mathbf X)J_0,\,0\bigr)
\]
for suitable smooth functions
$c_\omega,c_\nu:\mathbb R\times\mathcal U_+\to\mathbb R$. For $a=\omega,\nu$, set
\[
W_a(s;\mathbf X):=W_a\bigl(\mathtt u_{\mathfrak c_-(\mathbf X)}(s)\bigr),\qquad V_a(s;\mathbf X):=V_a\bigl(\mathtt u_{\mathfrak c_-(\mathbf X)}(s)\bigr).
\]
In particular,
\[
W_a'(s;\mathbf X)+\mathbf K_-(\kappa_{\mathfrak c_-(\mathbf X)}(s),\tau_{\mathfrak c_-(\mathbf X)}(s))W_a(s;\mathbf X)=c_a(s;\mathbf X)V_a(s;\mathbf X),\qquad a=\omega,\nu.
\]
Set
\[
\delta_\omega(s;\mathbf X)=\sqrt{\omega(\mathfrak c_-(\mathbf X))^2-\tau_{\mathfrak c_-(\mathbf X)}(s)^2},\quad \delta_\nu(s;\mathbf X)=\sqrt{\nu(\mathfrak c_-(\mathbf X))^2-\tau_{\mathfrak c_-(\mathbf X)}(s)^2}.
\]
Since the fourth component of $W_a(s;\mathbf X)$ vanishes for
$a=\omega,\nu$, comparing the fourth components of the preceding
identities and using \eqref{eq:torus-global-vectors}, we obtain
\[
\scalebox{0.99}{$\displaystyle
c_\omega(s;\mathbf X)\frac{\delta_\omega(s;\mathbf X)}{\sqrt{d(\mathfrak c_-(\mathbf X))}}=-\frac{\omega(\mathfrak c_-(\mathbf X))}{\sqrt{d(\mathfrak c_-(\mathbf X))}\,\delta_\omega(s;\mathbf X)},\quad c_\nu(s;\mathbf X)\frac{\delta_\nu(s;\mathbf X)}{\sqrt{d(\mathfrak c_-(\mathbf X))}}=-\frac{\nu(\mathfrak c_-(\mathbf X))}{\sqrt{d(\mathfrak c_-(\mathbf X))}\,\delta_\nu(s;\mathbf X)}.
$}
\]
Thus, by \eqref{eq:angular-functions},
\[
c_\omega(s;\mathbf X)=-\frac{\omega(\mathfrak c_-(\mathbf X))}{\delta_\omega(s;\mathbf X)^2}=\partial_s\theta_\omega(s;\mathfrak c_-(\mathbf X)),\quad c_\nu(s;\mathbf X)=-\frac{\nu(\mathfrak c_-(\mathbf X))}{\delta_\nu(s;\mathbf X)^2}=\partial_s\theta_\nu(s;\mathfrak c_-(\mathbf X)),
\]
which gives
\begingroup
\fontsize{9.5}{11.5}\selectfont
\begin{equation}\label{eq:W-differential}
\bm{\mathcal Q}'(s;\mathbf X)+\mathbf K_-\bigl(\kappa_{\mathfrak c_-(\mathbf X)}(s),\tau_{\mathfrak c_-(\mathbf X)}(s)\bigr)\bm{\mathcal Q}(s;\mathbf X)=\bm{\mathcal Q}(s;\mathbf X)\operatorname{diag}\bigl(\partial_s\theta_\omega(s;\mathfrak c_-(\mathbf X))J_0,\partial_s\theta_\nu(s;\mathfrak c_-(\mathbf X))J_0,0\bigr).
\end{equation}
\endgroup

Let $t\mapsto(\mathbf A(t),\mathtt u(t))$ be a Reeb trajectory with
moment $\mathbf X$, and write
\[
(\mathbf A(t),\mathtt u(t))=F_\sigma\bigl(\mathbf X,\phi_1(t),\phi_2(t),\vartheta(t)\bigr).
\]
Choosing local real lifts of the angular coordinates
$\phi_1(t),\phi_2(t),\vartheta(t)$, set
$s(t)=\widehat T(\mathfrak c_-(\mathbf X))\vartheta(t)$, so that
$\mathtt u(t)=\mathtt u_{\mathfrak c_-(\mathbf X)}(s(t))$. In what follows, dots
denote differentiation with respect to $t$. Since
$\mathtt u_{\mathfrak c_-(\mathbf X)}(s)$ is a solution of
\eqref{eq:reeb-dynamics}, the reduced Reeb equations imply
$\mathtt u_{\mathfrak c_-(\mathbf X)}'(s(t))\dot s(t)=\mathtt u_{\mathfrak c_-(\mathbf X)}'(s(t))$.
Since the reduced periodic orbit is regular,
$\mathtt u_{\mathfrak c_-(\mathbf X)}'$ never vanishes. Hence $\dot s=1$, and
therefore $\dot\vartheta=\widehat T(\mathfrak c_-(\mathbf X))^{-1}$. By
\eqref{eq:torus-trivialization}, the $K$-component of the trajectory is
\[
\mathbf C(t)=\mathbf B_\sigma(\mathbf X)^{-1}\mathbf A(t)\bm{\mathcal Q}(s(t);\mathbf X)=R_1(\phi_1(t),\phi_2(t)).
\]
Since $\mathbf X$ is constant and $\dot s=1$, differentiating this
identity and using \eqref{eq:W-differential} together with the Reeb
equation $\dot{\mathbf A}=\mathbf A\mathbf K_-(\kappa,\tau)$, we obtain
\[
\dot{\mathbf C}(t)=\mathbf C(t)\operatorname{diag}\bigl(\partial_s\theta_\omega(s(t);\mathfrak c_-(\mathbf X))J_0,\partial_s\theta_\nu(s(t);\mathfrak c_-(\mathbf X))J_0,0\bigr).
\]
Consequently,
\[
\dot\phi_1(t)=\partial_s\theta_\omega(s(t);\mathfrak c_-(\mathbf X)),\qquad \dot\phi_2(t)=\partial_s\theta_\nu(s(t);\mathfrak c_-(\mathbf X)).
\]
By \eqref{eq:modified-torus-trivialization}, the new angular coordinates
are $\phi_1(t)-\widetilde\theta_\omega(s(t);\mathfrak c_-(\mathbf X))$ and
$\phi_2(t)-\widetilde\theta_\nu(s(t);\mathfrak c_-(\mathbf X))$. Since
\[
\partial_s\widetilde\theta_a(s;\mathfrak c_-(\mathbf X))=\partial_s\theta_a(s;\mathfrak c_-(\mathbf X))-\Theta_a(\mathfrak c_-(\mathbf X)),\qquad \dot s=1,
\]
their velocities are $\Theta_\omega(\mathfrak c_-(\mathbf X))$ and
$\Theta_\nu(\mathfrak c_-(\mathbf X))$, respectively. Together with
$\dot\vartheta=\widehat T(\mathfrak c_-(\mathbf X))^{-1}$, this proves
\eqref{eq:reeb-linearization}.
\end{proof}

\begin{remark}
Let $K_{\max}\simeq SO(2)\times SO(3)$ be a maximal compact subgroup of
$O_+^\uparrow(2,3)$ containing $K$.
The Cartan decomposition $O_+^\uparrow(2,3)=\exp(\mathfrak p)K_{\max}$,
with $\mathfrak p\simeq\mathbb R^6$ (see \cite{Helgason2001}), and the
spectral construction in the proof of Theorem~\ref{thm:invariant-torus-fibration}
yield a global diffeomorphism
$\mathcal U_+\simeq\mathcal D\times(K_{\max}/K)\times\mathfrak p
\simeq\mathcal D\times S^2\times\mathbb R^6$.
Under this identification, the domains $\mathcal V_\sigma^+$ take the form
$\mathcal D\times W_\sigma$, with $W_\sigma\subset S^2\times\mathbb R^6$,
and the transition functions of $\widetilde F_\sigma$ depend only on the
transverse variables.

The principal $\mathbb T^3$-action associated with $\widetilde F_\sigma$
preserves $\alpha_-$. The product structure and the action-independent
transition functions allow one to construct, by radial integration in
the action variables, a smooth map $\mathbf t:\mathcal U_+\to\mathbb T^2$
such that the modified
trivializations\footnote{The dot denotes the principal right action of
$\mathbb T^3$ on $\widehat{\mathcal Y}_{-1}^\uparrow$.}
$\widehat F_\sigma(\mathbf X,\boldsymbol\phi)
=\widetilde F_\sigma(\mathbf X,\boldsymbol\phi)
\cdot(\mathbf t(\mathbf X),0)$ satisfy
\[
\widehat F_\sigma^*\alpha_-
=-\frac{\omega(h,E)}{2}\,d\phi_1
-\frac{\nu(h,E)}{2}\,d\phi_2
+c(h,E)\,d\vartheta+\beta_\sigma(h,E),
\]
where $c$ is smooth on $\mathcal D$ and $\beta_\sigma(h,E)$ is a smooth
family of $1$-forms on $W_\sigma$.
Since $(\omega,\nu)$ are global coordinates on $\mathcal D$, this is the
generalized contact action-angle normal form of Jovanovi\'c
\cite[Theorem~5.1 and Definition~5.2]{Jovanovic2012}.
\end{remark}

We now consider conformal extremals of negative helicity with
nonconstant periodic conformal torsion and parameters $(h,E)\in\mathcal D$.
Up to orientation-preserving conformal transformations, these are
precisely the projections onto the Einstein universe of the integral
curves of the Reeb vector field restricted to
$\widehat{\mathcal Y}_{-1}^\uparrow$. The linearization of the Reeb flow,
together with the spectral normalization, provides a closed-form
representation of these extremals.

\begin{proposition}\label{prop:explicit-integration}
Let $\gamma$ be a conformal extremal of negative helicity with
nonconstant periodic conformal torsion and parameters $(h,E)\in\mathcal D$.
Then, after a suitable translation of its natural parameter, there
exists a positively oriented pseudo-orthonormal basis
$\mathbf B=(B_1,\ldots,B_5)$ of $\mathbb R^{2,3}$ such that
$\gamma(s)=[\Gamma(s)]$, where $\Gamma:\mathbb R\to\mathbb R^{2,3}$
is defined, with the abbreviation $\theta_a(s)=\theta_a(s;h,E)$,
by the expression
\begingroup
\fontsize{10}{12}\selectfont
\begin{equation}\label{eq:gamma-explicit}
\frac{\sqrt{\omega^2-\tau(s)^2}}{\omega\sqrt d}\bigl(\cos\theta_\omega(s)\,B_1+\sin\theta_\omega(s)\,B_2\bigr)+\frac{\sqrt{\nu^2-\tau(s)^2}}{\nu\sqrt d}\bigl(\cos\theta_\nu(s)\,B_3+\sin\theta_\nu(s)\,B_4\bigr)+\frac{\tau(s)}{\omega\nu}B_5.
\end{equation}
\endgroup
Here $\omega=\omega(h,E)$, $\nu=\nu(h,E)$ are as in
\eqref{eq:normal-form-frequencies}, $d=\omega^2-\nu^2$, and
$\tau(s)=\tau_{h,E}(s)$ is given by \eqref{eq:periodic-torsion-curvature}.
The angular functions $\theta_\omega$ and $\theta_\nu$
are defined in \eqref{eq:angular-functions}. Conversely, for every
$(h,E)\in\mathcal D$ and every positively oriented pseudo-orthonormal
basis $\mathbf B$, the curve $\gamma(s)=[\Gamma(s)]$ defined by
\eqref{eq:gamma-explicit} is a conformal extremal of the prescribed type.
\end{proposition}

\begin{proof}
After a suitable translation of the natural parameter, the reduced
trajectory is $\mathtt u_{h,E}(s)$. Let $\mathbf A(s)$ be the canonical
frame along $\gamma$, and set
$\mathbf W(s;h,E)=\mathbf Q(\mathtt u_{h,E}(s))$.
By the Frenet equations and \eqref{eq:W-differential}, the matrix
$\mathbf S(s)=\mathbf A(s)\mathbf W(s;h,E)$ satisfies
\[
\mathbf S'=\mathbf S\operatorname{diag}\bigl(\partial_s\theta_\omega(s;h,E)J_0,\,\partial_s\theta_\nu(s;h,E)J_0,\,0\bigr).
\]
Integrating and using
$\theta_\omega(0;h,E)=\theta_\nu(0;h,E)=0$, we obtain
\[
\mathbf S(s)=\mathbf B\,R_1\bigl(\theta_\omega(s;h,E),\theta_\nu(s;h,E)\bigr),
\]
where the columns of $\mathbf B=\mathbf A(0)\mathbf W(0;h,E)$ form a
positively oriented pseudo-orthonormal basis. Consequently,
\[
\mathbf A(s)=\mathbf B\,R_1\bigl(\theta_\omega(s;h,E),\theta_\nu(s;h,E)\bigr)
\mathbf W(s;h,E)^{-1}.
\]
Using \eqref{eq:torus-global-vectors}--\eqref{eq:torus-global-normalization},
we obtain
\[
\mathbf W(s;h,E)^{-1}E_1=\bigl(\frac{(\omega^2-\tau(s)^2)^{1/2}}{\omega d^{1/2}},\,0,\,\frac{(\nu^2-\tau(s)^2)^{1/2}}{\nu d^{1/2}},\,0,\,\frac{\tau(s)}{\omega\nu}\bigr)^T.
\]
Taking the first column of $\mathbf A(s)$ now gives
\eqref{eq:gamma-explicit}.

For the converse, we may first assume that $\mathbf B$ is time-oriented,
since the general case follows by applying an orientation-preserving
conformal transformation that reverses the time-orientation. Set
\[
\mathbf A(s)
=\mathbf B\,R_1\bigl(\theta_\omega(s;h,E),\theta_\nu(s;h,E)\bigr)
\mathbf W(s;h,E)^{-1}.
\]
By \eqref{eq:W-differential} and the reduced Euler--Lagrange equations
\eqref{eq:reeb-dynamics}, the pair $(\mathbf A,\mathtt u_{h,E})$ is a
Reeb trajectory. The conclusion follows from the Reeb characterization
of conformal extremals.
\end{proof}

\begin{remark}
The angular functions $\theta_a$ and their mean frequencies $\Theta_a$
admit the representations
\begin{equation}\label{eq:theta-elliptic}
\resizebox{.84\linewidth}{!}{$\displaystyle
\theta_a(s;h,E)=-\frac{1}{ar_+}\Pi\bigl(\frac{r_-^2}{a^2};r_+s\mid m\bigr),\qquad
\Theta_a(h,E)=-\frac{1}{aK(m)}\Pi\bigl(\frac{r_-^2}{a^2}\mid m\bigr),\quad a=\omega,\nu.
$}
\end{equation}
Here $a=a(h,E)$, and $r_\pm$ and $m$ are evaluated at $(h,E)$.
The symbol $\Pi$ denotes the elliptic integral of the third kind in
Jacobi's notation (see \cite[\S~17.2, (17.2.16)]{AS}). Thus,
\eqref{eq:gamma-explicit} provides explicit parametrizations of these
conformal extremals in terms of Jacobi elliptic functions and elliptic
integrals of the third kind.
\end{remark}

\begin{numbereddefinition}\label{def:period-map}
Define the period map $\Psi=(\Psi_1,\Psi_2):\mathcal D\longrightarrow\mathbb R^2$ by
\begin{equation}\label{eq:integration:15}
\Psi_1(h,E)=\frac{2}{\pi\omega r_+}
\Pi\bigl(\frac{r_-^2}{\omega^2}\bigm|m\bigr),\qquad
\Psi_2(h,E)=\frac{2}{\pi\nu r_+}
\Pi\bigl(\frac{r_-^2}{\nu^2}\bigm|m\bigr).
\end{equation}
Here $\omega,\nu,r_\pm,m$ are evaluated at $(h,E)$.
\end{numbereddefinition}

The frame representation obtained in the proof of
Proposition~\ref{prop:explicit-integration} also determines the monodromy,
leading to the following closure criterion.

\begin{proposition}\label{prop:closing-criterion}
Let $\gamma$ be a conformal extremal of negative helicity with
nonconstant periodic conformal torsion and parameters $(h,E)\in\mathcal D$.
Then $\gamma$ is closed if and only if $\Psi(h,E)\in\mathbb Q^2$.
\end{proposition}

\begin{proof}
Put $T=T(h,E)$, the minimal common period of the conformal invariants.
Since $\mathtt u_{h,E}$ is $T$-periodic, so is
$\mathbf W(s;h,E)=\mathbf Q(\mathtt u_{h,E}(s))$. Moreover, by the
periodicity of $\partial_s\theta_a$ and the definition of the
mean angular frequencies, we have
\[
\theta_a(s+T;h,E)=\theta_a(s;h,E)+T\Theta_a(h,E),\qquad a=\omega,\nu.
\]
Thus, the expression for the canonical frame obtained in the proof of
Proposition~\ref{prop:explicit-integration} gives
$\mathbf A(s+T)=\mathbf M_\gamma\mathbf A(s)$, where
\[
\mathbf M_\gamma
=\mathbf B\,R_1\bigl(T\Theta_\omega(h,E),T\Theta_\nu(h,E)\bigr)
\mathbf B^{-1}
\]
is the monodromy of $\gamma$ (see Definition~\ref{def:monodromic}).
By Remark~\ref{rem:monodromic-properties}, $\gamma$ is closed if and
only if its monodromy has finite order. Since $R_1$ consists of two
planar rotations, this is equivalent to
$T\Theta_\omega(h,E)/(2\pi)\in\mathbb Q$ and
$T\Theta_\nu(h,E)/(2\pi)\in\mathbb Q$.
By \eqref{eq:theta-elliptic} and $T=4K(m)/r_+$,
\[
T\Theta_\omega(h,E)/(2\pi)=-\Psi_1(h,E),\qquad
T\Theta_\nu(h,E)/(2\pi)=-\Psi_2(h,E).
\]
Hence these rationality conditions are equivalent to
$\Psi(h,E)\in\mathbb Q^2$.
\end{proof}

\subsection{The global geometry of the period map and the classification of closed extremals}\label{subsec:period-map}

We now investigate the global properties of the period map and their
consequences for the classification of closed conformal extremals.

\begin{theorem}\label{thm:period-map}
The period map $\Psi:\D\longrightarrow\R^2$ is a real-analytic diffeomorphism of $\D$ onto the domain
\[
\widetilde{\mathcal D}=\left\{(x,y)\in\R^2:0<x<y,\quad x+y>1,\quad x^2+y^2<1\right\}.
\]
\end{theorem}

\begin{proof}
Let $I=(0,1)$ and consider the real-analytic diffeomorphism
\[
\mathcal F:(m,q)\in I\times I\longmapsto
\bigl(\frac{1+m}{q(1-m)},\frac{4m}{q^2(1-m)^2}\bigr)\in\D.
\]
Its inverse is given by $q=1/\sqrt{h^2-E}$ and $m=(h-\sqrt{h^2-E})/(h+\sqrt{h^2-E})$. Set $\widetilde\Psi=\Psi\circ\mathcal F$. We shall prove that $\widetilde\Psi:I\times I\longrightarrow\widetilde{\mathcal D}$ is a real-analytic diffeomorphism. Writing $\widetilde\Psi(m,q)=\bigl(\widetilde\Psi_1(m,q),\widetilde\Psi_2(m,q)\bigr)$, substitution of \eqref{eq:integration:15} gives
\[
\widetilde\Psi_j(m,q)
=\frac{\sqrt2\,q(1-m)}{\pi\sqrt{D_j}}
\Pi\!\bigl(\frac{2m}{D_j}\bigm|m\bigr),
\qquad j=1,2,
\]
where $D_1=1+m+(1-m)c$, $D_2=1+m-(1-m)c$, and $c=\sqrt{1-q^2}$. Since $0<c<1$, we have $D_1-2m=(1-m)(1+c)>0$ and $D_2-2m=(1-m)(1-c)>0$. Moreover, $D_j<2$, and hence $m<2m/D_j<1$, $j=1,2$.

Recall that $\Pi(n\mid m)$ is real analytic on $(0,1)\times(0,1)$ and admits a continuous extension to $[0,1)\times[0,1)$; see \cite[\S19.2(ii), Eqs.~(19.2.7)--(19.2.8)]{dlmf}.\footnote{We use the parameter $m=k^2$, whereas the NIST DLMF uses the modulus $k$; moreover, its characteristic $\alpha^2$ corresponds to our $n$.} In particular, $\Pi(n\mid0)=\pi/(2\sqrt{1-n})$ for $0\le n<1$; see \cite[\S19.6(i), Eq.~(19.6.3)]{dlmf}. It follows that $\widetilde\Psi$ extends continuously to the side $m=0$, and direct substitution gives
\begin{equation}\label{eq:integration:16}
\widetilde\Psi(0,q)
=\frac{1}{\sqrt2}\bigl(
(1-\sqrt{1-q^2})^{1/2},
(1+\sqrt{1-q^2})^{1/2}
\bigr),\qquad 0<q\le1.
\end{equation}
On the side $q=1$, one has $c=0$ and $D_1=D_2=1+m$, so that
\[
\widetilde\Psi_1(m,1)=\widetilde\Psi_2(m,1)
=\frac{\sqrt2(1-m)}{\pi\sqrt{1+m}}
\Pi\!\bigl(\frac{2m}{1+m}\bigm|m\bigr),
\qquad 0\le m<1.
\]
Thus $\widetilde\Psi$ extends continuously to $[0,1)\times(0,1]$.

We next compute its derivatives. Let $\Eell(m)$ denote the complete elliptic integral of the second kind \cite[p.~590, formulas~17.3.3--17.3.4]{AS}. We use the standard differentiation formulas for the complete elliptic integral of the third kind; see \cite[formulas 08.03.20.0001 and 08.03.20.0003]{wolfram-functions}:
\[
\begin{aligned}
\partial_n\Pi(n\mid m)
&=\frac{n\Eell(m)+(m-n)K(m)+(n^2-m)\Pi(n\mid m)}{2(m-n)(n-1)n},\\[6pt]
\partial_m\Pi(n\mid m)
&=\frac{\Eell(m)-(1-m)\Pi(n\mid m)}{2(1-m)(m-n)}.
\end{aligned}
\]
The chain rule gives
\[
\begin{aligned}
\partial_m\widetilde\Psi_j
&=\frac{\sqrt2\,q}{2\pi m\sqrt{D_j}}
\bigl(\Eell(m)-K(m)\bigr),\\[6pt]
\partial_q\widetilde\Psi_j
&=\frac{\sqrt2\,\sigma_j}{2\pi c\sqrt{D_j}}
\left[2\Eell(m)-(1-m)(1-\sigma_jc)K(m)\right],
\end{aligned}
\qquad j=1,2,
\]
where $\sigma_1=1$ and $\sigma_2=-1$.

The standard inequality $(1-m)K(m)<\Eell(m)<K(m)$ for $0<m<1$, which follows from \cite[\S19.9(i), Eq.~(19.9.8)]{dlmf}, implies that $\partial_m\widetilde\Psi_1<0$, $\partial_m\widetilde\Psi_2<0$, $\partial_q\widetilde\Psi_1>0$, and $\partial_q\widetilde\Psi_2<0$. These monotonicity properties imply injectivity: equality of the first components at two points with $m_1<m_2$ forces $q_1<q_2$, contradicting equality of the second components. Moreover,
\[
\det D\widetilde\Psi
=\frac{q}{\pi^2mc\sqrt{D_1D_2}}
\bigl(K(m)-\Eell(m)\bigr)
\bigl(2\Eell(m)-(1-m)K(m)\bigr)>0.
\]
Hence $\widetilde\Psi$ is a one-to-one real-analytic local diffeomorphism. To determine its behavior as $m\to1^-$, fix $0<q\le1$ and put $\delta=1-m$, $n_j=2m/D_j$, and
\[
a_{j,\delta}=\frac{1-n_j}{\delta}
=\frac{1+\sigma_jc}{D_j}
\longrightarrow a_j=\frac{1+\sigma_jc}{2}\in(0,1),
\qquad j=1,2.
\]
By the integral representation of the complete elliptic integral of the third kind, followed by the substitutions $x=\cot\varphi$ and $x=\sqrt\delta\,t$, we obtain
\[
\delta\,\Pi(n_j\mid1-\delta)
=\int_0^\infty
\frac{\sqrt{1+\delta t^2}}{(t^2+a_{j,\delta})\sqrt{1+t^2}}\,dt.
\]
For sufficiently small $\delta>0$, we have $a_{j,\delta}\ge a_j/2$ and $\delta<1$, so the integrand is bounded by the integrable function $(t^2+a_j/2)^{-1}$. Dominated convergence, followed by the substitution $u=t/\sqrt{1+t^2}$, therefore gives
\[
\lim_{\delta\to0^+}\delta\,\Pi(n_j\mid1-\delta)
=\int_0^\infty\frac{dt}{(t^2+a_j)\sqrt{1+t^2}}
=\frac{\arccos\sqrt{a_j}}{\sqrt{a_j(1-a_j)}}.
\]
Substituting $a_1=(1+c)/2$, $a_2=(1-c)/2$, and $c=\sqrt{1-q^2}$ into the expressions for $\widetilde\Psi_j$, we obtain
\begin{equation}\label{eq:integration:17}
\lim_{m\to1^-}\widetilde\Psi_1(m,q)=\pi^{-1}\arcsin q,
\qquad
\lim_{m\to1^-}\widetilde\Psi_2(m,q)=1-\pi^{-1}\arcsin q.
\end{equation}

The monotonicity of $\widetilde\Psi_1$ and $\widetilde\Psi_2$ with respect to $q$ implies that, for every $q_0\in(0,1]$,
\[
\lim_{\substack{m\to1^-\\q\to q_0}}\widetilde\Psi(m,q)
=\frac1\pi\bigl(\arcsin q_0,\pi-\arcsin q_0\bigr).
\]
Indeed, for $q_-<q<q_+$,
\[
\widetilde\Psi_1(m,q_-)\le\widetilde\Psi_1(m,q)\le\widetilde\Psi_1(m,q_+),\qquad
\widetilde\Psi_2(m,q_+)\le\widetilde\Psi_2(m,q)\le\widetilde\Psi_2(m,q_-).
\]
For fixed $q_-<q_0<q_+$, these inequalities hold whenever $q$ is sufficiently close to $q_0$. Applying \eqref{eq:integration:17} to the bounding terms and then letting $q_-\uparrow q_0$ and $q_+\downarrow q_0$ proves the asserted limit. The case $q_0=1$ follows similarly, with $q_+=1$ and $q\le1$.

We now determine the image of $\widetilde\Psi$. Let $(x,y)=\widetilde\Psi(m,q)$. Clearly $x>0$. Since $\widetilde\Psi_1$ is strictly increasing in $q$, while $\widetilde\Psi_2$ is strictly decreasing in $q$, and the two components coincide for $q=1$, we have $x<y$. Since both components are strictly decreasing in $m$, \eqref{eq:integration:17} gives $x>\pi^{-1}\arcsin q$ and $y>1-\pi^{-1}\arcsin q$, and therefore $x+y>1$. On the other hand, if $(x_0,y_0)=\widetilde\Psi(0,q)$, then \eqref{eq:integration:16} gives $x_0^2+y_0^2=1$. Since both components are positive and strictly decreasing in $m$, we have $x^2+y^2<x_0^2+y_0^2=1$. Hence $\widetilde\Psi(I\times I)\subset\widetilde{\mathcal D}$.

We next examine the remaining boundary behavior. For fixed $q\in I$, monotonicity with respect to $m$, together with \eqref{eq:integration:16} and \eqref{eq:integration:17}, gives
\[
0<\widetilde\Psi_1(m,q)<\frac{(1-\sqrt{1-q^2})^{1/2}}{\sqrt2},\;
1-\frac{\arcsin q}{\pi}<\widetilde\Psi_2(m,q)<\frac{(1+\sqrt{1-q^2})^{1/2}}{\sqrt2}.
\]
The corresponding weak inequalities hold also for $m=0,1$, by continuity. They imply that $\widetilde\Psi(m,q)\to(0,1)$ as $q\to0^+$, uniformly in $m\in[0,1]$. Hence $\widetilde\Psi$ extends continuously to $[0,1]^2$, with
\begin{equation}\label{eq:integration:18}
\widetilde\Psi(1,q)
=\pi^{-1}\left(\arcsin q,\,\pi-\arcsin q\right),
\qquad 0<q\le1,
\end{equation}
and $\widetilde\Psi(m,0)=(0,1)$.

By \eqref{eq:integration:16}, the side $m=0$ is mapped onto the circular arc $x^2+y^2=1$, $0\le x\le1/\sqrt2$, $y\ge x$. By \eqref{eq:integration:18}, the side $m=1$ is mapped onto the line segment $x+y=1$, $0\le x\le1/2$. On the side $q=1$, the two components coincide, and monotonicity in $m$, together with the endpoint values $1/\sqrt2$ and $1/2$, shows that this side is mapped onto the diagonal segment $x=y$, $1/2\le x\le1/\sqrt2$. Finally, the side $q=0$ is collapsed to $(0,1)$. Consequently, $\widetilde\Psi\bigl(\partial[0,1]^2\bigr)=\partial\widetilde{\mathcal D}$.

It remains to prove surjectivity. Let $U=\widetilde\Psi(I\times I)$. Since $\widetilde\Psi$ is a local diffeomorphism, $U$ is nonempty and open in $\widetilde{\mathcal D}$. The boundary identity established above also gives $U=\widetilde\Psi([0,1]^2)\cap\widetilde{\mathcal D}$. Since $\widetilde\Psi([0,1]^2)$ is compact, $U$ is relatively closed in $\widetilde{\mathcal D}$. By connectedness, $U=\widetilde{\mathcal D}$. Hence $\widetilde\Psi$ is bijective, and its inverse is real analytic because its Jacobian never vanishes. This completes the proof.
\end{proof}

As a consequence, we obtain the following.

\begin{theorem}\label{thm:closed-classification}
The equivalence classes of closed conformal extremals with nonconstant conformal torsion are in one-to-one correspondence with $\widetilde{\mathcal D}^{\,*}=\widetilde{\mathcal D}\cap\Q^2$.
\end{theorem}

\begin{proof}
Let $\xi\in\widetilde{\mathcal D}^{\,*}$. By Theorem~\ref{thm:period-map}, there exists a unique $(h,E)\in\D$ such that $\Psi(h,E)=\xi$. The conformal curvature and torsion given by \eqref{eq:periodic-torsion-curvature} satisfy the Euler--Lagrange equations. By Corollary~\ref{cor:EL-equations}, these invariants are realized by a conformal extremal of negative helicity, which is closed by Proposition~\ref{prop:closing-criterion}.

By Proposition~\ref{prop:monodromic-classification}, all nonconstant periodic solutions with fixed $(h,E)$ differ by a translation of the natural parameter. The uniqueness part of Proposition~\ref{prop:canonical-moving-frame} then implies that the corresponding extremals are equivalent. Thus each $\xi\in\widetilde{\mathcal D}^{\,*}$ determines a unique equivalence class of closed conformal extremals.

Conversely, let $\gamma$ be a closed conformal extremal with nonconstant conformal torsion. By Proposition~\ref{prop:spectral-restrictions}, its helicity is negative and its parameters $(h,E)$ belong to $\D$. Proposition~\ref{prop:closing-criterion} then implies that $\Psi(h,E)\in\widetilde{\mathcal D}^{\,*}$. Since $h$ and $E$ are invariants of conformal equivalence, distinct points of $\widetilde{\mathcal D}^{\,*}$ determine distinct equivalence classes. This proves the theorem.
\end{proof}

\begin{numbereddefinition}\label{def:characteristic-numbers}
Let $\gamma$ be a closed conformal extremal with parameters $(h,E)\in\D$. Write $\Psi_1(h,E)=p_1/\ell_1$ and $\Psi_2(h,E)=p_2/\ell_2$, where $p_j,\ell_j$ are relatively prime integers and $\ell_j>0$, $j=1,2$. The integers $p_1,\ell_1,p_2,\ell_2$ are called the characteristic numbers of $\gamma$. The geometric significance of the characteristic numbers will be discussed in Section~\ref{sec:geometry-closed-extremals}.
\end{numbereddefinition}

\section{Geometry of closed conformal extremals}\label{sec:geometry-closed-extremals}

In this section we study the global geometry of closed conformal extremals. The non-confinement result applies to all such curves; for the remaining results, we assume that the conformal torsion is nonconstant.

Let $\gamma$ be a closed conformal extremal with $\varepsilon=-1$ and $(h,E)\in\mathcal D$. By Proposition~\ref{prop:explicit-integration}, after a suitable translation of the natural parameter and up to conformal congruence, its representation in the $S^1\times S^2$-model associated with an adapted positively oriented and time-oriented pseudo-orthonormal basis $\mathbf B=(B_1,\ldots,B_5)$ is
\begin{equation}\label{eq:closed-gamma-hat}
\begin{aligned}
\widehat\gamma(s)
&=\bigl(\cos\theta_\omega(s),\sin\theta_\omega(s);\,r(s)\cos\theta_\nu(s),r(s)\sin\theta_\nu(s),z(s)\bigr),\\[5pt]
r(s)
&=\frac{\omega}{\nu}\sqrt{\frac{\nu^2-\tau(s)^2}{\omega^2-\tau(s)^2}},\qquad
z(s)=\frac{\sqrt{\omega^2-\nu^2}\,\tau(s)}{\nu\sqrt{\omega^2-\tau(s)^2}}.
\end{aligned}
\end{equation}

Let $T=T(h,E)$ be the fundamental period of the conformal torsion, given by \eqref{eq:fundamental-period}. By the quasiperiodicity relation established in the proof of Proposition~\ref{prop:closing-criterion} and the definition of the period map in \eqref{eq:integration:15}, we have $\theta_\omega(s+T)=\theta_\omega(s)-2\pi\Psi_1(h,E)$ and $\theta_\nu(s+T)=\theta_\nu(s)-2\pi\Psi_2(h,E)$. Consequently, $\widehat\gamma(s+T)=\widehat{\mathbf M}_\gamma\widehat\gamma(s)$, where
\[
\widehat{\mathbf M}_\gamma
=R_1\bigl(-2\pi\Psi_1(h,E),-2\pi\Psi_2(h,E)\bigr)
\in O_+^\uparrow(2,3)
\]
is the monodromy of $\gamma$ with respect to the adapted pseudo-orthonormal basis $\mathbf B$.

Put $\xi=(\xi_1,\xi_2)=\Psi(h,E)\in\widetilde{\mathcal D}^{\,*}$, where $\xi_1=p_1/\ell_1$ and $\xi_2=p_2/\ell_2$, with $p_1,\ell_1,p_2,\ell_2$ the characteristic numbers of $\gamma$ defined in Section~\ref{sec:integration}. Setting $N=\operatorname{lcm}(\ell_1,\ell_2)$, we have that $NT$ is the least period of $\gamma$ and $\widehat{\mathbf M}_\gamma^{\,N}=\mathbf I$.

The above representation has the following geometric consequence.

\begin{proposition}\label{prop:no-closed-chamber}
A closed conformal extremal is not contained in any chamber of $\Ein$ (see Definition~\ref{def:chambers-walls}).
\end{proposition}

\begin{proof}
Suppose first that the conformal invariants are constant. By Theorem~\ref{thm:closed-constant-invariants}, in a suitable adapted pseudo-orthonormal basis the curve has the null lift
\[
\Gamma_{m,n}(\theta)
=\cos(m\theta)B_1+\sin(m\theta)B_2
+\cos(n\theta)B_3+\sin(n\theta)B_4,
\]
where $m,n$ are positive integers. Hence $\int_0^{2\pi}\Gamma_{m,n}(\theta)\,d\theta=0$. If the curve were contained in a chamber, there would be a nonzero vector $V$ such that $\langle\Gamma_{m,n}(\theta),V\rangle>0$ for every $\theta$, contradicting this identity.

Assume now that the conformal torsion is nonconstant.
With the normalization adopted in Proposition~\ref{prop:explicit-integration}, we have $\tau(0)=0$ and $\theta_\omega(0)=\theta_\nu(0)=0$, so that $\widehat\gamma(0)=(1,0,1,0,0)^t$. Using $\widehat{\mathbf M}_\gamma=R_1(-2\pi\xi_1,-2\pi\xi_2)$, we obtain
\[
\widehat\gamma(jT)=(\cos(2\pi j\xi_1),-\sin(2\pi j\xi_1);\cos(2\pi j\xi_2),-\sin(2\pi j\xi_2),0)
\]
for every integer $j$. As $0<\xi_1<\xi_2<1$ and $N\xi_1,N\xi_2\in\mathbb Z$, the sums of the corresponding $N$-th roots of unity vanish. Consequently, $\sum_{j=0}^{N-1}\widehat\gamma(jT)=0$.

Suppose, by contradiction, that $\gamma$ is contained in a chamber $\mathcal T$. Then there exists a nonzero vector $V\in\mathbb R^{2,3}$, expressed in the adapted coordinates, such that $\mathcal T=\{[X]\in\Ein:\langle X,V\rangle>0\}$, where $\langle\cdot,\cdot\rangle$ denotes the pseudo-Euclidean scalar product in the adapted coordinates. Hence $\langle\widehat\gamma(s),V\rangle>0$ for every $s$. In particular, $\langle\widehat\gamma(jT),V\rangle>0$, $j=0,\ldots,N-1$. Summing over $j$, we obtain $\sum_{j=0}^{N-1}\langle\widehat\gamma(jT),V\rangle>0$, a contradiction. This completes the proof.
\end{proof}

\begin{remark}
Since three-dimensional Minkowski, anti-de Sitter, and de Sitter spaces are conformally equivalent to chambers of $\Ein$, Proposition~\ref{prop:no-closed-chamber} implies that none of these spaces admits nondegenerate closed conformal extremals.
\end{remark}

\subsection{Symmetries of closed conformal extremals}

In the coordinates determined by $\mathbf B$, the pseudo-orthogonal group is represented by $O(2,3)$, with the conventions established in Subsection~\ref{subsec:einstein-universe}. For $\mathbf X\in O(2,3)$, we denote by $\mathbf X\cdot P$ the induced action on the $S^1\times S^2$-model of $\Ein$, obtained by applying $\mathbf X$ to the corresponding null vector and then normalizing.

\begin{definition}
The \emph{full symmetry group} of $\widehat\gamma$ is $G_\gamma=\{\mathbf X\in O(2,3):\mathbf X\cdot|\widehat\gamma|=|\widehat\gamma|\}$. Its elements are called \emph{symmetries} of $\widehat\gamma$. We also set $G_{\gamma,+}=G_\gamma\cap O_+(2,3)$ and $G_{\gamma,+}^{\uparrow}=G_\gamma\cap O_+^\uparrow(2,3)$. Thus $G_{\gamma,+}$ consists of the orientation-preserving symmetries, while $G_{\gamma,+}^{\uparrow}$ consists of those preserving both the orientation and the time-orientation of $\Ein$.
\end{definition}

\begin{remark}
Let $\gamma:\mathbb R/L\mathbb Z\to\Ein$ be a generic real-analytic spacelike immersion, parametrized by its natural parameter, with conformal torsion not identically zero, and suppose that $\gamma$ is not a multiple covering of its trajectory. Then every symmetry $\mathbf X\in G_\gamma$ induces a unique isometry $\varrho_{\mathbf X}$ of $\mathbb R/L\mathbb Z$ such that $\mathbf X\cdot\gamma([s])=\gamma(\varrho_{\mathbf X}([s]))$. This defines a faithful representation $\varrho:G_\gamma\to\operatorname{Isom}(\mathbb R/L\mathbb Z)$. Each $\varrho_{\mathbf X}$ is of the form $\varrho_{\mathbf X}([s])=[\varepsilon_{\mathbf X}s+c_{\mathbf X}]$, where $\varepsilon_{\mathbf X}\in\{\pm1\}$ is uniquely determined and $c_{\mathbf X}\in\mathbb R$ is uniquely determined modulo $L\mathbb Z$. For the closed conformal extremal under consideration, $L=NT$.
\end{remark}

Set $\mathbf D_1=\operatorname{diag}(1,1,1,1,-1)\in O_-^\uparrow(2,3)$ and $\mathbf D_2=\operatorname{diag}(1,-1,1,-1,-1)\in O_-^\downarrow(2,3)$, and define $\mathbf U=R_1(-\pi\xi_1,-\pi\xi_2)\mathbf D_1$ and $\mathbf S=\mathbf U\mathbf D_2$. We use the convention that $D_k$ denotes the dihedral group of order $2k$.

\begin{proposition}\label{prop:closed-symmetries}
Let $\widehat\gamma$ be as above. Then $G_{\gamma,+}^\uparrow=\langle\widehat{\mathbf M}_\gamma\rangle\simeq\mathbb Z_N$, $G_{\gamma,+}=\langle\widehat{\mathbf M}_\gamma,\mathbf S\rangle\simeq D_N$, and $G_\gamma=\langle\mathbf U,\mathbf S\rangle\simeq D_{2N}$.
\end{proposition}

\begin{proof}
The identities $\tau(s+T/2)=-\tau(s)$ and $\tau(-s)=-\tau(s)$, together with $\theta_\omega(s+T/2)=\theta_\omega(s)-\pi\xi_1$, $\theta_\nu(s+T/2)=\theta_\nu(s)-\pi\xi_2$, $\theta_\omega(T/2-s)=-\pi\xi_1-\theta_\omega(s)$, and $\theta_\nu(T/2-s)=-\pi\xi_2-\theta_\nu(s)$, imply $\widehat\gamma(s+T/2)=\mathbf U\cdot\widehat\gamma(s)$ and $\widehat\gamma(T/2-s)=\mathbf S\cdot\widehat\gamma(s)$.

Thus $\mathbf U,\mathbf S\in G_\gamma$, and their induced isometries of the parameter circle are $\varrho_{\mathbf U}([s])=[s+T/2]$ and $\varrho_{\mathbf S}([s])=[T/2-s]$. The first is a translation of order $2N$, while the second is a reflection. They generate a dihedral group of order $4N$. Since $\varrho$ is faithful, we conclude that $\langle\mathbf U,\mathbf S\rangle\simeq D_{2N}$. Moreover, $\mathbf U^2=\widehat{\mathbf M}_\gamma$.

We now show that $G_\gamma=\langle\mathbf U,\mathbf S\rangle$. Since a conformal symmetry can change the sign of the conformal torsion but preserves its absolute value, its zero set is invariant under $G_\gamma$. As $\tau(s)=r_-\operatorname{sn}(r_+s\mid m)$, this set, regarded as a subset of the parameter circle, is $Z_\tau=\{[jT/2]:j=0,\ldots,2N-1\}\subset\mathbb R/NT\mathbb Z$.

The subgroup of $\operatorname{Isom}(\mathbb R/NT\mathbb Z)$ preserving $Z_\tau$ is precisely the dihedral group generated by $\varrho_{\mathbf U}$ and $\varrho_{\mathbf S}$. Therefore $\varrho(G_\gamma)\subseteq\operatorname{Stab}(Z_\tau)=\langle\varrho_{\mathbf U},\varrho_{\mathbf S}\rangle$. The reverse inclusion follows from $\mathbf U,\mathbf S\in G_\gamma$. Hence, by faithfulness, $G_\gamma=\langle\mathbf U,\mathbf S\rangle\simeq D_{2N}$.

Finally, since $R_1(\phi_1,\phi_2)\in O_+^\uparrow(2,3)$, we have $\mathbf U\in O_-^\uparrow(2,3)$ and $\mathbf S\in O_+^\downarrow(2,3)$. Consequently, an element $\mathbf U^j\mathbf S^\delta$, with $0\le j<2N$ and $\delta\in\{0,1\}$, preserves the orientation precisely when $j$ is even, and among these it preserves the time-orientation precisely when $\delta=0$. It follows that $G_{\gamma,+}=\langle\mathbf U^2,\mathbf S\rangle=\langle\widehat{\mathbf M}_\gamma,\mathbf S\rangle\simeq D_N$ and $G_{\gamma,+}^\uparrow=\langle\mathbf U^2\rangle=\langle\widehat{\mathbf M}_\gamma\rangle\simeq\mathbb Z_N$. This completes the proof.
\end{proof}

\subsection{Geometry of the spherical projection}

We now consider the spherical projection $\beta(s)=(r(s)\cos\theta_\nu(s),r(s)\sin\theta_\nu(s),z(s))$ of $\widehat\gamma$.

Throughout this subsection, $(x,y,z)$ denote the Cartesian coordinates on $\mathbb R^3$, corresponding to the last three vectors $B_3,B_4,B_5$ of the adapted basis.

We begin by describing the oscillations of $\beta$ between two parallels and its motion around the polar axis. Since $r>0$ and $\theta_\nu'<0$, the spherical projection is a regular curve. Moreover, $\beta$ avoids the two poles $N_\gamma=(0,0,1)$ and $S_\gamma=(0,0,-1)$, and
\[
-z_*\le z(s)\le z_*,
\qquad
z_*=\frac{r_-\sqrt{\omega^2-\nu^2}}{\nu\sqrt{\omega^2-r_-^2}}\in(0,1).
\]
Thus $\beta$ is contained in the spherical zone bounded by the parallels $z=\pm z_*$. The maximum and minimum of the height function occur at $s\equiv T/4\pmod T$ and $s\equiv3T/4\pmod T$, respectively.

Since $z$ is a strictly increasing function of the conformal torsion, the height increases strictly from $-z_*$ to $z_*$ on $[-T/4,T/4]$ and decreases strictly from $z_*$ to $-z_*$ on $[T/4,3T/4]$. Consequently, each of the corresponding arcs of $\beta$ intersects every intermediate parallel exactly once.

Furthermore, since the longitude $\theta_\nu(s)$ is strictly decreasing, $\beta$ winds clockwise around the polar axis, as viewed from the north pole.

\begin{proposition}\label{prop:spherical-projection-symmetries}
The spherical projection $\beta$ has least period $\ell_2T$ and $2\ell_2$ ordinary geodesic inflections, counted as parameter values over one least period. Its full rotational symmetry group is the dihedral group $D_{\ell_2}$ of order $2\ell_2$, generated by the axial rotation $\mathbf R_{-2\pi\xi_2}$ and the half-turn $\mathcal R_\pi$ about the equatorial $x$-axis.
\end{proposition}

\begin{proof}
With respect to the outward-normal orientation of $S^2$, a direct but lengthy computation of the signed geodesic curvature gives
\[
k_\beta(s)
=
\frac{\sqrt{\omega^2-\nu^2}}{\omega}\,\tau(s)
\Bigl(
\frac{\omega^2-\tau(s)^2}
{1+(\omega^2-\nu^2)(\omega^2-\tau(s)^2)}
\Bigr)^{3/2}.
\]
Consequently, $k_\beta(s)=0$ if and only if $s=jT/2$, $j\in\mathbb Z$. All these zeros are simple, since $\tau'\ne0$ whenever $\tau=0$. Thus $\beta$ has exactly two ordinary geodesic inflections per period of the conformal torsion, occurring at its crossings of the equatorial plane.

Let $\mathbf R_\alpha$ denote the rotation of $S^2$ through the angle $\alpha$ about the polar axis, and let $\rho(x,y,z)=(x,y,-z)$ be reflection across the equatorial plane. The half-period identities $r(s+T/2)=r(s)$, $z(s+T/2)=-z(s)$, and $\theta_\nu(s+T/2)=\theta_\nu(s)-\pi\xi_2$ imply $\beta(s+T/2)=\mathbf R_{-\pi\xi_2}\rho(\beta(s))$ and $\beta(s+T)=\mathbf R_{-2\pi\xi_2}\beta(s)$.

Since $T$ is the least period of $z$ and $\xi_2=p_2/\ell_2$ is in lowest terms, the least period of $\beta$ is $\ell_2T$. In particular, $\beta$ has $2\ell_2$ ordinary geodesic inflections, counted as parameter values over one least period.

During this period, the height function attains its maximum precisely at the $\ell_2$ distinct points $\beta(T/4+jT)$, $j=0,\ldots,\ell_2-1$, which form a regular orbit under the cyclic group $\langle\mathbf R_{-2\pi\xi_2}\rangle\simeq\mathbb Z_{\ell_2}$. Since any axial rotation preserving the trajectory of $\beta$ must permute these points, this is exactly its group of axial rotational symmetries.

Moreover, the parity properties of $r,z,\theta_\nu$ imply $\beta(-s)=\mathcal R_\pi\beta(s)$, where $\mathcal R_\pi(x,y,z)=(x,-y,-z)$. Since $\mathcal R_\pi\mathbf R_\alpha\mathcal R_\pi=\mathbf R_{-\alpha}$, the transformations $\mathbf R_{-2\pi\xi_2}$ and $\mathcal R_\pi$ generate a dihedral group isomorphic to $D_{\ell_2}$.

Finally, every rotational symmetry of the trajectory preserves its set of geodesic inflections. Since $p_2$ and $\ell_2$ are relatively prime and $\ell_2\ge3$, these include at least three distinct points of the equator. Every rotational symmetry must therefore preserve the equatorial plane. Such a rotation is either axial or the product of $\mathcal R_\pi$ with an axial rotation. Hence $D_{\ell_2}$ is the full rotational symmetry group of the trajectory of $\beta$.
\end{proof}

\subsection{Geometric interpretation of the characteristic numbers}

Consider the decomposition $\widehat\gamma=(\eta,\beta):\mathbb R/NT\mathbb Z\to S^1\times S^2$, where $\eta(s)=(\cos\theta_\omega(s),\sin\theta_\omega(s))$. Here $NT$ is the least period of $\widehat\gamma$, although it need not be the least period of either component.

The spherical winding number is understood with respect to the rotational axis of the adapted spherical projection; the fixed product decomposition of the Einstein universe is not itself conformally invariant.

\begin{proposition}\label{prop:characteristic-winding}
Let $w_\eta$ be the winding number of $\eta$, and let $w_\beta$ be the winding number of $\beta$ about its rotational axis, both computed with respect to the adapted orientations. Then
\[
|w_\eta|=Np_1/\ell_1,
\qquad
|w_\beta|=Np_2/\ell_2.
\]
The integer $|w_\eta|$ determines the homotopy class of the closed extremal up to inversion in $\pi_1(\Ein)\simeq\mathbb Z$. Both absolute winding numbers are independent of the adapted basis and invariant under the action of $SO(2,3)$. Moreover, the characteristic numbers of $\gamma$ can be recovered from them and the order of $G_{\gamma,+}^{\uparrow}$.
\end{proposition}

\begin{proof}
Since the Einstein universe is diffeomorphic to $S^1\times S^2$, we have $\pi_1(\Ein)\simeq\mathbb Z$. Thus, with respect to the adapted identification, the free homotopy class of $\widehat\gamma$ is represented by $\deg(\eta)$. Changing the identification can reverse its sign. With the counterclockwise orientation of $S^1$ in the adapted coordinates, $w_\eta:=\deg(\eta)=(\theta_\omega(NT)-\theta_\omega(0))/(2\pi)=-N\xi_1=-Np_1/\ell_1$.

The spherical component $\beta$ avoids the poles $N_\gamma,S_\gamma$ of its rotational axis. The radial retraction onto the equator maps $\beta(s)$ to $(\cos\theta_\nu(s),\sin\theta_\nu(s))$. Orienting the equator counterclockwise in the adapted $(x_3,x_4)$-plane, we obtain $w_\beta=(\theta_\nu(NT)-\theta_\nu(0))/(2\pi)=-N\xi_2=-Np_2/\ell_2$.

Taking absolute values gives the asserted identities. Since $N$ and $\xi=\Psi(h,E)$ are invariant under the action of $SO(2,3)$, these formulas show that $|w_\eta|$ and $|w_\beta|$ are independent of the adapted basis and invariant under this action. Moreover, $\xi_1=|w_\eta|/N$ and $\xi_2=|w_\beta|/N$.

By Proposition~\ref{prop:closed-symmetries}, $N$ is the order of $G_{\gamma,+}^{\uparrow}$. The characteristic numbers $(p_1,\ell_1)$ and $(p_2,\ell_2)$ are therefore recovered by reducing these two fractions to lowest terms.
\end{proof}

\subsection{Toroidal model and examples}

We introduce a toroidal model for visualizing closed conformal extremals and illustrate it with three numerical examples.

Throughout this subsection, $\gamma$ denotes a closed conformal extremal with nonconstant conformal torsion and $\widehat\gamma$ its representation \eqref{eq:closed-gamma-hat} in the ray model determined by the adapted pseudo-orthonormal basis $\mathbf B$. The typewriter coordinates $(\mathtt{x},\mathtt{y},\mathtt{z})$ are Euclidean coordinates for the toroidal visualization, whereas $t$ denotes the height coordinate on the spherical factor; along $\widehat\gamma$, $t=z(s)$. Fix $0<\lambda<\sqrt2-1$ and set $\mathfrak d(\mathtt{x},\mathtt{y},\mathtt{z})=\bigl((\sqrt{\mathtt{x}^2+\mathtt{y}^2}-\sqrt2)^2+\mathtt{z}^2\bigr)^{1/2}$. Define
\[
\mathscr S_\lambda
=\bigl\{(\mathtt{x},\mathtt{y},\mathtt{z})\in\mathbb R^3:
1-\lambda<\mathfrak d(\mathtt{x},\mathtt{y},\mathtt{z})<1+\lambda\bigr\}.
\]
The boundary components of $\mathscr S_\lambda$ are tori with major radius $\sqrt2$ and minor radii $1-\lambda$ and $1+\lambda$, while $\mathfrak d=1$ is the middle torus. The restriction $\lambda<\sqrt2-1$ ensures that all these tori are embedded ring tori. For $u,v\in\mathbb R/2\pi\mathbb Z$ and $t\in(-1,1)$, define
\[
\Phi_\lambda(u,v,t)=
\bigl(
[\sqrt2+(1+\lambda t)\cos v]\cos u,\,
[\sqrt2+(1+\lambda t)\cos v]\sin u,\,
(1+\lambda t)\sin v
\bigr).
\]
Then $\mathfrak d(\Phi_\lambda(u,v,t))=1+\lambda t$. Hence $\Phi_\lambda$ is a diffeomorphism from $S^1\times S^1\times(-1,1)$ onto $\mathscr S_\lambda$.

Put $\Ein_*=\Ein\setminus(S^1\times\{N_\gamma,S_\gamma\})$. The adapted $S^1\times S^2$-model identifies $\Ein_*$ with $S^1\times S^1\times(-1,1)$, so that $\Phi_\lambda$ induces a diffeomorphism $\chi_\lambda:\Ein_*\to\mathscr S_\lambda$.

Under $\chi_\lambda$, the anti-de Sitter wall $\mathcal T_\infty(B_5)$ corresponds to the middle torus $\mathfrak d=1$, while $\mathcal T_+(B_5)$ and $\mathcal T_-(B_5)$, with their polar circles removed, correspond to $\mathfrak d>1$ and $\mathfrak d<1$, respectively. Since $\widehat\gamma$ avoids the polar circles, its toroidal portrait $\mathcal P_{\gamma,\lambda}=\chi_\lambda\circ\widehat\gamma$ is parametrized by
\[
\mathcal P_{\gamma,\lambda}(s)=\Phi_\lambda\bigl(\theta_\omega(s),\theta_\nu(s),z(s)\bigr),
\qquad s\in\mathbb R/NT\mathbb Z.
\]
The parameter $\lambda$ controls only the Euclidean thickness of the shell and has no conformal-geometric significance. We use the shell only as a \emph{topological visualization}: its Euclidean metric does not represent the conformal structure of the Einstein universe, and $\chi_\lambda$ does not, in general, turn conformal symmetries into Euclidean isometries.

Given characteristic numbers $(p_1,\ell_1,p_2,\ell_2)$, with $\xi=(p_1/\ell_1,p_2/\ell_2)\in\widetilde{\mathcal D}^{\,*}$, inversion of the period map gives $(m,q)$ and hence $(h,E)$, the fundamental period $T$, the height $z(s)$, and the angles $\theta_\omega(s),\theta_\nu(s)$. The integration formula reconstructs the extremal up to conformal congruence, and the formula above gives its toroidal portrait for any admissible value of $\lambda$, with $N=\operatorname{lcm}(\ell_1,\ell_2)$.

For the three examples below, the period map was inverted using Python 3.12.14, NumPy 2.3.5, and SciPy 1.17.0. The complete elliptic integral of the third kind was evaluated through Carlson's symmetric integrals~\cite{Carlson1995}, and the nonlinear system was solved with \texttt{scipy.optimize.root}. The Python code was written by ChatGPT (OpenAI). For the toroidal portraits displayed below, we take $\lambda=\frac34(\sqrt2-1)$.

\Needspace{7\baselineskip}\medskip\noindent\textbf{Example 1.} For $\xi=(1/2,3/4)$, the characteristic numbers are $(1,2,3,4)$ and $N=4$. The symmetry groups are $G_{\gamma,+}^{\uparrow}\simeq\mathbb Z_4$, $G_{\gamma,+}\simeq D_4$, and $G_\gamma\simeq D_8$. Inversion gives \mbox{$(m,q)\simeq(0.4123944,\!0.9283029)$} and~\mbox{$(h,E)\simeq(2.5892881,\!5.5439786)$}. Here $T\simeq3.7313529$ and $NT=4T\simeq14.9254118$. Over the least period, the circular projection makes two clockwise turns and has least period $2T$, while the spherical projection winds three times around its polar axis, has least period $4T$, and rotational symmetry group $D_4$. The toroidal portrait has no self-intersections; since $z_*\simeq0.6119569$, it lies between $\mathfrak d\simeq0.8098894$ and $\mathfrak d\simeq1.1901106$ and crosses $\mathfrak d=1$ at the eight distinct points corresponding to $s=jT/2$, $j=0,\ldots,7$.

\begin{figure}[H]
\centering
\makebox[\linewidth][c]{%
\includegraphics[width=0.437\linewidth]{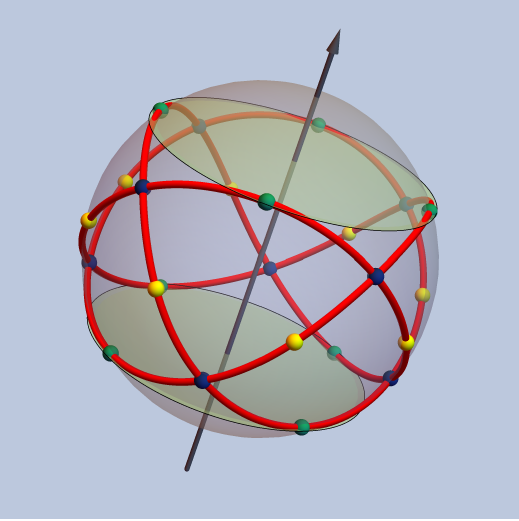}%
\hspace{0.078\linewidth}%
\includegraphics[width=0.437\linewidth]{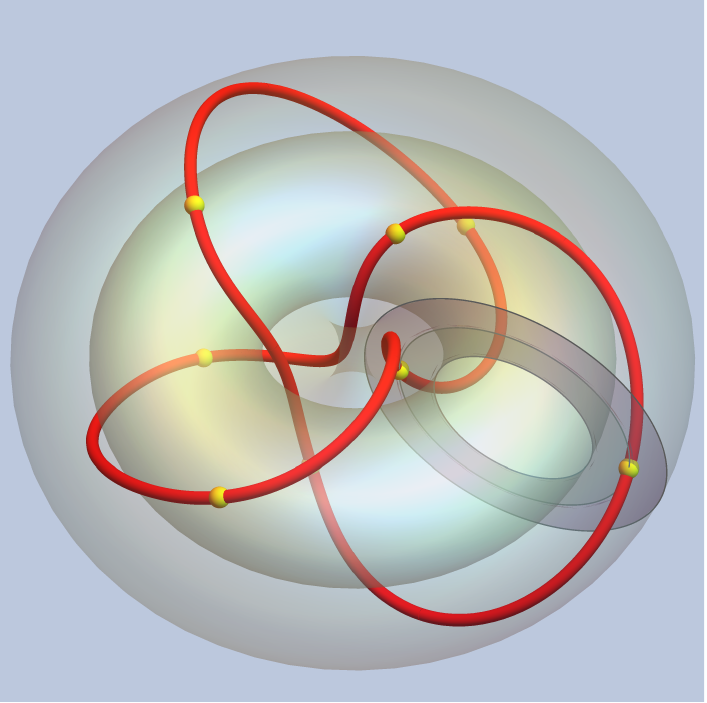}%
}
\caption{Spherical projection (left) and toroidal portrait (right) for $\xi=(1/2,3/4)$. On the spherical projection, yellow marks equatorial inflections, blue double points, and green height extrema; on the toroidal portrait, yellow marks the eight wall crossings.}
\label{fig:example-one}
\end{figure}

\medskip\noindent\textbf{Example 2.} For $\xi=(2/5,4/5)$, the characteristic numbers are $(2,5,4,5)$ and $N=5$. The symmetry groups are $G_{\gamma,+}^{\uparrow}\simeq\mathbb Z_5$, $G_{\gamma,+}\simeq D_5$, and $G_\gamma\simeq D_{10}$. Inversion gives \mbox{$(m,q)\simeq(0.4801361,\!0.8163340)$} and~\mbox{$(h,E)\simeq(3.4877398,\!10.6637316)$}. Here $T\simeq3.3859083$ and $NT=5T\simeq16.9295414$. Over the least period, the circular projection makes two clockwise turns and has least period $5T/2$, while the spherical projection winds four times around its polar axis, has least period $5T$, and rotational symmetry group $D_5$. The toroidal portrait has five ordinary double points on $\mathfrak d=1$, arising from the ten values $s=jT/2$, $j=0,\ldots,9$, paired as $(j,j+5)$; since $z_*\simeq0.7719780$, it lies between $\mathfrak d\simeq0.7601772$ and $\mathfrak d\simeq1.2398228$.

\begin{figure}[H]
\centering
\makebox[\linewidth][c]{%
\includegraphics[width=0.399\linewidth]{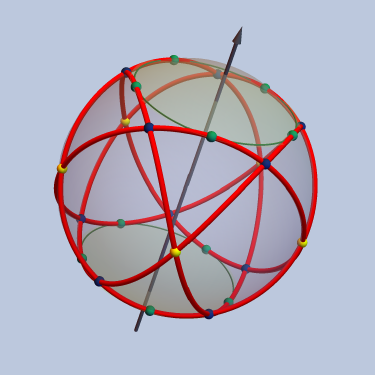}%
\hspace{0.078\linewidth}%
\includegraphics[width=0.399\linewidth]{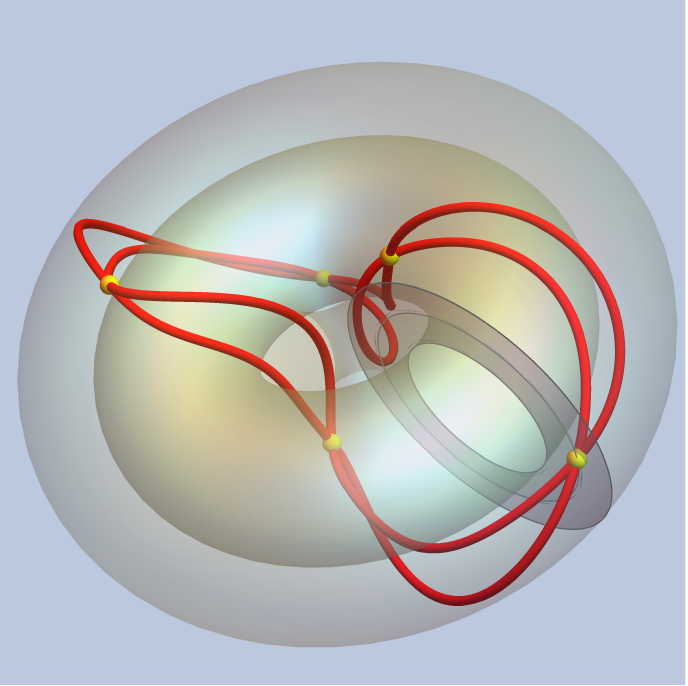}%
}
\caption{Spherical projection (left) and toroidal portrait (right) for $\xi=(2/5,4/5)$. On the spherical projection, yellow marks equatorial inflections, blue non-equatorial double points, and green height extrema; on the toroidal portrait, yellow marks the five double points on the wall.}
\label{fig:example-two}
\end{figure}

\Needspace{7\baselineskip}\medskip\noindent\textbf{Example 3.} For $\xi=(3/7,5/7)$, the characteristic numbers are $(3,7,5,7)$ and $N=7$. The symmetry groups are $G_{\gamma,+}^{\uparrow}\simeq\mathbb Z_7$, $G_{\gamma,+}\simeq D_7$, and $G_\gamma\simeq D_{14}$. Inversion gives \mbox{$(m,q)\simeq(0.7024817,\!0.9024808)$} and~\mbox{$(h,E)\simeq(6.3406057,\!38.9754902)$}. Here $T\simeq3.0470750$ and $NT=7T\simeq21.3295248$. Over the least period, the circular projection makes three clockwise turns and the spherical projection winds five times around its polar axis; both have least period $7T$, and the rotational symmetry group of the spherical projection is $D_7$. Since $z_*\simeq0.7330323$, the toroidal portrait lies between $\mathfrak d\simeq0.7722761$ and $\mathfrak d\simeq1.2277239$ and crosses $\mathfrak d=1$ at the fourteen distinct points corresponding to $s=jT/2$, $j=0,\ldots,13$.

\begin{figure}[H]
\centering
\makebox[\linewidth][c]{%
\includegraphics[width=0.399\linewidth]{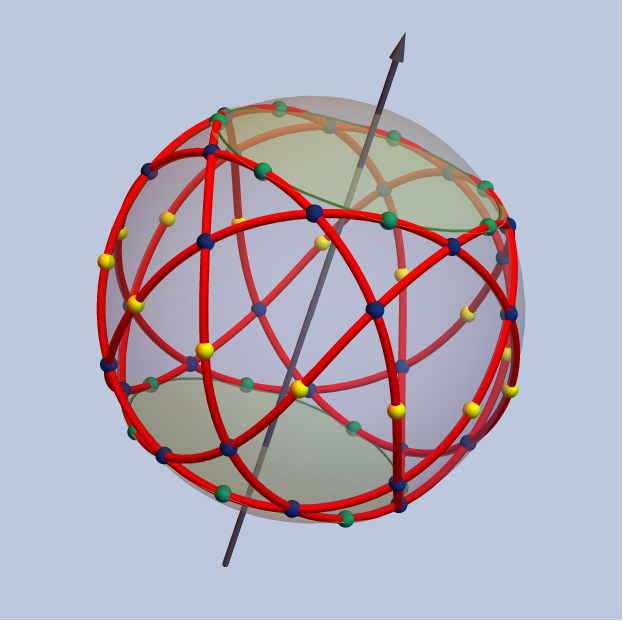}%
\hspace{0.078\linewidth}%
\includegraphics[width=0.399\linewidth]{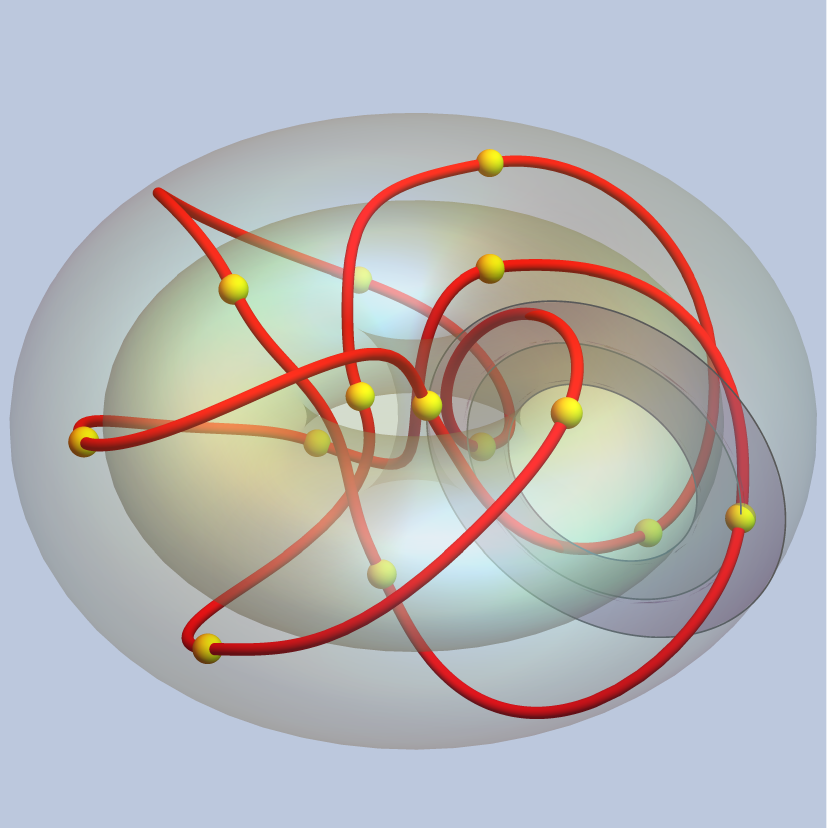}%
}
\caption{Spherical projection (left) and toroidal portrait (right) for $\xi=(3/7,5/7)$. On the spherical projection, yellow marks equatorial inflections, blue double points, and green height extrema; on the toroidal portrait, yellow marks the fourteen wall crossings.}
\label{fig:example-three}
\end{figure}

\par\vspace{-3pt}
\Needspace{10\baselineskip}
\section*{Acknowledgments}
\vspace{-2pt}
\begingroup\small
The author acknowledges the use of ChatGPT (OpenAI) during the preparation of this paper for mathematical discussion, critical review of arguments, assistance with symbolic and numerical computations, and the preparation and revision of the LaTeX source files of the manuscript. ChatGPT was also used in the preparation of some of the graphical material and of the numerical code described in Section~4.4. The author takes full responsibility for the mathematical statements, proofs, computations, and conclusions of the paper.
\endgroup

{\fontsize{8.7}{10}\selectfont

}
\vspace{.5em}
{\small
\noindent\textsc{(E. Musso)} Dipartimento di Scienze Matematiche L.~Lagrange,
Politecnico di Torino, Corso Duca degli Abruzzi 24,
I-10129 Torino, Italy\\
E-mail address: \href{mailto:emilio.musso@polito.it}{\texttt{emilio.musso@polito.it}}\\
ORCID iD: \href{https://orcid.org/0000-0002-0899-819X}{0000-0002-0899-819X}
\par}
\end{document}